\title{Spectral measures associated to rank two Lie groups \\ and finite subgroups of $GL(2,\mathbb{Z})$}
\author{
{\sc David E.\ Evans and Mathew Pugh}\\
 {\footnotesize School of Mathematics, Cardiff University,}\\  {\footnotesize Senghennydd Road, Cardiff CF24 4AG, Wales, U.K.}
}
\date{\today}

\documentclass[12pt]{article}

\usepackage{amsmath,amssymb,amsthm}
\usepackage{graphicx}
\usepackage{framed}
\usepackage{color}
\usepackage[usenames,dvipsnames]{xcolor}
\usepackage{multirow}
\usepackage{cancel}
\usepackage{array}
\usepackage{arydshln}

\theoremstyle{definition}
\newtheorem{Def}{Definition}[section]

\newtheorem{Lemma}[Def]{Lemma}

\newtheorem{Thm}[Def]{Theorem}
\newtheorem{Rem}[Def]{Remark}
\newtheorem{Conj}[Def]{Conjecture}

\begin{document}
\maketitle

\begin{abstract}
Spectral measures for fundamental representations of the rank two Lie groups $SU(3)$, $Sp(2)$ and $G_2$ have been studied. Since these groups have rank two, these spectral measures can be defined as measures over their maximal torus $\mathbb{T}^2$ and are invariant under an action of the corresponding Weyl group, which is a subgroup of $GL(2,\mathbb{Z})$.
Here we consider spectral measures invariant under an action of the other finite subgroups of $GL(2,\mathbb{Z})$. These spectral measures are all associated with fundamental representations of other rank two Lie groups, namely $\mathbb{T}^2=U(1) \times U(1)$, $U(1) \times SU(2)$, $U(2)$, $SU(2) \times SU(2)$, $SO(4)$ and $PSU(3)$.
\end{abstract}

{\footnotesize
\tableofcontents
}

\section{Introduction} \label{sect:intro}

In \cite{evans/pugh:2009v, evans/pugh:2010i, evans/pugh:2012i, evans/pugh:2012ii, evans/pugh:2012iii} the authors have studied spectral measures for fundamental representations of the rank two Lie groups $SU(3)$, $Sp(2)$ and $G_2$.
The spectral measure of an operator is a particular compactly supported probability measure on the spectrum of that operator. A representation graph of $G$ is the fusion graph for an irreducible character of $G$. By a spectral measure for the Lie group $G$ we will mean the spectral measure of (the adjacency matrix of) such a fusion graph.
In particular, the joint spectral measure for the two fundamental representations of the rank two Lie group were studied. The joint spectral measure has support given by the joint spectrum $\sigma(\Delta_1,\Delta_2)$ of the operators given by the adjacency matrices $\Delta_1, \Delta_2$ of the representation graphs for the two fundamental representations $\rho_1,\rho_2$ of $G$. The pushforward of these joint spectral measures under the projection on the spectrum $\sigma(\Delta_j)$ of a single fundamental representation $\rho_j$ then yields the spectral measure for $\rho_j$ on $\sigma(\Delta_j)$.

Since these groups have rank two, these spectral measures can also be defined as measures over the two-torus $\mathbb{T}^2$ (which is isomorphic to the maximal torus of the rank two Lie group).
The characters of the fundamental representations define a surjection from $\mathbb{T}^2$ to the joint spectrum $\sigma(\Delta_1,\Delta_2)$, which is invariant under the action of the Weyl group $W(G)$ of $G$ on the torus. The Jacobian for the change of variables from $\mathbb{T}^2$ to  $\sigma(\Delta_1,\Delta_2)$ clearly plays a key role in the spectral measures for $G$ over both $\mathbb{T}^2$ and the joint spectrum.
This approach fits with the spectral measure blowup philosophy of \cite{banica/bichon:2014}.

In conformal field theories built from these Lie groups, one considers the Verlinde algebra at a finite level $k$, represented by a non-degenerately braided system ${}_N \mathcal{X}_N$ of irreducible endomorphisms on a type $\mathrm{III}_1$ factor $N$, whose fusion rules $\{ \mathcal{N}_{\lambda \nu}^{\mu} \}$ reproduce exactly those of the positive energy representations of the loop group of the Lie group $G$ at level $k$,
$\mathcal{N}_{\lambda} \mathcal{N}_{\mu} = \sum_{\nu} \mathcal{N}_{\lambda \nu}^{\mu} \mathcal{N}_{\nu}$.
The statistics generators $S$, $T$ obtained from the braided tensor category ${}_N \mathcal{X}_N$ match exactly those of the Ka\u{c}-Peterson modular $S$, $T$ matrices which perform the conformal character transformations \cite{wassermann:1998} (see also footnote 2 in \cite{bockenhauer/evans:2001}).
The fusion graph for these irreducible endomorphisms are truncated versions of the representation graphs of $G$ itself.

A braided subfactor is an inclusion $N \subset M$ where the dual canonical endomorphism decomposes as a finite combination of elements in ${}_N \mathcal{X}_N$, and yields a modular invariant partition function through the procedure of $\alpha$-induction \cite{bockenhauer/evans/kawahigashi:1999, bockenhauer/evans:2000, evans:2003}.
The action of the $N$-$N$ sectors ${}_N \mathcal{X}_N$ on the $M$-$N$ sectors ${}_M \mathcal{X}_N$ and produces a nimrep (non-negative integer matrix representation of the original fusion rules)
$\mathcal{G}_{\lambda} \mathcal{G}_{\mu} = \sum_{\nu} \mathcal{N}_{\lambda \nu}^{\mu} \mathcal{G}_{\nu}$
whose spectrum reproduces exactly the diagonal part of the modular invariant. In the case of the trivial embedding of $N$ in itself, the nimrep $\mathcal{G}$ is simply $\mathcal{N}$. The joint spectrum of the nimrep graphs for the fundamental generators of the system ${}_N \mathcal{X}_N$ is again contained in the joint spectrum $\sigma(\Delta_1,\Delta_2)$ of $G$.
One can then determine the (joint) spectral measure for these nimrep graphs over both $\mathbb{T}^2$ and $\sigma(\Delta_1,\Delta_2)$ -- see \cite{evans/pugh:2009v, evans/pugh:2010i, evans/pugh:2012i, evans/pugh:2012iii} for more details in the cases of the rank two Lie groups $SU(3)$, $Sp(2)$ and $G_2$.

For a rank two Lie group $G$, the spectrum of the McKay graphs (or representation graphs) of a finite subgroup $H \subset G$ are also contained in the joint spectrum $\sigma(\Delta_1,\Delta_2)$ of $G$.
One can thus also determine the (joint) spectral measure for these graphs over both $\mathbb{T}^2$ and $\sigma(\Delta_1,\Delta_2)$ \cite{evans/pugh:2012ii}.

The compact semisimple rank two Lie algebras are $SU(2) \times SU(2)$, $SU(3)$, $Sp(2)$ and $G_2$ (they are in fact simple, apart from $SU(2) \times SU(2)$), and are all connected. Their Weyl groups are the dihedral groups $D_4$, $D_6 \equiv S_3$, $D_8$ and $D_{12}$ respectively, where $D_n$ is the dihedral group of order $n$, which are all subgroups of $GL(2,\mathbb{Z})$.
Spectral measures for $SU(3)$, $Sp(2)$ and $G_2$ were considered in \cite{evans/pugh:2009v, evans/pugh:2010i, evans/pugh:2012i, evans/pugh:2012ii, evans/pugh:2012iii}.
Spectral measures for $SU(2) \times SU(2)$ are products of the spectral measures for $SU(2)$, considered in \cite{banica/bisch:2007, evans/pugh:2009v}.

In this paper we consider spectral measures associated to all the finite subgroups $\Gamma$ of $GL(2,\mathbb{Z})$ which are not also finite subgroups of $SL(2,\mathbb{Z})$. The spectral measures considered are invariant under an action of the group.
We consider two types of representation graphs $\mathcal{G}_{\rho}$ and $\mathcal{H}_{\rho}$, and hence determine spectral measures for two types of operators. The representation graph $\mathcal{G}_{\rho}$ (with adjacency matrix denoted by $\Delta_{\rho}$) is given by the fusion rules for characters with respect to multiplication by the character $\chi_{\rho}$ for the representation $\rho$.
The representation graph $\mathcal{H}_{\rho}$ (with adjacency matrix denoted by ${}^{\Gamma} \hspace{-1mm} \Delta_{\rho}$) is given by the fusion rules for characters of $\mathbb{T}^2$, with respect to multiplication by the restriction of $\chi_{\rho}$ to $\mathbb{T}^2$.
(Joint) spectral measures for the representation graphs associated with fundamental representations of other rank two Lie groups, namely $\mathbb{T}^2=U(1) \times U(1)$, $U(1) \times SU(2)$, $U(2)$, $SU(2) \times SU(2)$, $SO(4)$ and $PSU(3)$ are studied.

The paper is organised as follows. In Section \ref{sect:preliminaries} we discuss some preliminary material, beginning with the finite subgroups of $GL(2,\mathbb{Z})$. Then in Section \ref{sect:orbit_functions} we discuss orbit functions \cite{klimyk/patera:2006, klimyk/patera:2007} which have been the objects of much attention in the last decade, see e.g. \cite{nesterenko/patera/tereszkiewicz:2011}.
For an irreducible representation of a compact semisimple Lie group $G$ with Weyl group $W(G)$, these orbit functions are the contribution to the character from a single orbit under the action of $W(G)$. We present an analogous definition in the case of an orbit under the action on $\mathbb{T}^2$ of any reflection group which is a finite subgroup $\Gamma$ of $GL(2,\mathbb{Z})$. The orbit functions are used to define formal characters (which coincide with the characters of $G$ for $\Gamma$ the Weyl group $W(G)$), which in turn are used in Section \ref{sect:representation_graphs} to define families of representation graphs associated with the finite subgroup $\Gamma$. The first family is given by the fusion graphs $\mathcal{G}$ for these formal characters, and the second by the fusion graphs $\mathcal{H}$ for the action of these formal characters on $\mathbb{T}^2$.
A discussion on the relation between spectral measures over certain different domains is given in Section \ref{sect:measures-different_domains}.

In Section \ref{sect:measures_over_T2} we discuss the spectral measure for the representation graphs $\mathcal{H}$, $\mathcal{G}$ over $\mathbb{T}^2$ (for all finite subgroups of $GL(2,\mathbb{Z})$). Then in Sections \ref{sect:Z0}-\ref{sect:D6(1)} we determine the (joint) spectral measures for the representation graphs over the (joint) spectrum for each finite subgroup of $GL(2,\mathbb{Z})$. This will essentially be obtained by determining the Jacobian $J = J_{\Gamma}$ of a particular change of variable for each group $\Gamma$.
As one consequence, we prove a conjecture from the Online Encyclopedia of Integer Sequences (OEIS) \cite{OEIS:2010}, namely that the number of walks on $\mathbb{N}^2$ starting and ending at $(0,0)$ and consisting of $2n$ steps taken from $\{(-1, -1), (-1, 1), (1, -1), (1, 1)\}$, is given by the squared Catalan numbers A001246 (see Remark \ref{Rem:OEIS-conjecture}).

\section{Preliminaries} \label{sect:preliminaries}
\subsection{Finite subgroups of $GL(2,\mathbb{Z})$} \label{sect:subgroups-GL(2,Z)}

There are 13 finite subgroups $\Gamma \subset GL(2,\mathbb{Z})$, up to conjugacy in $GL(2,\mathbb{Z})$ \cite{newman:1972}:
\begin{eqnarray*}
(i) & \Gamma \subset SL(2,\mathbb{Z}): & \;\; \mathbb{Z}_0, \; \mathbb{Z}_2^{(1)}, \; \mathbb{Z}_3, \; \mathbb{Z}_4, \; \mathbb{Z}_6, \\
(ii) & \Gamma \not\subset SL(2,\mathbb{Z}): & \;\; \mathbb{Z}_2^{(2)}, \; \mathbb{Z}_2^{(3)}, \; D_4^{(1)}, \; D_4^{(2)}, \; D_6^{(1)}, \; D_6^{(2)}, \; D_8, \; D_{12},
\end{eqnarray*}
where $\Gamma^{(i)}$, $\Gamma^{(j)}$ denote non-conjugate embeddings of $\Gamma$ in $GL(2,\mathbb{Z})$ for $i \neq j$.
The first five subgroups are also finite subgroups of $SL(2,\mathbb{Z})$. We denote by $\mathfrak{G}$ the set of all finite subgroups of $GL(2,\mathbb{Z})$ (up to conjugacy in $GL(2,\mathbb{Z})$) which are not finite subgroups of $SL(2,\mathbb{Z})$, i.e. those listed in $(ii)$.
Generators for these groups are given below:
\begin{eqnarray*}
& \mathbb{Z}_0: \;\; I, \qquad \quad \mathbb{Z}_2^{(1)}: \;\; -I, \qquad \quad \mathbb{Z}_3: \;\; T_3 = \left( \begin{array}{cc} 0 & -1 \\ 1 & -1 \end{array} \right), \qquad \quad \mathbb{Z}_4: \;\; T_4 = \left( \begin{array}{cc} 0 & -1 \\ 1 & 0 \end{array} \right), \\
& \mathbb{Z}_6: \;\; T_6 = \left( \begin{array}{cc} 0 & 1 \\ -1 & 1 \end{array} \right), \qquad \quad \mathbb{Z}_2^{(2)}: \;\; T_2 = \left( \begin{array}{cc} 1 & 0 \\ 0 & -1 \end{array} \right), \qquad \quad \mathbb{Z}_2^{(3)}: \;\; T_2' = \left( \begin{array}{cc} 0 & 1 \\ 1 & 0 \end{array} \right), \\
& D_4^{(1)}: \;\; -I, T_2, \qquad \quad D_4^{(2)}: \;\; -I, T_2', \qquad \quad D_6^{(1)}: \;\; T_3, T_2', \\
& D_6^{(2)}: \;\; T_3, -T_2', \qquad \quad D_8: \;\; T_4, T_2', \qquad \quad D_{12}: \;\; T_6, T_2'.
\end{eqnarray*}

There is an obvious action of $\Gamma \subset GL(2,\mathbb{Z})$ on $\mathbb{R}^2$, that is, $T(m,n) = (a_{11}m+a_{12}n,a_{21}m+a_{22}n)$, for $m,n \in \mathbb{R}$, which drops to the quotient $\mathbb{R}^2/\mathbb{Z}^2 \cong \mathbb{T}^2$.

For a finite subgroup $\Gamma \in \mathfrak{G}$, we denote by $P_+$ the fundamental domain of the quotient $\mathbb{Z}^2/\Gamma$ such that $(\lambda_1,\lambda_2) \in P_+$ for all $0 \leq \lambda_2 \leq \lambda_1$, and by $P_{++}$ the set $P_{++} = \{ \lambda \in P_+ | \, \gamma\lambda \neq \lambda \textrm{ for any } \gamma \in \Gamma \}$, i.e. if $\lambda \in P_+ \setminus P_{++}$ then $\lambda$ lies on the boundary of $P_+$.

\begin{figure}[tb]
\begin{center}
 \includegraphics[width=80mm]{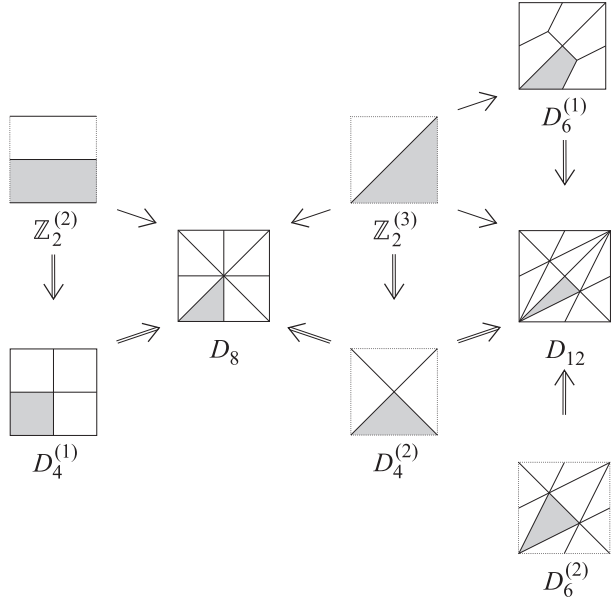} \label{Fig-subgroup_embeddings}
 \caption{Inclusions of finite subgroups of $GL(2,\mathbb{Z})$}
\end{center}
\end{figure}

We have inclusions of these subgroups as illustrated in Figure \ref{Fig-subgroup_embeddings}. The lines between these subgroups indicate the inclusions of one subgroup in another, where a double line denotes that one is a normal subgroup of the other. Above each subgroup $\Gamma$ is a diagram illustrating the lines of reflection (the solid lines in each diagram) given by the action of $\Gamma$ on $\mathbb{T}^2$, and the shaded region indicates a fundamental domain of $\mathbb{T}^2/\Gamma$. The individual diagrams are given in more detail in Sections \ref{sect:Z2^2}-\ref{sect:D6(1)}.

We are interested in the finite subgroups of $GL(2,\mathbb{Z})$ which could appear as the Weyl group for a rank two Lie group.
The finite subgroups of $SL(2,\mathbb{Z})$ are generated by rotations of $\mathbb{Z}^2$, whilst the subgroups in $\mathfrak{G}$ are generated by reflections of $\mathbb{Z}^2$.
As Weyl groups are finite reflection groups, we will restrict our attention to the groups in $\mathfrak{G}$, i.e. the finite subgroups of $GL(2,\mathbb{Z})$ that are not subgroups of $SL(2,\mathbb{Z})$.
We note however that $\mathbb{Z}_n \in SL(2,\mathbb{Z})$ is a normal subgroup of a dihedral group $D_{2n} \cong \mathbb{Z}_n \rtimes \mathbb{Z}_2 \in GL(2,\mathbb{Z})$. More precisely, $\mathbb{Z}_2^{(1)} \rtimes \mathbb{Z}_2^{(j)} \cong D_4^{(j-1)}$ for $j=2,3$, $\mathbb{Z}_3 \rtimes \mathbb{Z}_2^{(3)} \cong D_6^{(1)}$, $\mathbb{Z}_4 \rtimes \mathbb{Z}_2^{(j)} \cong D_8$ for $j=2,3$, and $\mathbb{Z}_6 \rtimes \mathbb{Z}_2^{(3)} \cong D_{12}$.

\subsection{Rank two Lie groups} \label{sect:rank2LieGroups}

The simply-connected compact simple rank two Lie groups are $SU(3)$, $Sp(2)$ (the second order symplectic group, the set of $4\times4$ unitary matrices $U$ such that $U^T JU = J$, where $J = I_2 \otimes T_4$, for $I_2$ the $2\times2$ identity matrix), and $G_2$, whilst the only other simply-connected compact semisimple rank two Lie group is $SU(2) \times SU(2)$. Their Weyl groups are the dihedral groups $D_6 \equiv S_3$, $D_8$, $D_{12}$ and $D_4$ respectively, where $D_n$ is the dihedral group of order $n$, which are all subgroups of $GL(2,\mathbb{Z})$.
All compact connected rank two Lie groups are given by quotients of $SU(3)$, $Sp(2)$, $G_2$ or of products of $SU(2)$, $\mathbb{T} = U(1)$, by finite subgroups of their centers. The centers of $SU(3)$, $Sp(2)$, $G_2$, $SU(2)$ are $\mathbb{Z}_3$, $\mathbb{Z}_2$, $\mathbb{Z}_0$, $\mathbb{Z}_2$ respectively. Thus all connected but non-simply-connected compact rank two Lie groups are given by $\mathbb{T}^2 = U(1) \times U(1)$, $\mathbb{T} \times SU(2) = U(1) \times SU(2)$,
the double covered groups $(\mathbb{T} \times SU(2))/\mathbb{Z}_2 = U(2)$, $(SU(2) \times SU(2))/\mathbb{Z}_2 = SO(4)$, $\mathbb{T} \times (SU(2)/\mathbb{Z}_2) = U(1) \times SO(3)$, $SU(2) \times (SU(2)/\mathbb{Z}_2) = SU(2) \times SO(3)$, $Sp(2)/\mathbb{Z}_2 = SO(5)$,
the triple covered group $SU(3)/\mathbb{Z}_3 = PSU(3)$,
and the quadruple covered group $(SU(2)/\mathbb{Z}_2) \times (SU(2)/\mathbb{Z}_2) = SO(3) \times SO(3)$.
These Lie groups have Weyl groups 0, $\mathbb{Z}_2$,
$\mathbb{Z}_2$, $D_4$, $\mathbb{Z}_2$, $D_4$, $D_8$,
$D_6$
and $D_4$ respectively, which are also all finite subgroups of $GL(2,\mathbb{Z})$. There are two non-conjugate embeddings each of $\mathbb{Z}_2$, $D_4$ and $D_6$ in $GL(2,\mathbb{Z})$, as described in Section \ref{sect:subgroups-GL(2,Z)}. By comparing the action of the Weyl group on the maximal torus with the action of $\Gamma \subset GL(2,\mathbb{Z})$ on $\mathbb{T}^2$ described in Section \ref{sect:subgroups-GL(2,Z)}, we find that both embeddings in $GL(2,\mathbb{Z})$ of each of $\mathbb{Z}_2$, $D_4$ and $D_6$ appear, as described in Table \ref{Table:Lie-Weyl}. Thus the Weyl groups for the connected compact rank two Lie groups listed above in fact exhaust all finite subgroups of $GL(2,\mathbb{Z})$.

\begin{table}[t]
\begin{center}
\begin{tabular}{|c|c|c|c|c|} \hline
Compact Lie group $G$ & Lie algebra & $W_G$ & $\pi_0(G)$ & $\pi_1(G)$ \\
\hline 
\hline
$SU(3)$ & $A_2 = \mathfrak{su}(3)$ & $D_6^{(2)}$ & 0 & 0 \\
$Sp(2)$ & $C_2 = \mathfrak{sp}(2)$ & $D_8$ & 0 & 0 \\
$G_2$ & $\mathfrak{g}_2$ & $D_{12}$ & 0 & 0 \\
$SU(2) \times SU(2)$ & $A_1 \times A_1 = \mathfrak{su}(2) \times \mathfrak{su}(2)$ & $D_4^{(1)}$ & 0 & 0 \\
\hline 
$\mathbb{T}^2=U(1) \times U(1)$ & $\mathfrak{t}^2$ & 0 & 0 & $\mathbb{Z}^2$ \\
$\mathbb{T} \times SU(2)$ & $\mathfrak{t} \times \mathfrak{su}(2)$ & $\mathbb{Z}_2^{(2)}$ & 0 & $\mathbb{Z}$ \\
\hline 
$U(2)$ & $\mathfrak{u}(2) \cong \mathfrak{t} \times \mathfrak{su}(2)$ & $\mathbb{Z}_2^{(3)}$ & 0 & $\mathbb{Z}$ \\
$SO(4)$ & $D_2 = \mathfrak{so}(4) \cong \mathfrak{su}(2) \times \mathfrak{su}(2)$ & $D_4^{(2)}$ & 0 & $\mathbb{Z}_2$ \\
$PSU(3)$ & $\mathfrak{su}(3)$ & $D_6^{(1)}$ & 0 & $\mathbb{Z}_3$ \\
\hdashline
$\mathbb{T} \times SO(3)$ & $\mathfrak{t} \times \mathfrak{su}(2)$ & $\mathbb{Z}_2^{(2)}$ & 0 & $\mathbb{Z} \times \mathbb{Z}_2$ \\
$SU(2) \times SO(3)$ & $\mathfrak{su}(2) \times \mathfrak{su}(2)$ & $D_4^{(1)}$ & 0 & $\mathbb{Z}_2$ \\
$SO(3) \times SO(3)$ & $\mathfrak{su}(2) \times \mathfrak{su}(2)$ & $D_4^{(1)}$ & 0 & $\mathbb{Z}_2 \times \mathbb{Z}_2$ \\
$SO(5)$ & $B_2 = \mathfrak{so}(5) \cong \mathfrak{sp}(2)$ & $D_8$ & 0 & $\mathbb{Z}_2$ \\
\hline 
$\mathbb{T} \times O(2) = \mathbb{T} \times (\mathbb{T} \rtimes \mathbb{Z}_2)$ & $\mathfrak{t}^2$ & $\mathbb{Z}_2^{(2)}$ & $\mathbb{Z}_2$ & - \\
$O(2) \times O(2)$ & $\mathfrak{t}^2$ & $D_4^{(1)}$ & $\mathbb{Z}_2 \times \mathbb{Z}_2$ & - \\
$O(2) \times SU(2)$ & $\mathfrak{t} \times \mathfrak{su}(2)$ & $D_4^{(1)}$ & $\mathbb{Z}_2$ & - \\
$O(2) \times SO(3)$ & $\mathfrak{t} \times \mathfrak{su}(2)$ & $D_4^{(1)}$ & $\mathbb{Z}_2$ & - \\
$\mathbb{T}^2 \rtimes \mathbb{Z}_2^{(3)}$ & $\mathfrak{t}^2$ & $\mathbb{Z}_2^{(3)}$ & $\mathbb{Z}_2$ & - \\
$\left(SU(2) \times SU(2)\right) \rtimes \mathbb{Z}_2$ & $\mathfrak{su}(2) \times \mathfrak{su}(2)$ & $D_8$ & $\mathbb{Z}_2$ & - \\
$\left(O(2) \times O(2)\right) \rtimes \mathbb{Z}_2$ & $\mathfrak{t}^2$ & $D_8$ & $D_8$ & - \\
\hline
\end{tabular} \\
\caption{Compact Lie groups $G$ of rank two and their corresponding Lie algebra, Weyl group $W_G \subset GL(2,\mathbb{Z})$, group of components $\pi_0(G)$ and fundamental group $\pi_1(G)$.} \label{Table:Lie-Weyl}
\end{center}
\end{table}

The Lie group $G$ is connected if and only if its group of components $\pi_0(G)$ is trivial, and simply-connected if and only if its fundamental group $\pi_1(G)$ is trivial.
The compact rank two Lie Groups in Table \ref{Table:Lie-Weyl} are grouped into blocks as follows.
The first block of four Lie groups are all semi-simple, connected, simply-connected compact Lie groups -- these are the only simply-connected compact rank two Lie groups. The next block consists of the products of $\mathbb{T}$ and $SU(2)$ (excluding $SU(2) \times SU(2)$ which was included in the first block). Note that $\mathbb{T} \times SU(2)$ is semisimple, but $\mathbb{T}^2$ is not usually regarded as semisimple since it is abelian. The groups listed in the third block are all compact connected rank two Lie groups which are (double-, triple-, or quadruple-)covered by groups in the first two blocks. The dashed line separates the Lie groups for which the Weyl group $W_G \subset GL(2,\mathbb{Z})$ do not already appear in the first two blocks, and those for which $W_G$ does appear. The latter Lie groups are not considered in this paper.

In the final block we have listed non-connected compact rank two Lie groups. We do not list Lie groups which are products $G \times H$ where $G$ is a compact, connected rank two Lie group and $H$ is an arbitrary finite group, as the Weyl group in this case is just $W_G$, the Weyl group for $G$, and its group of connected components is $H$. The only semidirect product $\mathbb{T}^2 \rtimes \Gamma$, $\Gamma \subset GL(2,\mathbb{Z})$, which appears in this list is for $\Gamma = \mathbb{Z}_2^{(3)} \subset GL(2,\mathbb{Z})$. This is because all other subgroups $\Gamma \subset GL(2,\mathbb{Z})$ either act trivially by conjugation or the only non-trivial action of $\Gamma$ is that given by $\mathbb{Z}_2^{(3)}$ (as is the case for $\mathbb{Z}_4$ and $D_4^{(2)}$). Thus the semidirect product $\mathbb{T}^2 \rtimes \Gamma$ reduces to a product of $\mathbb{T}^2$ (or $\mathbb{T}^2 \rtimes \mathbb{Z}_2^{(3)}$) by some $\Gamma' \subset GL(2,\mathbb{Z})$. Finally, for the last two groups listed, the $\mathbb{Z}_2$-action interchanges the two components in the product.

In this paper we determine spectral measures associated to the compact connected rank two Lie groups listed in the first three blocks, as these are sufficient to exhaust all finite subgroups of $GL(2,\mathbb{Z})$, up to conjugation in $GL(2,\mathbb{Z})$.
For our purposes it is sufficient only to know the embedding of the Weyl group as a subgroup of $GL(2,\mathbb{Z})$. The irreducible characters and corresponding representation graphs will be constructed in the proceeding sections from knowledge of $\Gamma \subset GL(2,\mathbb{Z})$.
For one of the groups $G$ contained in the fourth block in Table \ref{Table:Lie-Weyl}, one could obtain the spectral measures in this case by considering its corresponding covering group $H$ in the first two blocks, and determining the spectral measures for representation graphs corresponding to those irreducible representations of $H$ which are fundamental generators (see Section \ref{sect:representation_graphs}) of $G$. Such spectral measures will not be determined in this paper.

The framework used here, which is described in Sections \ref{sect:orbit_functions}-\ref{sect:measures-different_domains}, cannot be used for subgroups $\Gamma \subset SL(2,\mathbb{Z})$, since a fundamental domain $P_+$ of $\mathbb{Z}^2/\Gamma$ is not uniquely defined.
It is still possible to associate a graph $\mathcal{G}$, and hence a spectral measure, to $\Gamma$. One way of doing this is described in Remark \ref{Rem:Z2^1=2.D4^1} for the case of $\mathbb{Z}_2^{(1)}$. However the joint spectrum $\mathfrak{D}$ of the graph $\mathcal{G}$ is equal to the joint spectrum of $D_4^{(1)} \cong \mathbb{Z}_2^{(1)} \rtimes \mathbb{Z}_2^{(2)}$, and the spectral measure over $\mathfrak{D}$ is twice that for $D_4^{(1)}$.
Similar statements can be made for the other subgroups $\Gamma \subset SL(2,\mathbb{Z})$, with the spectral measure over the joint spectrum $\mathfrak{D}$ for $\Gamma$ being twice that for the corresponding dihedral group $D \cong \Gamma \rtimes \mathbb{Z}_2$.
Thus from this perspective the finite subgroups of $SL(2,\mathbb{Z})$ do not give anything new compared with the subgroups in $\mathfrak{G}$.

\subsection{Orbit Functions and characters} \label{sect:orbit_functions}

Symmetric, anti-symmetric orbit functions $C_{\lambda}$, $S_{\lambda}$ respectively are defined for a compact, semisimple Lie group $G$ and are closely related to the Weyl group $W_G$ of $G$ \cite{klimyk/patera:2006, klimyk/patera:2007} (see also \cite{nesterenko/patera/tereszkiewicz:2011}). They are also called $C$-, $S$-functions respectively, since when defined for $SU(2)$ these functions coincide with cosine, sine functions respectively.
When $G$ is a Lie group of rank $n$, then $C_{\lambda}$, $S_{\lambda}$ are functions of $n$ variables which are the set of distinct points in $\mathbb{R}^n$ generated by the action of $W_G$ on $\lambda$. The $n$-tuples $\lambda$ are usually taken to be in the set $P_+ = \{ \sum_{i=1}^n \lambda_i \Lambda_i | \, 0 \leq \lambda_i \in \mathbb{Z} \}$, where $\Lambda_i$ are the fundamental weights of $G$. However, the definition extends to $\lambda \in P = \{ \sum_{i=1}^n \lambda_i \Lambda_i | \, \lambda_i \in \mathbb{Z} \} = \mathbb{Z}^n$.

We will define $C$-, $S$-functions for any finite subgroup $\Gamma$ of $GL(2,\mathbb{Z})$, by replacing the Weyl group $W_G$ with $\Gamma$. As discussed in Section \ref{sect:rank2LieGroups}, such a $\Gamma$ is in fact the Weyl group for a compact, connected rank two Lie group $G$, although the existence of such a Lie group is not necessary for the definition. Such a Lie group $G$ will not necessarily be semisimple, and thus extends the definitions of \cite{klimyk/patera:2006, klimyk/patera:2007} of orbit functions to non-semisimple compact, connected Lie groups.
For $\lambda \in \mathbb{Z}^n$ and $\theta \in \mathbb{R}^n$,
\begin{equation} \label{def:C-fns}
C_{\lambda}(\theta) := \sum_{\gamma \in \Gamma} e^{2 \pi i \langle \gamma \lambda, \theta \rangle}, \qquad S_{\lambda}(\theta) := \sum_{\gamma \in \Gamma} \mathrm{det}(\gamma) e^{2 \pi i \langle \gamma \lambda, \theta \rangle},
\end{equation}
where $\langle \, \cdot \, , \, \cdot \, \rangle$ is the Euclidean inner product on $\mathbb{R}^n$. For $\Gamma \subset SL(2,\mathbb{Z})$, the definitions of $C_{\lambda}$ and $S_{\lambda}$ coincide, since $\mathrm{det}(\gamma) = 1$ for all $\gamma \in \Gamma$.

In the case where $\Gamma = W_G$ is the Weyl group for some compact, semisimple Lie group $G$, and $\lambda \in P_+$, the definition of $C_{\lambda}$ given above is not the usual definition of $C$-function used by Patera et al., which is $\widetilde{C}_{\lambda}(x) = |\textrm{Stab}_{\lambda}|^{-1} \sum_{w \in W_G} e^{2 \pi i \langle w \lambda, x \rangle}$.
Here $\textrm{Stab}_{\lambda} = \{ w \in W_G | w \lambda = \lambda \}$ is the stabilizer group of $\lambda$. Note however that for $\lambda \in P_{++} = \{ \lambda = \sum_{i=1}^n \lambda_i \Lambda_i | \, 0 < \lambda_i \in \mathbb{Z} \}$, i.e. for $\lambda$ in the interior of $P_+$, $|\textrm{Stab}_{\lambda}|=1$ so that in that case $\widetilde{C}_{\lambda} = C_{\lambda}$.

For any $\Gamma \subset GL(2,\mathbb{Z})$ these orbit functions are orthogonal over $\mathbb{T}^n$ (c.f. \cite{moody/patera:2006} in the case where $\Gamma$ is the Weyl group for some compact, semisimple Lie group):
\begin{equation} \label{eqn:orthogonality-CS}
\int_{\mathbb{T}^n} C_{\lambda}(\theta) \overline{C_{\mu}(\theta)} \mathrm{d}\theta \;\; = \;\; \delta_{\lambda,\mu} |\Gamma| \;\; = \;\; \int_{\mathbb{T}^n} S_{\lambda}(\theta) \overline{S_{\mu}(\theta)} \mathrm{d}\theta,
\end{equation}
for $\lambda,\mu \in P_+$, where $\theta = (\theta_1,\theta_2,\ldots,\theta_n)$, $e^{2\pi i \theta} \in \mathbb{T}^n$ and $\mathrm{d}\theta = \mathrm{d}\theta_1 \cdots \mathrm{d}\theta_n$ for $\mathrm{d}\theta_i$ the uniform Lebesque measure over $\mathbb{T}$.
These equalities follow from the orthogonality $\int_{\mathbb{T}} u^m \mathrm{d}u = \delta_{m,0}$ of $\mathbb{T}$, since
$$\int_{\mathbb{T}^n} C_{\lambda}(\theta) \overline{C_{\mu}(\theta)} \mathrm{d}\theta \;\; = \;\; \sum_{\gamma,\gamma' \in \Gamma} \int_{\mathbb{T}^n} e^{2\pi i \langle \gamma\lambda-\gamma'\mu, \theta \rangle} \mathrm{d}\theta \;\; = \;\; \sum_{\gamma,\gamma' \in \Gamma} \, \delta_{\gamma\lambda,\gamma'\mu} \;\; = \;\; |\Gamma| \, \delta_{\lambda,\mu},$$
and
\begin{eqnarray*}
\int_{\mathbb{T}^n} S_{\lambda}(\theta) \overline{S_{\mu}(\theta)} \mathrm{d}\theta & = & \sum_{\gamma \in \Gamma} \mathrm{det}(\gamma) \sum_{\gamma' \in \Gamma} \mathrm{det}(\gamma') \int_{\mathbb{T}^n} e^{2\pi i \langle \gamma\lambda-\gamma'\mu, \theta \rangle} \mathrm{d}\theta \\
& = & \sum_{\gamma \in \Gamma} \mathrm{det}(\gamma) \sum_{\gamma' \in \Gamma} \mathrm{det}(\gamma') \, \delta_{\gamma\lambda,\gamma'\mu} \;\; = \;\; |\Gamma| \, \delta_{\lambda,\mu},
\end{eqnarray*}
where we also use the fact that for $\lambda,\mu \in P_+$, $\gamma \lambda = \gamma' \mu \Leftrightarrow \lambda = \mu$ and $\gamma = \gamma'$.

For $\Gamma \in \mathfrak{G}$, we have $S_{\lambda} = 0$ for $\lambda \in P_+ \setminus P_{++}$, i.e. for $\lambda$ lying on the boundary of $P_+$, since for every element $\gamma$ with $\mathrm{det}(\gamma)=1$ (i.e. $\gamma$ is a rotation), there exists an element $\gamma'$ with $\mathrm{det}(\gamma')=-1$ (i.e. $\gamma'$ is a reflection) such that $\gamma' \gamma \lambda = \gamma \lambda$. Thus in the summation for $S_{\lambda}$ the terms for $\gamma' \gamma$ and $\gamma$ have opposite signs and thus cancel.

When $\Gamma \neq \mathbb{Z}_2^{(3)}$, we denote by $\varrho$ the point in $P_{++}$ such that $\langle \varrho, \varrho \rangle \leq \langle \lambda, \lambda \rangle$ for all $\lambda \in P_{++}$. For $\mathbb{Z}_2^{(3)}$ there are two points $(1,0)$, $(0,-1)$ in $P_{++}$ which satisfy this condition, and we take $\varrho$ to be $(1,0)$.
Then we define a formal character $\chi_{\lambda}$ for $\lambda \in P_+$ by
\begin{equation} \label{eqn:character=S/S}
\chi_{\lambda}(\omega_1,\omega_2) = S_{\lambda + \varrho}(\theta_1,\theta_2)/S_{\varrho}(\theta_1,\theta_2),
\end{equation}
where $\omega_j = e^{2\pi i \theta_j} \in \mathbb{T}$, $j=1,2$. Note that $\chi_{\lambda}$ is non-zero only for $\lambda \in P_{++}$, since for $\lambda' \in P_+ \setminus P_{++}$, $\chi_{\lambda'} = 0$ (since $S_{\lambda'} = 0$).

We observe that since $\Gamma = W_G$ is the Weyl group for a compact, connected Lie group $G$ listed in the first three blocks of Table \ref{Table:Lie-Weyl}, equation (\ref{eqn:character=S/S}) is just the Weyl character formula for $G$.
In the case of the other connected Lie groups $G$ listed in Table \ref{Table:Lie-Weyl}, since these groups are each covered by a group $H$ from the first two blocks, (\ref{eqn:character=S/S}) still yields the character for irreducible representations of the Lie group, but only for the subset of $P_+$ corresponding to those irreducible representations of $H$ which are irreducible representations of $G$.

For non-connected Lie groups, which are all semi-direct products $G = N \rtimes \mathbb{Z}_2$ where $N$ is one of the groups listed in the first two blocks of Table \ref{Table:Lie-Weyl} and is normal in $G$, the irreducible representations $\lambda$ and $\nu(\lambda)$ (for $\nu$ the non-trivial element of $\mathbb{Z}_2$) are either isomorphic or not. In the case where they are, one obtains two irreducible representation $(\lambda,0)$, $(\lambda,1)$ of $G$ whose characters on $\mathbb{T}^2$ are both given by $\chi_{(\lambda,0)}(t) = \chi_{(\lambda,1)}(t) = \chi_{\lambda}(t)$ for $t \in \mathbb{T}^2$. In the case where they are not isomorphic, one obtains a single irreducible representation $\lambda'$ of $G$ whose character on $\mathbb{T}^2$ is given by $\chi_{\lambda'}(t) = \chi_{\lambda}(t) + \chi_{\nu(\lambda)}(t)$.
Thus equation (\ref{eqn:character=S/S}) is the Weyl character formula in the case where $\lambda \cong \nu(\lambda)$, whilst for $\lambda \not \cong \nu(\lambda)$, the Weyl character formula reads $\chi_{\lambda'}(\omega_1,\omega_2) = (S_{\lambda+\varrho}(\theta_1,\theta_2) + S_{\nu(\lambda+\varrho)}(\theta_1,\theta_2))/S_{\varrho}(\theta_1,\theta_2)$.

It was noted in \cite[$\S$5]{evans/pugh:2012i} in the context of the Lie group $G_2$ that there is a connection between the orbit function $S_{\lambda+\varrho}(x)$ and the modular $S$-matrix for the conformal field theory associated to $G_2$ at finite level $k$. More precisely, if we let $x = ((\mu_1+1)/3(k+4),-(\mu_2+1)/(k+4))$, then up to a common scalar multiple, $S_{\lambda+\varrho}(x) = S_{\lambda,\mu}$. Here we are instead using Dynkin labels for $\lambda,\mu \in P_+$. Similarly, if we let $x = ((\mu_1+1)/2(k+2),(\mu_2+1)/2(k+2)), ((2\mu_1+\mu_2+3)/3(k+3),(\mu_1+2\mu_2+3)/3(k+4)), ((\mu_1+\mu_2+2)/2(k+3),(\mu_1+2\mu_2+3)/2(k+3))$ respectively, then up to a common scalar multiple, $S_{\lambda+\varrho}(x) = S_{\lambda,\mu}$ for $SU(2) \times SU(2)$, $SU(3)$, $Sp(2)$ respectively, at finite level $k$, again using Dynkin labels for $\lambda,\mu \in P_+$. Thus for all the semi-simple, connected, simply-connected compact rank two Lie groups the modular $S$-matrix is given by the orbit $S$-function, up to some common scalar multiple which ensures that the $S$-matrix has norm 1.

\subsection{Representation graphs} \label{sect:representation_graphs}

We construct two families of graphs for each subgroup $\Gamma \subset \mathfrak{G}$.
One family is the McKay graphs $\mathcal{G}^{\Gamma}$ for the irreducible representations of a compact, connected Lie group $G$ with Weyl group $\Gamma$, the other is the McKay graphs $\mathcal{H}^{\Gamma}$ for the action of the irreducible representations of $G$ on the irreducible representations of the torus $\mathbb{T}^2$.

For any $\lambda,\mu \in P_+$, $\chi_{\lambda} \chi_{\mu}$ decomposes into a finite sum of characters $\chi_{\nu}$ for $\nu \in P_+$ since these are characters of a compact, connected Lie group $G$. More explicity,
\begin{equation} \label{eqn:fusion-chi}
\chi_{\lambda} \chi_{\mu} = \sum_{\nu} m^{\lambda}_{\mu\nu} \chi_{\nu},
\end{equation}
where $m^{\lambda}_{\mu\nu} \in \mathbb{N}$ such that $\sum_{\nu} m^{\lambda}_{\mu\nu} < \infty$ for all $\lambda,\mu \in P_+$.

Thus we may form the (infinite) graphs $\mathcal{G}^{\Gamma}_{\lambda}$, the McKay graphs for the irreducible representations $\lambda$ of the compact, connected Lie group $G$, whose vertices correspond to the characters and edges correspond to multiplication by $\chi_{\lambda}$.
The normal matrices $\Delta_{\lambda} = (m^{\lambda}_{\mu\nu})_{\mu,\nu}$, the adjacency matrices of the graphs $\mathcal{G}^{\Gamma}_{\lambda}$, commute since the characters $\chi_{\lambda}$ do. In the case where $\chi_{\lambda}(x)$ is real-valued for all $x \in \mathbb{T}^2$, the matrix $\Delta_{\lambda}$ is self-adjoint.

Now let $\{ \sigma_{(\mu_1,\mu_2)} \}_{\mu_1,\mu_2 \in \mathbb{Z}}$ be the irreducible characters of $\mathbb{T}^2$, where $\sigma_{(\mu_1,\mu_2)}(t_1,t_2) = t_1^{\mu_1}t_2^{\mu_2}$, for $t_i \in \mathbb{T}$, $\mu_1,\mu_2 \in \mathbb{Z}$.
The characters $\chi_{\lambda}$ of $G$ decompose as a finite sum of characters $\sigma_{\nu}$ of $\mathbb{T}^2$ for $\nu \in \mathbb{Z}^2$, and hence $\chi_{\lambda} \sigma_{\mu}$ decomposes into a finite sum of characters $\sigma_{\nu}$, $\nu \in \mathbb{Z}^2$,
\begin{equation} \label{eqn:fusion-chiT}
\chi_{\lambda} \, \sigma_{\mu} = \sum_{\nu} n^{\lambda}_{\mu\nu} \sigma_{\nu},
\end{equation}
where $n^{\lambda}_{\mu\nu} \in \mathbb{N}$ such that $\sum_{\nu} n^{\lambda}_{\mu\nu} < \infty$ for all $\lambda,\mu \in \mathbb{Z}^2$.
By considering the action of $\chi_{\lambda}$ on the $\sigma_{\mu}$ we obtain a second pair of (infinite) graphs, $\mathcal{H}^{\Gamma}_{\lambda}$, whose vertices correspond to the irreducible characters of $\mathbb{T}^2$, and edges correspond to multiplication by $\chi_{\lambda}$. We will label the vertex corresponding to $\sigma_{\nu}$ by $\nu \in \mathbb{Z}^2$.
Again, the normal matrices ${}^{\Gamma} \hspace{-1mm} \Delta_{\lambda} = (n^{\lambda}_{\mu\nu})_{\mu,\nu}$, the adjacency matrices of the graphs $\mathcal{H}^{\Gamma}_{\lambda}$, commute, and in the case where $\chi_{\lambda}(x)$ is real-valued for all $x \in \mathbb{T}^2$, the matrix ${}^{\Gamma} \hspace{-1mm} \Delta_{\lambda}$ is self-adjoint.

We will specify a pair of points $\rho_1,\rho_2 \in P_+$ which we will call fundamental generators, since their characters are generators in the sense that $\chi_{\lambda}$ for any other $\lambda \in P_+$ appears in the decomposition of the product of some powers of $\chi_{\rho_1}$ and $\chi_{\rho_2}$. In fact, $\{ \chi_{\rho_1}, \chi_{\rho_2} \}$ is a system of generators for the algebra of characters in all cases, except for the groups $\mathbb{Z}_2^{(3)}$ and $D_6^{(1)}$. For $\mathbb{Z}_2^{(3)}$ one needs to also use the additional information that $\chi_{\overline{\lambda}}(\omega_1,\omega_2) = \chi_{\lambda}(\overline{\omega_1},\omega_2) \; \left( = \chi_{\lambda}(\omega_1,\overline{\omega_2})\right)$, where $\overline{\lambda} = (\lambda_1,-\lambda_2)$ for $\lambda=(\lambda_1,\lambda_2)$, which corresponds to the automorphism of the graph $\mathcal{G}^{\mathbb{Z}_2^{(3)}}_{\mu}$ given by reflecting the graph about the line $y=-x$ and reversing all orientations (see Figures \ref{Fig-Graph_Z2^3-G10}, \ref{Fig-Graph_Z2^3-G11} for the cases where $\mu=\rho_1,\rho_2$ respectively). For $D_6^{(1)}$ one needs to also use the additional information that $\chi_{\overline{\lambda}}(\omega_1,\omega_2) = \chi_{\lambda}(\overline{\omega_1},\overline{\omega_2})$, where $\overline{\lambda} = (\lambda_1+\lambda_2,-\lambda_2)$, which corresponds to the automorphism of the graph $\mathcal{G}^{D_6^{(1)}}_{\mu}$ given by reflecting the graph about the $x$-axis and reversing all orientations (see Figures \ref{Fig-Graph_D6^1-G10}, \ref{Fig-Graph_D6^1-G11} for the cases where $\mu=\rho_1,\rho_2$ respectively).

We will study (joint) spectral measures for the pair of graphs $(\mathcal{G}^{\Gamma}_{\rho_1},\mathcal{G}^{\Gamma}_{\rho_2})$, and similarly for the pair of graphs $(\mathcal{H}^{\Gamma}_{\rho_1},\mathcal{H}^{\Gamma}_{\rho_2})$.
Equation (\ref{eqn:fusion-chi}) can be interpreted as meaning that the matrix $\Delta_{\rho_i}$ has eigenvector $(\chi_{\nu}(\theta))_{\nu}$ for eigenvalue $\chi_{\rho_i}(\theta)$, $\theta \in [0,2\pi]^2$.
Thus the spectrum of $\Delta_{\rho_i}$ is given by $\chi_{\rho_i}(\mathbb{T}^2)$.
Similarly, from equation (\ref{eqn:fusion-chiT}) we see that the spectrum of ${}^{\Gamma} \hspace{-1mm} \Delta_{\rho_i}$ is also given by $\chi_{\rho_i}(\mathbb{T}^2)$.

\subsection{Spectral measures over different domains} \label{sect:measures-different_domains}

Suppose $A$ is a unital $C^{\ast}$-algebra with state $\varphi$.
If $b \in A$ is a normal operator then there exists a compactly supported probability measure $\nu_b$ on the spectrum $\sigma(b) \subset \mathbb{C}$ of $b$, uniquely determined by its moments
\begin{equation} \label{eqn:moments_normal_operator}
\varphi(b^m b^{\ast n}) = \int_{\sigma(b)} z^m \overline{z}^n \mathrm{d}\nu_b (z),
\end{equation}
for non-negative integers $m$, $n$.
If $a$ is self-adjoint (\ref{eqn:moments_normal_operator}) reduces to
\begin{equation} \label{eqn:moments_selfadjoint}
\varphi(a^m) = \int_{\sigma(a)} x^m \mathrm{d}\nu_a (x),
\end{equation}
with $\sigma(a) \subset \mathbb{R}$, for any non-negative integer $m$.

One can also consider more general measures over the joint spectrum $\sigma(a,b) \subset \sigma(a) \times \sigma(b) \subset \mathbb{C}^2$ of commuting normal operators $a$ and $b$. The abelian $C^{\ast}$-algebra $B$ generated by $a$, $b$ and the identity 1 is isomorphic to $C(X)$, where $X$ is the spectrum of $B$. The joint spectrum is defined as $\sigma(a,b) = \{ (a(x), b(x)) | \, x \in X \}$. In fact, one can identify the spectrum $X$ with its image $\sigma(a,b)$ in $\mathbb{C}^2$, since the map $x \mapsto (a(x), b(x))$ is continuous and injective, and hence a homeomorphism since $X$ is compact \cite{takesaki:2002}.
In the case where the operators $a$, $b$ act on a finite-dimensional Hilbert space, this is the set of all pairs of real numbers $(\lambda_a,\lambda_b)$ for which there exists a non-zero vector $\phi$ such that $a\phi = \lambda_a \phi$, $b\phi = \lambda_b \phi$.
Then there exists a compactly supported probability measure $\widetilde{\nu}_{a,b}$ on $\sigma(a,b)$, the joint spectral measure of $a$, $b$, which, for $a \neq b$, is uniquely determined by its cross moments
\begin{equation} \label{eqn:cross_moments_sa_operators}
\varphi(a^{m_1} a^{\ast n_1} b^{m_2} b^{\ast n_2}) = \int_{\sigma(a,b)} w^{m_1} \overline{w}^{n_1} z^{m_2} \overline{z}^{n_2} \mathrm{d}\widetilde{\nu}_{a,b} (w,z),
\end{equation}
for all non-negative integers $m_i$, $n_i$. In the case where $a$ or $b$ is self-adjoint, it is sufficient to consider the cross moments with $n_1=0$ or $n_2=0$ respectively in (\ref{eqn:cross_moments_sa_operators}) to determine the measure $\widetilde{\nu}_{a,b}$.
For $a,b$ both self-adjoint, the spectral measure for $a$ is given by the pushforward $(p_a)_{\ast}(\widetilde{\nu}_{a,b})$ of the joint spectral measure $\widetilde{\nu}_{a,b}$ under the orthogonal projection $p_a$ onto the spectrum $\sigma(a)$.
In particular, $\nu_a$, $\nu_b$ can be determined by additionally setting $m_2=n_2=0$, $m_1=n_1=0$ respectively.

Let $x_{\lambda} = \chi_{\lambda}(\omega_1,\omega_2)$ and let $\Psi_{\lambda,\mu}$ be the map $(\omega_1,\omega_2) \mapsto (x_{\lambda},x_{\mu})$. We denote by $\mathfrak{D}_{\lambda,\mu}$ the image of $\Psi_{\lambda,\mu}(\mathbb{T}^2)$.
Note that $\mathfrak{D}_{\lambda,\mu} \cong \mathfrak{D}_{\mu,\lambda}$.
Then any $\Gamma$-invariant measure $\varepsilon_{\lambda,\mu}$ on $\mathbb{T}^2$ produces a probability measure $\widetilde{\nu}_{\lambda,\mu}$ on $\mathfrak{D}$ by
\begin{equation} \label{eqn:measures-T2-D_Gamma}
\int_{\mathfrak{D}_{\lambda,\mu}} \psi(x_{\lambda},x_{\mu}) \mathrm{d}\widetilde{\nu}_{\lambda,\mu}(x_{\lambda},x_{\mu}) = \int_{\mathbb{T}^2} \psi(\chi_{\lambda}(\omega_1,\omega_2),\chi_{\mu}(\omega_1,\omega_2)) \mathrm{d}\varepsilon_{\lambda,\mu}(\omega_1,\omega_2),
\end{equation}
for any continuous function $\psi:\mathfrak{D}_{\lambda,\mu} \rightarrow \mathbb{C}$.
Any such measure is uniquely determined by its cross moments $\varsigma_{m_1,n_1,m_2,n_2} = \int_{\mathfrak{D}_{\lambda,\mu}} x_{\lambda}^{m_1} \overline{x_{\lambda}}^{n_1} x_{\mu}^{m_2} \overline{x_{\mu}}^{n_2} \mathrm{d}\widetilde{\nu}_{\lambda,\mu}(x_{\lambda},x_{\mu})$.
Since the $\chi_{\lambda}$ are invariant under the action of $\Gamma$ on $\mathbb{T}^2$, $\mathfrak{D}_{\lambda,\mu}$ is isomorphic to a quotient of $\mathbb{T}^2$ by $\Gamma$. We denote by $C$ a fundamental domain of $\mathbb{T}^2$ under the action of $\Gamma$. The torus then contains $|\Gamma|$ copies of $C$, so that
\begin{equation} \label{eqn:measureT2=|Gamma|C}
\int_{\mathbb{T}^2} \phi(\omega_1,\omega_2) \mathrm{d}\varepsilon_{\lambda,\mu}(\omega_1,\omega_2) = |\Gamma| \int_{C} \phi(\omega_1,\omega_2) \mathrm{d}\varepsilon_{\lambda,\mu}(\omega_1,\omega_2),
\end{equation}
for any $\Gamma$-invariant function $\phi:\mathbb{T}^2 \rightarrow \mathbb{C}$.

As discussed above any probability measure on $\mathfrak{D}_{\lambda,\mu}$ yields a probability measure on $\sigma_{\lambda}$, given by the pushforward $(p_{\lambda})_{\ast}(\widetilde{\nu}_{\lambda,\mu})$ of the joint spectral measure $\widetilde{\nu}_{\lambda,\mu}$ under the orthogonal projection $p_{\lambda}$ onto the spectrum $\sigma(\lambda)$. In particular, when $\psi(x_{\lambda},x_{\mu}) = \widetilde{\psi}(x_{\lambda})$ is only a function of one variable $x_{\lambda}$, then
$$\int_{\mathfrak{D}_{\lambda,\mu}} \widetilde{\psi}(x_{\lambda}) \mathrm{d}\widetilde{\nu}_{\lambda,\mu}(x_{\lambda},x_{\mu}) = \int_{\sigma_{\lambda}} \widetilde{\psi}(x_{\lambda}) \int_{\mathfrak{D}_{\lambda,\mu}(x_{\lambda})} \mathrm{d}\widetilde{\nu}_{\lambda,\mu}(x_{\lambda},x_{\mu}) = \int_{\sigma_{\lambda}} \widetilde{\psi}(x_{\lambda}) \mathrm{d}\nu_{\lambda}(x_{\lambda})$$
where the measure $\mathrm{d}\nu_{\lambda}(x_{\lambda}) = \int_{x_{\mu} \in \mathfrak{D}_{\lambda,\mu}(x_{\lambda})} \mathrm{d}\widetilde{\nu}_{\lambda,\mu}(x_{\lambda},x_{\mu})$ is given by the integral over $x_{\mu} \in \mathfrak{D}_{\lambda,\mu}(x_{\lambda}) = \{ x_{\mu} \in \sigma_{\mu} | \, (x_{\lambda},x_{\mu}) \in \mathfrak{D}_{\lambda,\mu} \}$.
Note that $\mathfrak{D}_{\lambda,\overline{\lambda}} = \sigma_{\lambda}$, thus for non-self-adjoint $\lambda$ the joint spectral measure of $\lambda$, $\overline{\lambda}$ (over $\mathfrak{D}_{\lambda,\overline{\lambda}}$) is in fact the spectral measure of both $\lambda$ and $\overline{\lambda}$.

\section{Joint spectral measures for rank two Lie groups} \label{sect:measures_over_T2}
\subsection{Joint spectral measure for ${}^{\Gamma} \hspace{-1mm} \Delta_{\rho_i}$} \label{sect:measures_over_T2-H}

The adjacency matrices ${}^{\Gamma} \hspace{-1mm} \Delta_{\rho_i}$ can be identified with operators on $\ell^2(\mathbb{Z}) \otimes \ell^2(\mathbb{Z})$, where if $\chi_{\rho_i}$ decomposes into irreducible characters of $\mathbb{T}^2$ as $\chi_{\rho_i} = \sum_{\nu} p^i_{\nu} \sigma_{\nu}$, then ${}^{\Gamma} \hspace{-1mm} \Delta_{\rho_i} = \sum_{\nu} p^i_{\nu} s^{\nu_1} \otimes s^{\nu_2}$, where $\nu = (\nu_1,\nu_2)$, $s$ is the bilateral shift on $\ell^2(\mathbb{Z})$. For $\Omega = (\delta_{j,0})_{j\in\mathbb{Z}}$, we regard $\Omega \otimes \Omega$ as corresponding to the vertex $(0,0)$ whilst $(s^{\mu_1} \otimes s^{\mu_2})(\Omega \otimes \Omega)$ corresponds to the vertex $(\mu_1,\mu_2)$ of $\mathcal{H}_{\rho_i}^{\Gamma}$.
We define a state $\varphi$ on $C^{\ast}(v_Z^1,v_Z^2)$ by $\varphi( \, \cdot \, ) = \langle \, \cdot \, (\Omega \otimes \Omega), \Omega \otimes \Omega \rangle$. Then $\varphi(s^{\lambda_1} \otimes s^{\lambda_2}) = \langle (s^{\lambda_1} \otimes s^{\lambda_2})(\Omega \otimes \Omega), \Omega \otimes \Omega \rangle = \delta_{\lambda_1,0} \, \delta_{\lambda_2,0}$.

Then we have the following result for the joint spectral measure over $\mathbb{T}^2$ of $({}^{\Gamma} \hspace{-1mm} \Delta_{\rho_1}, {}^{\Gamma} \hspace{-1mm} \Delta_{\rho_2})$:

\begin{Thm} \label{thm:measureT2}
The joint spectral measure $\varepsilon$ (over $\mathbb{T}^2$) for the pair of graphs $(\mathcal{H}^{\Gamma}_{\rho_1}, \mathcal{H}^{\Gamma}_{\rho_2})$ is given by the uniform Lebesgue measure $\mathrm{d}\varepsilon(\omega_1,\omega_2) = \mathrm{d}\omega_1 \, \mathrm{d}\omega_2$.
\end{Thm}
\emph{Proof:}
The result follows from the fact that $\int_{\mathbb{T}} \omega^{\lambda} \, \mathrm{d}\omega := \int_0^1 e^{2 \pi i \theta \lambda} \, \mathrm{d}\theta = \delta_{\lambda,0}$, and thus
\begin{eqnarray*}
\lefteqn{ \int_{\mathbb{T}^2} \chi_{\rho_1}(\omega_1,\omega_2)^{m_1} \overline{\chi_{\rho_1}(\omega_1,\omega_2)}^{n_1} \chi_{\rho_2}(\omega_1,\omega_2)^{m_2} \overline{\chi_{\rho_2}(\omega_1,\omega_2)}^{n_2} \mathrm{d}\omega_1 \mathrm{d}\omega_2 } \\
&& \qquad \qquad \qquad \qquad \;\; = \;\; \varphi(({}^{\Gamma} \hspace{-1mm} \Delta_{\rho_1})^{m_1} ({}^{\Gamma} \hspace{-1mm} \Delta_{\rho_1}^T)^{n_1} ({}^{\Gamma} \hspace{-1mm} \Delta_{\rho_2})^{m_2} ({}^{\Gamma} \hspace{-1mm} \Delta_{\rho_2}^T)^{n_2}).
\end{eqnarray*}
See e.g. \cite[Theorem 2]{evans/pugh:2009v} and \cite[Theorem 3.1]{evans/pugh:2012i} for explicit details in the cases of $SU(3)$ and $G_2$ respectively.
\hfill
$\Box$

In fact, by a similar proof the measure $\varepsilon$ given in Theorem \ref{thm:measureT2} is the joint spectral measure over $\mathbb{T}^2$ for the pair of representation graphs graphs $(\mathcal{H}^{\Gamma}_{\lambda}, \mathcal{H}^{\Gamma}_{\mu})$ for any pair $\lambda,\mu \in P_+$ (note that $\lambda,\mu$ are irreducible representations of a connected Lie group $G$ from the first three blocks of Table \ref{Table:Lie-Weyl} such that the Weyl group of $G$ is $\Gamma$). Thus the spectral measure over $\mathbb{T}^2$ is independent of the choice of $\lambda,\mu \in P_+$.

We now consider the joint spectral measure $\widetilde{\nu}$ over $\mathfrak{D} := \mathfrak{D}_{\rho_1,\rho_2}$, the joint spectrum of the commuting normal operators ${}^{\Gamma} \hspace{-1mm} \Delta_{\rho_1}$, ${}^{\Gamma} \hspace{-1mm} \Delta_{\rho_2}$.
We denote by $J_{\Gamma}$ be the Jacobian $J_{\Gamma} = \mathrm{det}(\partial(x,y)/\partial(\theta_1,\theta_2))$ for the change of variables $x := x_{\rho_1}$, $y := x_{\rho_2}$.
Over $\mathfrak{D}$, we thus obtain
\begin{equation} \label{eqn:integral_overD}
\int_{C} \psi(\chi_{\rho_1}(\omega_1,\omega_2),\chi_{\rho_2}(\omega_1,\omega_2)) \mathrm{d}\omega_1 \, \mathrm{d}\omega_2 = \int_{\mathfrak{D}} \psi(x,y) |J_{\Gamma}(x,y)|^{-1} \mathrm{d}x \, \mathrm{d}y,
\end{equation}
Then from Theorem \ref{thm:measureT2} and (\ref{eqn:measureT2=|Gamma|C}) we obtain

\begin{Thm} \label{thm:joint_measure-D}
The joint spectral measure $\widetilde{\nu}$ (over $\mathfrak{D}$) for the pair of graphs $(\mathcal{H}^{\Gamma}_{\rho_1}, \mathcal{H}^{\Gamma}_{\rho_2})$ is
$$\mathrm{d}\widetilde{\nu}(x,y) = \frac{|\Gamma|}{|J_{\Gamma}(x,y)|} \, \mathrm{d}x \, \mathrm{d}y.$$
\end{Thm}

We determine the Jacobian $J = J_{\Gamma}$ for each group $\Gamma$ in the following sections.

\subsection{Joint spectral measure for $\Delta_{\rho_i}$} \label{sect:measures_over_T2-G}

The joint spectral measure for $\Delta_{\rho_i}$ (over $\mathbb{T}^2$) for $SU(3)$, $Sp(2)$, $G_2$ was shown to be given by $\mathrm{d}\varepsilon(\omega_1,\omega_2) = a_{\Gamma} |J_{\Gamma}(\theta_1,\theta_2)|^2/16 \pi^4 \, \mathrm{d}\omega_1 \, \mathrm{d}\omega_2$, where $\mathrm{d}\omega$ is the uniform Lebesque measure over $\mathbb{T}$ and $a_{\Gamma}=1$ for $\Gamma=D_8, D_{12}$ (the Weyl groups of $Sp(2),G_2$ respectively), and $a_{D_6^{(2)}}=2$ for $SU(3)$ \cite{evans/pugh:2009v, evans/pugh:2010i, evans/pugh:2012i, evans/pugh:2012ii, evans/pugh:2012iii}. Over $\mathfrak{D}$, the joint spectral measure is thus given by $\mathrm{d}\widetilde{\nu}(x,y) = a_{\Gamma} |J_{\Gamma}(x,y)|/16 \pi^4 \, \mathrm{d}x \, \mathrm{d}y$. The same results hold for $SU(2) \times SU(2)$, where $a_{\Gamma} = 1$ \cite{banica/bisch:2007, evans/pugh:2009v}. It can be shown from Sections \ref{sect:Z0}-\ref{sect:Z2(3)} that the same results hold for the rank two Lie groups $\mathbb{T}^2=U(1) \times U(1)$, $\mathbb{T} \times SU(2)$ and $U(2)$, with $a_{\Gamma}=1$ in each case. However, for $SO(4)$, considered in Section \ref{sect:D4^2}, the joint spectral measure for $\Delta_{\rho_i}$ (over $\mathfrak{D}$) is given by $(1+y)^{-2}|J_{\Gamma}(x,y)|/16\pi^2 \, \mathrm{d}x \, \mathrm{d}y$.

This leads to the question of why there is a difference for $SO(4)$ (and indeed for $SU(3)$ where $a_{\Gamma}=2$) and whether there is a consistent description for the spectral measure for $\Delta_{\rho_i}$ (over $\mathbb{T}^2$ or $\mathfrak{D}$) for any rank 2 Lie groups in terms of some object which is naturally associated to the group. It turns out that computing the orbit function $S_{\varrho}(\theta)$ for each group we obtain that $4\pi^2|S_{\varrho}(\theta)| = a_{\Gamma}|J_{\Gamma}(\theta_1,\theta_2)|$ for all the Lie groups $G$ discussed in the preceding paragraph, except for $SO(4)$ where we have $4\pi^2|S_{\varrho}(\theta)| = (1+y)^{-1}|J_{\Gamma}(\theta_1,\theta_2)|$. Thus, for all the rank two Lie groups discussed in the preceding paragraph, the spectral measure for $\Delta_{\rho_i}$ over $\mathbb{T}^2$ is given by $\mathrm{d}\varepsilon(\omega_1,\omega_2) = |\Gamma|^{-1} \, |S_{\varrho}(\theta)|^2 \, \mathrm{d}\omega_1 \, \mathrm{d}\omega_2$. The spectral measure over $\mathfrak{D}$ is $\mathrm{d}\widetilde{\nu}(x,y) = |S_{\varrho}(\theta)|^2 \, |J_{\Gamma}(x,y)|^{-1} \, \mathrm{d}x \, \mathrm{d}y$, where we now write $S_{\varrho}(\theta)$ in terms of the $\Gamma$-invariant variables $x$, $y$.
This leads naturally to the following conjecture:

\begin{Conj} \label{conj:spectral_measure-orbitS}
The joint spectral measure $\varepsilon$ (over $\mathbb{T}^2$) for the pair of graphs $(\mathcal{G}^{\Gamma}_{\rho_1}, \mathcal{G}^{\Gamma}_{\rho_2})$ is
$$\mathrm{d}\varepsilon(\omega_1,\omega_2) = \frac{1}{|\Gamma|} \, |S_{\varrho}(\theta)|^2 \, \mathrm{d}\omega_1 \, \mathrm{d}\omega_2.$$
The joint spectral measure $\widetilde{\nu}$ (over $\mathfrak{D}$) for the pair of graphs $(\mathcal{G}^{\Gamma}_{\rho_1}, \mathcal{G}^{\Gamma}_{\rho_2})$ is
$$\mathrm{d}\widetilde{\nu}(x,y) = \frac{|S_{\varrho}(\theta)|^2}{|J_{\Gamma}(x,y)|} \, \mathrm{d}x \, \mathrm{d}y.$$
\end{Conj}

\section{$\mathbb{Z}_0$: $\mathbb{T}^2=U(1) \times U(1)$} \label{sect:Z0}

\begin{figure}[tb]
\begin{minipage}[t]{7.9cm}
\begin{center}
  \includegraphics[width=40mm]{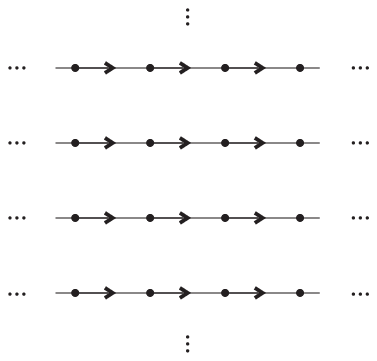}
 \caption{Infinite graph $\mathcal{G}^{\mathbb{Z}_0}_{\rho_1}$ for $\mathbb{T}$} \label{Fig-Graph_Z0-G10}
\end{center}
\end{minipage}
\hfill
\begin{minipage}[t]{7.9cm}
\begin{center}
  \includegraphics[width=40mm]{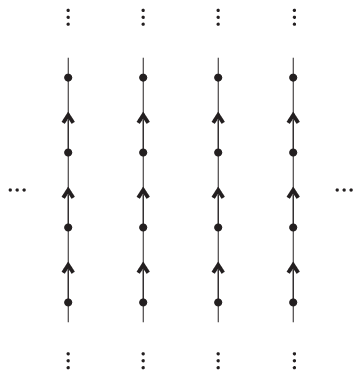}
 \caption{Infinite graph $\mathcal{G}^{\mathbb{Z}_0}_{\rho_2}$ for $\mathbb{T}$} \label{Fig-Graph_Z0-G01}
\end{center}
\end{minipage}
\end{figure}

The trivial subgroup $\mathbb{Z}_0 \in GL(2,\mathbb{Z})$ is generated by the identity matrix. We choose fundamental generators $\rho_1 = (1,0)$, $\rho_2 = (0,1)$.
Since the action of $\mathbb{Z}_0$ on $\mathbb{T}^2$ is trivial, the graphs $\mathcal{G}^{\mathbb{Z}_0}_{\rho}$ and
$\mathcal{H}^{\mathbb{Z}_0}_{\rho}$ are the same. The graphs $\mathcal{G}^{\mathbb{Z}_0}_{\rho}$ for $\rho = \rho_1, \rho_2$ are illustrated in Figures \ref{Fig-Graph_Z0-G10}-\ref{Fig-Graph_Z0-G01}.
Let
$x := x_{\rho_1} = \omega_1$, $y := x_{\rho_2} = \omega_2$, and denote by $\Psi$ be the map $\Psi_{\rho_1,\rho_2}: (\omega_1,\omega_2) \mapsto (x,y)$.
Then $\mathfrak{D} := \mathfrak{D}_{\rho_1,\rho_2} = \mathbb{T}^2$ is the joint spectrum $\sigma(\Delta_{\rho_1},\Delta_{\rho_2})$ of the commuting normal operators $\Delta_{\rho_1}$, $\Delta_{\rho_2}$.
Under the identification $x = \omega_1 = e^{2\pi i \theta_1}$, $y = \omega_2 = e^{2\pi i \theta_2}$, the Jacobian is $J_{\mathbb{Z}_0} = 4 \pi^2 e^{2\pi i(\theta_1+\theta_2)} = 4 \pi^2 xy$.
By integrating $|\Gamma| \, |J(x,y)_{\mathbb{Z}_0}^{-1}| = 1/4\pi^2$ over $y$, $x$ respectively we obtain the spectral measure $\nu_{\rho_1}$, $\nu_{\rho_2}$ respectively for $\Delta_{\rho_1}$, $\Delta_{\rho_2}$.
Since $\int_{\mathbb{T}} x^{m} \mathrm{d}x = 2\pi \delta_{m,0}$, we obtain that the spectral measure $\nu_{\rho_j}$ (over $\chi_{\rho_j}(\mathbb{T}^2) = \mathbb{T}$) for the graph $\mathcal{G}^{\mathbb{Z}_0}_{\rho_j}$, $j=1,2$, is given by
$\mathrm{d}\nu_{\rho_j}(x) = (2\pi)^{-1} \, \mathrm{d}x$.

The graph $\mathcal{G}^{\mathbb{Z}_0}_{\rho_j}$, $j=1,2$, is given by an infinite number of copies of the representation graph $\mathcal{G}^0_{\rho} = \mathcal{H}^0_{\rho}$ for $\mathbb{T}$ (which has Weyl group 0), where $\rho$ is the fundamental representation of $\mathbb{T}$. In particular the connected component of the distinguished vertex $(0,0)$ of $\mathcal{G}^{\mathbb{Z}_0}_{\rho_j}$ is the representation graph $\mathcal{G}^0_{\rho}$ for $\mathbb{T}$.

\section{$\mathbb{Z}_2^{(2)}$: $\mathbb{T} \times SU(2)$} \label{sect:Z2^2}

\begin{figure}[tb]
\begin{minipage}[t]{7.9cm}
\begin{center}
  \includegraphics[width=40mm]{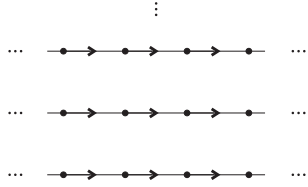}
 \caption{Infinite graph $\mathcal{G}^{\mathbb{Z}_2^{(2)}}_{\rho_1}$ for $\mathbb{T} \times SU(2)$} \label{Fig-Graph_Z2^2-G10}
\end{center}
\end{minipage}
\hfill
\begin{minipage}[t]{7.9cm}
\begin{center}
  \includegraphics[width=40mm]{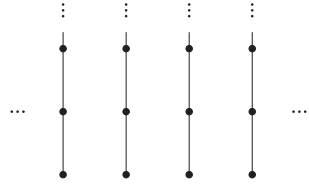}
 \caption{Infinite graph $\mathcal{G}^{\mathbb{Z}_2^{(2)}}_{\rho_2}$ for $\mathbb{T} \times SU(2)$} \label{Fig-Graph_Z2^2-G01}
\end{center}
\end{minipage}
\end{figure}

\begin{figure}[tb]
\begin{minipage}[t]{7.9cm}
\begin{center}
  \includegraphics[width=40mm]{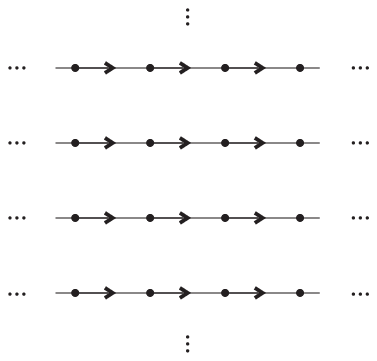}
 \caption{Infinite graph $\mathcal{H}^{\mathbb{Z}_2^{(2)}}_{\rho_1}$ for $\mathbb{T} \times SU(2)$} \label{Fig-Graph_Z2^2-WG10}
\end{center}
\end{minipage}
\hfill
\begin{minipage}[t]{7.9cm}
\begin{center}
  \includegraphics[width=40mm]{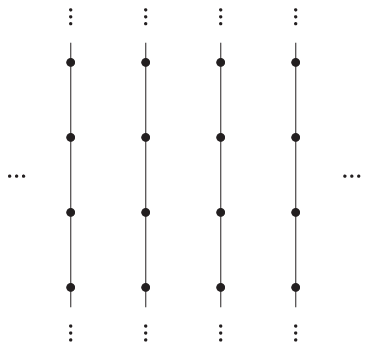}
 \caption{Infinite graph $\mathcal{H}^{\mathbb{Z}_2^{(2)}}_{\rho_2}$ for $\mathbb{T} \times SU(2)$} \label{Fig-Graph_Z2^2-WG01}
\end{center}
\end{minipage}
\end{figure}

\begin{figure}[tb]
\begin{minipage}[t]{7.9cm}
\begin{center}
  \includegraphics[width=55mm]{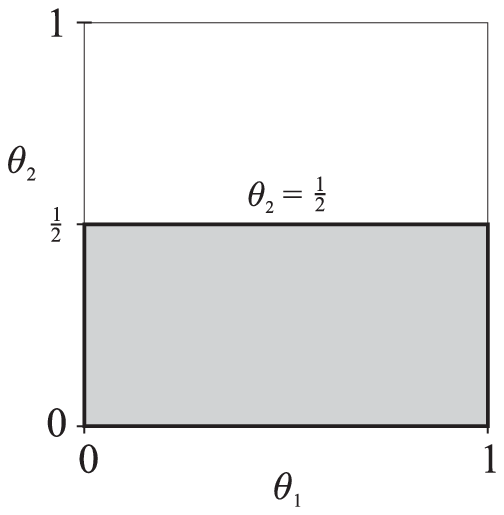}\\
 \caption{A fundamental domain of $\mathbb{T}^2/\mathbb{Z}_2^{(2)}$ for $\mathbb{T} \times SU(2)$.} \label{fig:DomainC-Z2^2}
\end{center}
\end{minipage}
\hfill
\begin{minipage}[t]{7cm}
\begin{center}
  \includegraphics[width=45mm]{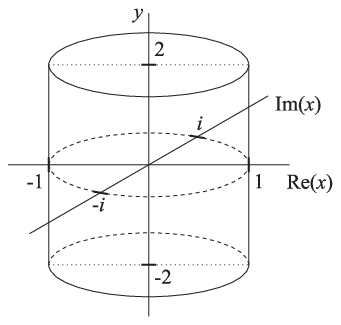}\\
 \caption{The domain $\mathfrak{D} = \Psi(C)$ for \mbox{$\mathbb{T} \times SU(2)$}.} \label{fig:DomainD-Z2^2}
\end{center}
\end{minipage}
\end{figure}

The subgroup $\mathbb{Z}_2^{(2)} \in GL(2,\mathbb{Z})$ is generated by the matrix $T_2$ of Section \ref{sect:intro}.
It is the Weyl group of the connected compact Lie group $\mathbb{T} \times SU(2) = U(1) \times SU(2)$. The spectral measures for $SU(2)$ were studied in \cite{banica/bisch:2007, evans/pugh:2009v}.
We choose fundamental generators $\rho_1 = (1,0)$, $\rho_2 = (0,1)$.
The graphs $\mathcal{G}^{\mathbb{Z}_2^{(2)}}_{\rho}$, $\mathcal{H}^{\mathbb{Z}_2^{(2)}}_{\rho}$ for $\rho = \rho_1, \rho_2$ are illustrated in Figures \ref{Fig-Graph_Z2^2-G10}-\ref{Fig-Graph_Z2^2-WG01}.
Let
\begin{equation} \label{eqn:x,y-Z2^2}
x := x_{\rho_1} = \omega_1, \qquad y := x_{\rho_2} = \omega_2 + \omega_2^{-1},
\end{equation}
and denote by $\Psi$ be the map $\Psi_{\rho_1,\rho_2}: (\omega_1,\omega_2) \mapsto (x,y)$.
A fundamental domain of $\mathbb{T}^2/\mathbb{Z}_2^{(2)}$ is illustrated in Figure \ref{fig:DomainC-Z2^2}.
Then $\mathfrak{D} := \mathfrak{D}_{\rho_1,\rho_2} = \mathbb{T} \times [-2,2]$ is the joint spectrum $\sigma(\Delta_{\rho_1},\Delta_{\rho_2})$ of the commuting normal operators $\Delta_{\rho_1}$, $\Delta_{\rho_2}$. Similarly, $\mathfrak{D}$ is also the joint spectrum of the commuting normal operators ${}^{\mathbb{Z}_2^{(2)}} \hspace{-1mm} \Delta_{\rho_1}$, ${}^{\mathbb{Z}_2^{(2)}} \hspace{-1mm} \Delta_{\rho_2}$. Note that $\Delta_{\rho_2}$, ${}^{\mathbb{Z}_2^{(2)}} \hspace{-1mm} \Delta_{\rho_2}$ are in fact self-adjoint.
Under the change of variables $x = \omega_1$, $y = \omega_2 + \omega_2^{-1}$, the Jacobian is given by
$J_{\mathbb{Z}_2^{(2)}} = -8 \pi^2 i e^{2 \pi i \theta_1} \sin(2 \pi \theta_2) = -4 \pi^2 \omega_1 (\omega_2 - \overline{\omega_2})$. The Jacobian is complex-valued and vanishes in $\mathbb{T}^2$ only on the boundaries of the images of the fundamental domain $C$ under $\mathbb{Z}_2^{(2)}$.
Now $J_{\mathbb{Z}_2^{(2)}}^2 = 16 \pi^4 \omega_1^2 (\omega_2^2 -2 + \overline{\omega_2}^2)$, which is invariant under the action of $\mathbb{Z}_2^{(2)}$ on $\mathbb{T}^2$. We write $J_{\mathbb{Z}_2^{(2)}}^2$ in terms of the $\mathbb{Z}_2^{(2)}$-invariant elements $x$, $y$ as $J_{\mathbb{Z}_2^{(2)}}^2 = - 16 \pi^4 x^2 (4-y^2)$. Since $4-y^2 \geq 0$ for all $y \in [-2,2]$, we can write $J_{\mathbb{Z}_2^{(2)}}$ in terms of $x$, $y$ as
$J_{\mathbb{Z}_2^{(2)}} = 4 \pi^2 i x \sqrt{4-y^2}$.

\subsection{Spectral measure for $\mathcal{H}^{\mathbb{Z}_2^{(2)}}_{\rho}$ for $\mathbb{T} \times SU(2)$} \label{sect:measure-rho-Z2^2-H}

By integrating $|\Gamma| \, |J_{\mathbb{Z}_2^{(2)}}(x,y)^{-1}| = (2\pi^2 \sqrt{4-y^2})^{-1}$ over $y$, $x$ respectively we obtain the spectral measure $\nu_{\rho_1}$, $\nu_{\rho_2}$ respectively for ${}^{\mathbb{Z}_2^{(2)}} \hspace{-1mm} \Delta_{\rho_1}$, ${}^{\mathbb{Z}_2^{(2)}} \hspace{-1mm} \Delta_{\rho_2}$.
Since $\int_{-2}^2 \sqrt{4-y^2}^{-1} \mathrm{d}y = \pi$, the spectral measure $\nu_{\rho_1}$ (over $\chi_{\rho_1}(\mathbb{T}^2) = \mathbb{T}$) for the graph $\mathcal{H}^{\mathbb{Z}_2^{(2)}}_{\rho_1}$ for $\mathbb{T} \times SU(2)$ is given by
$\mathrm{d}\nu_{\rho_1}(x) = (2\pi)^{-1} \, \mathrm{d}x$.
Since $\int_{\mathbb{T}} 1 \, \mathrm{d}y = 2\pi$, the spectral measure $\nu_{\rho_2}$ (over $\chi_{\rho_2}(\mathbb{T}^2) = [-2,2]$) for the graph $\mathcal{H}^{\mathbb{Z}_2^{(2)}}_{\rho_1}$ for $\mathbb{T} \times SU(2)$ is given by
$\mathrm{d}\nu_{\rho_2}(y) = (\pi\sqrt{4-y^2})^{-1} \, \mathrm{d}y$.
The graph $\mathcal{H}^{\mathbb{Z}_2^{(2)}}_{\rho_1}$ is given by an infinite number of copies of the representation graph $\mathcal{G}^0_{\rho} = \mathcal{H}^0_{\rho}$ for $\mathbb{T}$.

\subsection{Spectral measure for $\mathcal{G}^{\mathbb{Z}_2^{(2)}}_{\rho}$ for $\mathbb{T} \times SU(2)$}

The graph $\mathcal{G}^{\mathbb{Z}_2^{(2)}}_{\rho_1}$ is also given by an infinite number of copies of the representation graph $\mathcal{G}^0_{\rho}$ for $\mathbb{T}$, thus the spectral measure $\nu_{\rho_1}$ (over $\chi_{\rho_1}(\mathbb{T}^2) = \mathbb{T}$) for the graph $\mathcal{G}^{\mathbb{Z}_2^{(2)}}_{\rho_1}$ for $\mathbb{T} \times SU(2)$ is given by
$\mathrm{d}\nu_{\rho_1}(x) = (2\pi)^{-1} \, \mathrm{d}x$.
The graph $\mathcal{G}^{\mathbb{Z}_2^{(2)}}_{\rho_2}$ is given by an infinite number of copies of the representation graph $\mathcal{G}^{\langle -1 \rangle}_{\rho}$ for $SU(2)$ (which has Weyl group $\mathbb{Z}_2 = \langle -1 \rangle \subset GL(1,\mathbb{Z})$). In particular the connected component of the distinguished vertex $\ast$ of $\mathcal{G}^{\mathbb{Z}_2^{(2)}}_{\rho_2}$, the vertex with lowest Perron-Frobenius weight (which in this case is the apex vertex, i.e. the vertex in the bottom left corner in Figure \ref{Fig-Graph_Z2^2-G01}), is the representation graph $\mathcal{G}^{\langle -1 \rangle}_{\rho}$ for $SU(2)$. Thus the spectral measure $\nu_{\rho_2}$ (over $\chi_{\rho_2}(\mathbb{T}^2) = [-2,2]$) for the graph $\mathcal{G}^{\mathbb{Z}_2^{(2)}}_{\rho_1}$ for $\mathbb{T} \times SU(2)$ is given \cite{dykema/voiculescu/nica:1992} by
$\mathrm{d}\nu_{\rho_2}(y) = \frac{1}{2\pi}\sqrt{4-y^2} \, \mathrm{d}y$.

\section{$\mathbb{Z}_2^{(3)}$: $U(2)$} \label{sect:Z2(3)}

\begin{figure}[tb]
\begin{minipage}[t]{7.9cm}
\begin{center}
  \includegraphics[width=40mm]{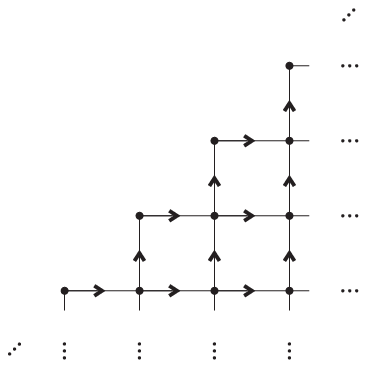}
 \caption{Infinite graph $\mathcal{G}^{\mathbb{Z}_2^{(3)}}_{\rho_1}$ for $U(2)$} \label{Fig-Graph_Z2^3-G10}
\end{center}
\end{minipage}
\hfill
\begin{minipage}[t]{7.9cm}
\begin{center}
  \includegraphics[width=40mm]{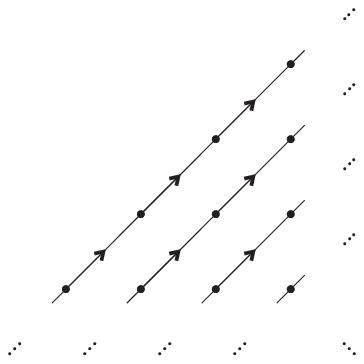}
 \caption{Infinite graph $\mathcal{G}^{\mathbb{Z}_2^{(3)}}_{\rho_2}$ for $U(2)$} \label{Fig-Graph_Z2^3-G11}
\end{center}
\end{minipage}
\end{figure}

\begin{figure}[tb]
\begin{minipage}[t]{7.9cm}
\begin{center}
  \includegraphics[width=40mm]{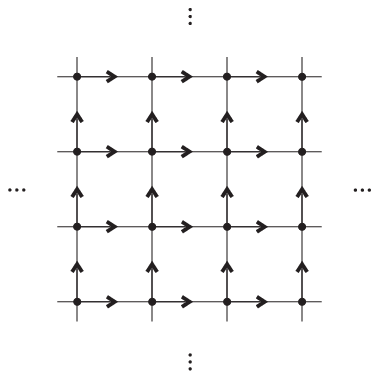}
 \caption{Infinite graph $\mathcal{H}^{\mathbb{Z}_2^{(3)}}_{\rho_1}$ for $U(2)$} \label{Fig-Graph_Z2^3-WG10}
\end{center}
\end{minipage}
\hfill
\begin{minipage}[t]{7.9cm}
\begin{center}
  \includegraphics[width=40mm]{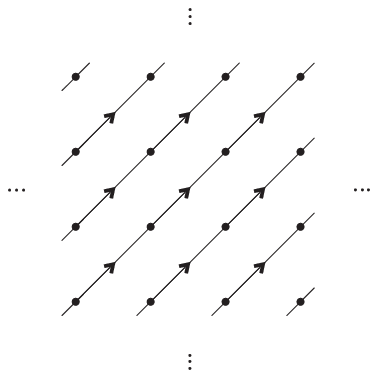}
 \caption{Infinite graph $\mathcal{H}^{\mathbb{Z}_2^{(3)}}_{\rho_2}$ for $U(2)$} \label{Fig-Graph_Z2^3-WG11}
\end{center}
\end{minipage}
\end{figure}

\begin{figure}[tb]
\begin{minipage}[t]{7cm}
\begin{center}
  \includegraphics[width=55mm]{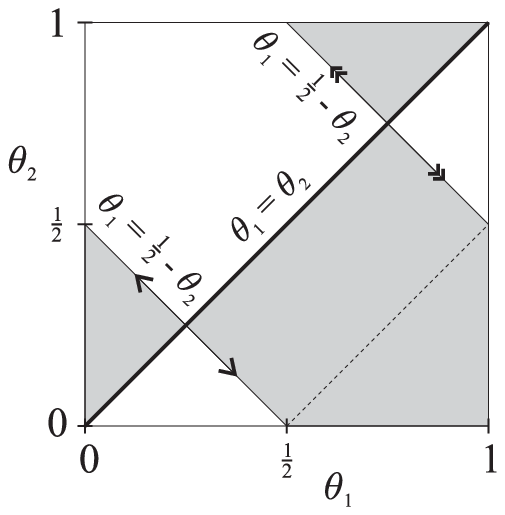}\\
 \caption{A fundamental domain $C$ \mbox{of $\mathbb{T}^2/\mathbb{Z}_2^{(3)}$} for $U(2)$.} \label{fig:DomainC-Z2^3}
\end{center}
\end{minipage}
\hfill
\begin{minipage}[t]{7.9cm}
\begin{center}
  \includegraphics[width=40mm]{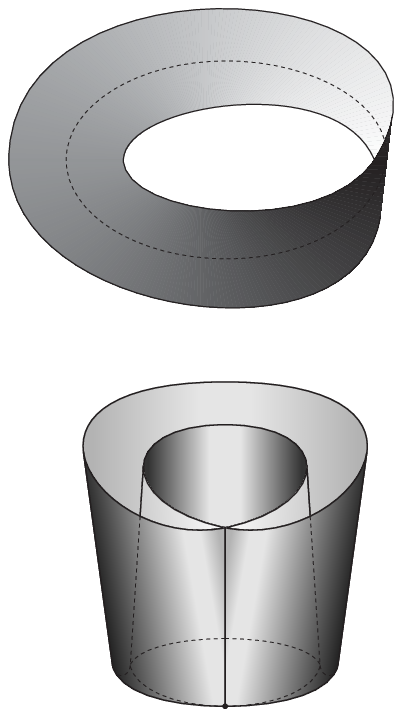}\\
 \caption{The surface $\mathbb{T}^2/\mathbb{Z}_2^{(3)}$ for $U(2)$.} \label{fig:DomainD-Z2^3}
\end{center}
\end{minipage}
\end{figure}

The subgroup $\mathbb{Z}_2^{(3)} \in GL(2,\mathbb{Z})$ is generated by the matrix $T_2'$ of Section \ref{sect:intro}.
It is the Weyl group of the connected compact Lie group $U(2)$.
We choose fundamental generators $\rho_1 = (1,0)$, $\rho_2 = (1,1)$.
The graphs $\mathcal{G}^{\mathbb{Z}_2^{(3)}}_{\rho}$, $\mathcal{H}^{\mathbb{Z}_2^{(3)}}_{\rho}$ for $\rho = \rho_1, \rho_2$ are illustrated in Figures \ref{Fig-Graph_Z2^3-G10}-\ref{Fig-Graph_Z2^3-WG11}.
Let
\begin{equation} \label{eqn:x,y-Z2^3}
x := x_{\rho_1} = \omega_1 + \omega_2, \qquad y := x_{\rho_2} = \omega_1\omega_2,
\end{equation}
and denote by $\Psi$ be the map $\Psi_{\rho_1,\rho_2}: (\omega_1,\omega_2) \mapsto (x,y)$.
A fundamental domain of $\mathbb{T}^2/\mathbb{Z}_2^{(3)}$ is illustrated in Figure \ref{fig:DomainC-Z2^3}, where the boundaries marked by arrows are identified. Then $\mathbb{T}^2/\mathbb{Z}_2^{(3)}$ is the M\"{o}bius strip, as in the top figure in Figure \ref{fig:DomainD-Z2^3}, where the dashed line $\theta_1=\theta_2+1/2$ in Figure \ref{fig:DomainC-Z2^3} is identified with the dashed line around the centre of the M\"{o}bius strip, and $\mathfrak{D} := \Psi(C)$ is an embedding of the M\"{o}bius strip in $\mathbb{C}^2$.
Under the change of variables $x = \omega_1 + \omega_2$, $y = \omega_1\omega_2$, the Jacobian is given by
$J_{\mathbb{Z}_2^{(3)}} = 4 \pi^2 \omega_1\omega_2(\omega_2 - \omega_1)$.
The Jacobian is complex-valued and vanishes in $\mathbb{T}^2$ only on the boundary $\theta_1 = \theta_2$ of the fundamental domain $C$.
Now $J_{\mathbb{Z}_2^{(3)}}^2 = 16\pi^4 \omega_1^2\omega_2^2(\omega_1^2-2\omega_1\omega_2+\omega_2^2)$ is invariant under the action of $\mathbb{Z}_2^{(3)}$ on $\mathbb{T}^2$, and we write $J_{\mathbb{Z}_2^{(3)}}^2$ in terms of the $\mathbb{Z}_2^3$-invariant elements $x$, $y$ as $J_{\mathbb{Z}_2^{(3)}}^2 = 16\pi^4 y^2(x^2-4y)$. It is easy to check that $|J_{\mathbb{Z}_2^{(3)}}| = 4 \pi^2 \sqrt{4-|x|^2}$.

\begin{Rem} \label{Rem:Z2^2-Z2^3}
Although the groups $\mathbb{Z}_2^{(2)}$ and $\mathbb{Z}_2^{(3)}$ are not conjugate in $GL(2,\mathbb{Z})$, they are however conjugate in $GL(2,\mathbb{R})$, where the conjugating matrix is
$H=\left( \begin{array}{cc} 1 & 1 \\ 1 & -1 \end{array} \right)$
which has inverse $H^{-1}=\left( \begin{array}{cc} 1/2 & 1/2 \\ 1/2 & -1/2 \end{array} \right) = \frac{1}{2} H$.
Thus $H$ intertwines $\mathbb{Z}_2^{(2)}$ and $\mathbb{Z}_2^{(3)}$ with $H\mathbb{Z}_2^{(2+i)} = \mathbb{Z}_2^{(3-i)}H$ for both $i=0,1$.
The geometric effect of the action of $H$ on $\mathbb{Z}^2$ is to reflect about the line $2y=x$ (or equivalently, to rotate the plane clockwise by $\pi/4$ and then reflect about the $x$-axis) and scale by $\sqrt{2}$. Applying this action twice simply has the effect of scaling $\mathbb{Z}^2$ by 2, as $H^2=2I$.

The origin of this relationship is the fact that the compact, connected Lie group $U(1) \times SU(2)$ (which has Weyl group $\mathbb{Z}_2^{(2)}$) is a double cover of the compact, non-connected Lie group $U(2)$ (which has Weyl group $\mathbb{Z}_2^{(3)}$). Thus not all irreducible representations of $U(1) \times SU(2)$ are irreducible representations of $U(2)$, but only the ``even'' ones, that is, those indexed by $\lambda = (\lambda_1,\lambda_2) \in \mathbb{Z}^2$ such that $\lambda_1 + \lambda_2 \equiv 0$ mod(2), that is, the irreducible representations $\mu$ of $U(2)$ yield the ``even'' irreducible representations $H\mu$ of $U(1) \times SU(2)$.

Re-drawing the M\"{o}bius strip in Figure \ref{fig:DomainD-Z2^3} by ``folding'' it about the dashed line around its centre we obtain the bottom figure in Figure \ref{fig:DomainD-Z2^3}, which should now be regarded as a surface in $\mathbb{R}^4$. It has a line of self-intersection when drawn in $\mathbb{R}^3$ (given by the vertical line down the centre of the surface), however in $\mathbb{R}^4$ the surface only intersects along this line at the point $\bullet$ at the bottom of the figure. The boundary of the surface is the top edge which is isomorphic to the circle (note that this circle does not self-intersect in $\mathbb{R}^4$ even though it appears to when drawn in $\mathbb{R}^3$).

For new variables $x_{H(0,1)} = x_{(1,-1)} = \omega_1\omega_2^{-1} + \omega_1^{-1}\omega_2$, $y_1=y\;(=x_{H(1,0)}=x_{(1,1)})$, in the joint spectrum $\mathfrak{D}_{(1,1),(1,-1)}$ the ``inner wall'' and ``outer wall'' in the figure at the bottom of Figure \ref{fig:DomainD-Z2^3} are identified, thus $\mathfrak{D}_{(1,1),(1,-1)}$ is the cylinder  $\mathbb{T} \times [-2,2]$, which is the joint spectrum $\mathfrak{D}^{\mathbb{Z}_2^{(2)}}$ for $\mathbb{Z}_2^{(2)}$ in Section \ref{sect:Z2^2}.
\end{Rem}

\subsection{Spectral measure for $\mathcal{H}^{\mathbb{Z}_2^{(3)}}_{\rho}$ for $U(2)$}

By integrating $|\Gamma| \, |J_{\mathbb{Z}_2^{(3)}}(x,y)^{-1}| = (4\pi^2 \sqrt{4-|x|^2})^{-1}$ over $y$, $x$ respectively we obtain the spectral measure $\nu_{\rho_1}$, $\nu_{\rho_2}$ respectively for ${}^{\mathbb{Z}_2^{(2)}} \hspace{-1mm} \Delta_{\rho_1}$, ${}^{\mathbb{Z}_2^{(2)}} \hspace{-1mm} \Delta_{\rho_2}$. Alternatively, we can determine these measures independently as follows.

We determine first the spectral measure for $\mathcal{H}^{\mathbb{Z}_2^{(3)}}_{\rho_1}$ over the spectrum $\chi_{\rho_1}(\mathbb{T}^2)$, which is uniquely determined by the moments $\varphi({}^{\mathbb{Z}_2^{(3)}} \hspace{-1mm} \Delta_{\rho_1}^m ({}^{\mathbb{Z}_2^{(3)}} \hspace{-1mm} \Delta_{\rho_1}^{\ast})^n)$ of ${}^{\mathbb{Z}_2^{(3)}} \hspace{-1mm} \Delta_{\rho_1}$. The spectrum $\sigma({}^{\mathbb{Z}_2^{(3)}} \hspace{-1mm} \Delta_{\rho_1}) = \chi_{\rho_1}(\mathbb{T}^2)$ of ${}^{\mathbb{Z}_2^{(3)}} \hspace{-1mm} \Delta_{\rho_1}$ is given by the disc $2\mathbb{D} = \{ z \in \mathbb{C} | \, |z| \leq 2 \}$ with radius 2.

\begin{figure}[htb]
\begin{center}
  \includegraphics[width=70mm]{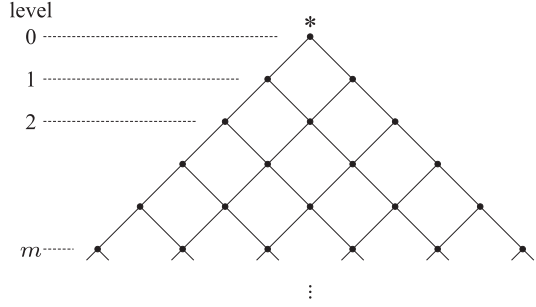}
 \caption{The Bratteli diagram for Dynkin diagram $A_{\infty,\infty}$} \label{Fig-Bratteli_diag-Ainftyinfty}
\end{center}
\end{figure}

The $m,n^{\mathrm{th}}$ moment $\varphi({}^{\mathbb{Z}_2^{(3)}} \hspace{-1mm} \Delta_{\rho_1}^m ({}^{\mathbb{Z}_2^{(3)}} \hspace{-1mm} \Delta_{\rho_1}^{\ast})^n)$ counts the number of paths of length $m+n$ on $\mathcal{H}^{\mathbb{Z}_2^{(3)}}_{\rho_1}$ and its opposite graph $(\mathcal{H}^{\mathbb{Z}_2^{(3)}}_{\rho_1})^{\mathrm{op}}$ (that is, the graph $\mathcal{H}^{\mathbb{Z}_2^{(3)}}_{\rho_1}$ with the orientation of all edges reversed) which start at any choice of distinguished vertex $\ast$, where the first $m$ edges are on $\mathcal{H}^{\mathbb{Z}_2^{(3)}}_{\rho_1}$ and the next $n$ edges are on its opposite graph. Due to the orientation of the edges of $\mathcal{H}^{\mathbb{Z}_2^{(3)}}_{\rho_1}$, the paths of length $m$ on $\mathcal{H}^{\mathbb{Z}_2^{(3)}}_{\rho_1}$ are given by the number of paths from $\ast$ to the $m^{\mathrm{th}}$ level of the graph in Figure \ref{Fig-Bratteli_diag-Ainftyinfty}.

It is thus easy to see that $\varphi({}^{\mathbb{Z}_2^{(3)}} \hspace{-1mm} \Delta_{\rho_1}^m ({}^{\mathbb{Z}_2^{(3)}} \hspace{-1mm} \Delta_{\rho_1}^{\ast})^n)$ will be zero unless $m=n$. When $m=n$, $\varphi({}^{\mathbb{Z}_2^{(3)}} \hspace{-1mm} \Delta_{\rho_1}^m ({}^{\mathbb{Z}_2^{(3)}} \hspace{-1mm} \Delta_{\rho_1}^{\ast})^n)$ counts all pairs of paths which start at $\ast$ and which both end at the same vertex at the $m^{\mathrm{th}}$ level of the graph in Figure \ref{Fig-Bratteli_diag-Ainftyinfty}, that is, $\varphi({}^{\mathbb{Z}_2^{(3)}} \hspace{-1mm} \Delta_{\rho_1}^m ({}^{\mathbb{Z}_2^{(3)}} \hspace{-1mm} \Delta_{\rho_1}^{\ast})^n)$ is the dimension of the finite-dimensional algebra given by the $m^{\mathrm{th}}$ level of the Bratteli diagram in Figure \ref{Fig-Bratteli_diag-Ainftyinfty}. These dimensions are given by the central binomial coefficient $C^{2m}_m$ (c.f. \cite[$\S$2.1]{evans/pugh:2009v}), thus $\varphi({}^{\mathbb{Z}_2^{(3)}}\Delta_{\rho_1}^m({}^{\mathbb{Z}_2^{(3)}}\Delta_{\rho_1}^{\ast})^n) = \delta_{m,n} C^{2m}_m$. The spectral measure (over [-2,2]) for the Dynkin diagram $A_{\infty,\infty}$ is given by $(\pi\sqrt{4-r^2})^{-1} \mathrm{d}r$.
Thus the spectral measure (over $2\mathbb{D}$) for $\mathcal{H}^{\mathbb{Z}_2^{(3)}}_{\rho_1}$ is given by $(\pi^2\sqrt{4-|x|^2})^{-1} \mathrm{d}x$, since
\begin{align*}
\int_{2\mathbb{D}} \frac{x^m \overline{x}^n}{\pi^2\sqrt{4-|x|^2}} \, \mathrm{d}x & := \int_0^2 \int_0^{2\pi} \frac{r^{m+n} e^{i\theta(m-n)}}{\pi^2\sqrt{4-r^2}} \, \mathrm{d}\theta \mathrm{d}r = \int_0^2 \frac{r^{m+n}}{\pi^2\sqrt{4-r^2}} \, \mathrm{d}r \; \int_0^{2\pi} e^{i\theta(m-n)} \mathrm{d}\theta \\
& =  \frac{1}{2} \int_{-2}^2 \frac{r^{2m}}{\pi^2\sqrt{4-r^2}} \, \mathrm{d}r \; \delta_{m,n} 2\pi = C^{2m}_m \delta_{m,n},
\end{align*}
where the penultimate equality follows since the product is zero for $m \neq n$, whilst for $m=n$ we have the integral of $r^{2m}/\sqrt{4-r^2}$ which is an even function.

The graph $\mathcal{H}^{\mathbb{Z}_2^{(3)}}_{\rho_2}$, illustrated in Figure \ref{Fig-Graph_Z2^3-WG11}, is given by an infinite number of copies of the representation graph $\mathcal{G}^0_{\rho}$ for $\mathbb{T}$, and in particular the connected component of the distinguished vertex $(0,0)$ of $\mathcal{H}^{\mathbb{Z}_2^{(3)}}_{\rho_2}$ is the representation graph for $\mathbb{T}$. Thus we have

\begin{Thm} \label{thm:measure-rho-Z2^3-H}
The spectral measure $\nu_{\rho_1}$ (over $2\mathbb{D}$) for the graph $\mathcal{H}^{\mathbb{Z}_2^{(3)}}_{\rho_1}$ for $U(2)$ is given by
$$\mathrm{d}\nu_{\rho_1}(x) = \frac{1}{\pi^2 \sqrt{4-|x|^2}} \, \mathrm{d}x, \qquad x \in 2\mathbb{D},$$
whilst the spectral measure $\nu_{\rho_2}$ (over $\mathbb{T}$) for the graph $\mathcal{H}^{\mathbb{Z}_2^{(3)}}_{\rho_2}$ for $U(2)$ is given by
$\mathrm{d}\nu_{\rho_2}(y) = (2\pi)^{-1} \, \mathrm{d}y$, $y \in \mathbb{T}$.
\end{Thm}

\subsection{Spectral measure for $\mathcal{G}^{\mathbb{Z}_2^{(3)}}_{\rho}$ for $U(2)$}

\begin{figure}[htb]
\begin{center}
  \includegraphics[width=40mm]{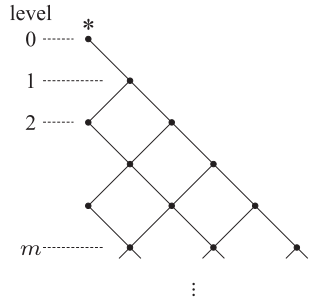}
 \caption{The Bratteli diagram for Dynkin diagram $A_{\infty}$} \label{Fig-Bratteli_diag-Ainfty}
\end{center}
\end{figure}

We turn now to the graphs $\mathcal{G}^{\mathbb{Z}_2^{(3)}}_{\rho_j}$, $j=1,2$. The adjacency matrix $\Delta_{\rho_1}$ for the graph $\mathcal{G}^{\mathbb{Z}_2^{(3)}}_{\rho_1}$ is normal, and thus its spectral measure (over $\chi_{\rho_1}(\mathbb{T}^2)$) is uniquely determined by its moments $\varphi(\Delta_{\rho_1}^m(\Delta_{\rho_1}^{\ast})^n)$. Its spectrum $\sigma(\Delta_{\rho_1}) = \chi_{\rho_1}(\mathbb{T}^2)$ is again given by the disc $2\mathbb{D}$ of radius 2. We define a state $\varphi$ by $\varphi( \, \cdot \, ) = \langle \cdot \Omega, \Omega \rangle$, where $\Omega$ is vector in $\ell^2(\mathcal{G}^{\mathbb{Z}_2^{(3)}}_{\rho_1})$ corresponding to a distinguished vertex $\ast$, which is chosen to be one of the vertices which is only the source (or range) of one edge. Thus the $m,n^{\mathrm{th}}$ moment $\varphi(\Delta_{\rho_1}^m(\Delta_{\rho_1}^{\ast})^n)$ counts the number of paths of length $m+n$ on $\mathcal{G}^{\mathbb{Z}_2^{(3)}}_{\rho_1}$ and its opposite graph $(\mathcal{G}^{\mathbb{Z}_2^{(3)}}_{\rho_1})^{\mathrm{op}}$, where the first $m$ edges are on $\mathcal{G}^{\mathbb{Z}_2^{(3)}}_{\rho_1}$ and the next $n$ edges are on its opposite graph. Due to the orientation of the edges of $\mathcal{G}^{\mathbb{Z}_2^{(3)}}_{\rho_1}$, the paths of length $m$ on $\mathcal{G}^{\mathbb{Z}_2^{(3)}}_{\rho_1}$ are given by the number of paths from $\ast$ to the $m^{\mathrm{th}}$ level of the graph in Figure \ref{Fig-Bratteli_diag-Ainfty}.

It is thus easy to see again that $\varphi(\Delta_{\rho_1}^m(\Delta_{\rho_1}^{\ast})^n)$ will be zero unless $m=n$. When $m=n$, $\varphi(\Delta_{\rho_1}^m(\Delta_{\rho_1}^{\ast})^m)$ counts all pairs of paths which start at $\ast$ and which both end at the same vertex at the $m^{\mathrm{th}}$ level of the graph in Figure \ref{Fig-Bratteli_diag-Ainfty}, that is, $\varphi(\Delta_{\rho_1}^m(\Delta_{\rho_1}^{\ast})^m)$ is the dimension of the finite-dimensional algebra given by the $m^{\mathrm{th}}$ level of the Bratteli diagram in Figure \ref{Fig-Bratteli_diag-Ainfty}. These dimensions are well known to be given by the Catalan numbers $c_m = C^{2m}_m/(m+1)$ \cite[Aside 5.1.1]{jones:1983}, thus $\varphi(\Delta_{\rho_1}^m(\Delta_{\rho_1}^{\ast})^n) = \delta_{m,n} c_m$. The spectral measure (over [-2,2]) for the Dynkin diagram $A_{\infty}$ is the semi-circle measure $\sqrt{4-r^2} \mathrm{d}r/2\pi$. Thus the spectral measure (over $2\mathbb{D}$) for $\mathcal{G}^{\mathbb{Z}_2^{(3)}}_{\rho_1}$ is given by $\sqrt{4-|x|^2} \mathrm{d}x/2\pi^2$, since
\begin{align*}
\int_{2\mathbb{D}} x^m \overline{x}^n \sqrt{4-|x|^2} \, \mathrm{d}x & := \int_0^2 \int_0^{2\pi} r^{m+n} e^{i\theta(m-n)} \sqrt{4-r^2} \, \mathrm{d}\theta \mathrm{d}r \\
& = \int_0^2 r^{m+n} \sqrt{4-r^2} \, \mathrm{d}r \; \int_0^{2\pi} e^{i\theta(m-n)} \mathrm{d}\theta \\
& = \frac{1}{2} \int_{-2}^2 r^{2m} \sqrt{4-r^2} \, \mathrm{d}r \; \delta_{m,n} 2\pi = 2\pi^2 c_m \delta_{m,n}.
\end{align*}

The graph $\mathcal{G}^{\mathbb{Z}_2^{(3)}}_{\rho_2}$ is again given by an infinite number of copies of the representation graph $\mathcal{G}^0_{\rho}$ for $\mathbb{T}$, thus we have

\begin{Thm} \label{thm:measure-rho-Z2^3-G}
The spectral measure $\nu_{\rho_1}$ (over $2\mathbb{D}$) for the graph $\mathcal{G}^{\mathbb{Z}_2^{(3)}}_{\rho_1}$ for $U(2)$ is given by
$$\mathrm{d}\nu_{\rho_1}(x) = \frac{1}{2\pi^2} \sqrt{4-|x|^2} \, \mathrm{d}x, \qquad x \in 2\mathbb{D},$$
whilst the spectral measure $\nu_{\rho_2}$ (over $\mathbb{T}$) for the graph $\mathcal{G}^{\mathbb{Z}_2^{(3)}}_{\rho_2}$ for $U(2)$ is given by
$\mathrm{d}\nu_{\rho_2}(y) = (2\pi)^{-1} \, \mathrm{d}y$, $y \in \mathbb{T}$.
\end{Thm}

\section{$D_4^{(1)}$: $SU(2) \times SU(2)$}

\begin{figure}[tb]
\begin{minipage}[t]{7cm}
\begin{center}
  \includegraphics[width=25mm]{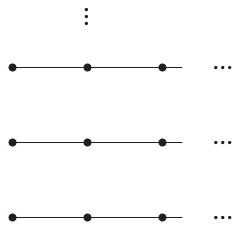}
 \caption{Infinite graph $\mathcal{G}^{D_4^{(1)}}_{\rho_1}$ for \mbox{$SU(2) \times SU(2)$}} \label{Fig-Graph_D4^1-G10}
\end{center}
\end{minipage}
\hfill
\begin{minipage}[t]{7cm}
\begin{center}
  \includegraphics[width=25mm]{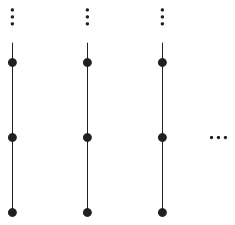}
 \caption{Infinite graph $\mathcal{G}^{D_4^{(1)}}_{\rho_2}$ for \mbox{$SU(2) \times SU(2)$}} \label{Fig-Graph_D4^1-G01}
\end{center}
\end{minipage}
\end{figure}

\begin{figure}[tb]
\begin{minipage}[t]{7cm}
\begin{center}
  \includegraphics[width=40mm]{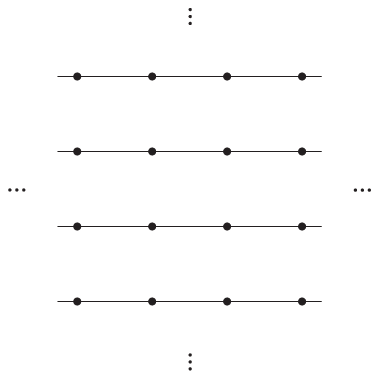}
 \caption{Infinite graph $\mathcal{H}^{D_4^{(1)}}_{\rho_1}$ for \mbox{$SU(2) \times SU(2)$}} \label{Fig-Graph_D4^1-WG10}
\end{center}
\end{minipage}
\hfill
\begin{minipage}[t]{7cm}
\begin{center}
  \includegraphics[width=40mm]{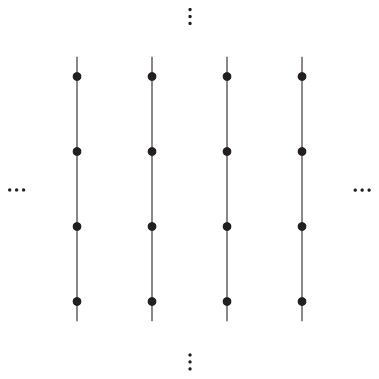}
 \caption{Infinite graph $\mathcal{H}^{D_4^{(1)}}_{\rho_2}$ for \mbox{$SU(2) \times SU(2)$}} \label{Fig-Graph_D4^1-WG01}
\end{center}
\end{minipage}
\end{figure}

\begin{figure}[tb]
\begin{minipage}[t]{7cm}
\begin{center}
  \includegraphics[width=55mm]{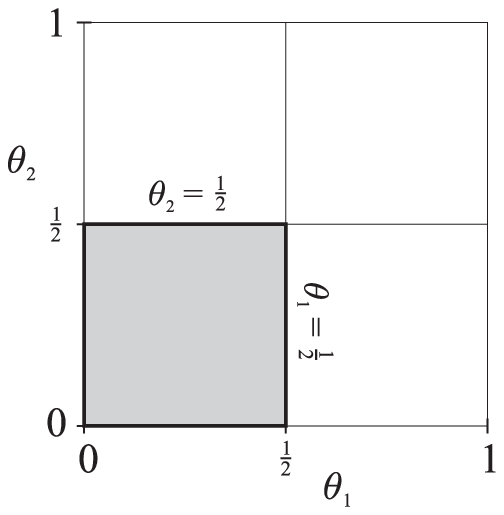}\\
 \caption{A fundamental domain $C$ \mbox{of $\mathbb{T}^2/D_4^{(1)}$} for $SU(2) \times SU(2)$.} \label{fig:DomainC-D4^1}
\end{center}
\end{minipage}
\hfill
\begin{minipage}[t]{7cm}
\begin{center}
  \includegraphics[width=55mm]{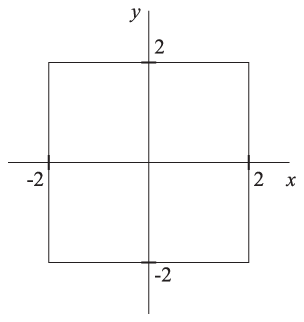}\\
 \caption{The domain $\mathfrak{D} = \Psi(C)$ \mbox{for $SU(2) \times SU(2)$}.} \label{fig:DomainD-D4^1}
\end{center}
\end{minipage}
\end{figure}

The subgroup $D_4^{(1)} \in GL(2,\mathbb{Z})$ is generated by $-I$ and the matrix $T_2$ of Section \ref{sect:intro}. It is the Weyl group for the connected, simply-connected, semisimple compact Lie group $SU(2) \times SU(2)$ (c.f. Section \ref{sect:Z2^2}).
We choose fundamental generators $\rho_1 = (1,0)$, $\rho_2 = (0,1)$.
The graphs $\mathcal{G}^{D_4^{(1)}}_{\rho}$, $\mathcal{H}^{D_4^{(1)}}_{\rho}$ for $\rho = \rho_1, \rho_2$ are illustrated in Figures \ref{Fig-Graph_D4^1-G10}-\ref{Fig-Graph_D4^1-WG01}.
Let
\begin{equation} \label{eqn:x,y-D4^1}
x := x_{\rho_1} = \omega_1 + \omega_1^{-1} = \cos(2\pi\theta_1), \qquad y := x_{\rho_2} = \omega_2 + \omega_2^{-1} = \cos(2\pi\theta_2),
\end{equation}
and denote by $\Psi$ be the map $\Psi_{\rho_1,\rho_2}: (\omega_1,\omega_2) \mapsto (x,y)$.
A fundamental domain $C$ of $\mathbb{T}^2/D_4^{(1)}$ is illustrated in Figure \ref{fig:DomainC-D4^1}.
Under the change of variables $x = \omega_1 + \omega_1^{-1}$, $y = \omega_2 + \omega_2^{-1}$, the Jacobian is given by
$J_{D_4^{(1)}} = 16 \pi^2 \sin(2 \pi \theta_1) \sin(2 \pi \theta_2) = 4 \pi^2 (\omega_1\overline{\omega_2} + \overline{\omega_1}\omega_2 - \omega_1\omega_2 - \overline{\omega_1}\overline{\omega_2})$.
The Jacobian is real and vanishes in $\mathbb{T}^2$ only on the boundaries of the images of the fundamental domain $C$ under $D_4^{(1)}$.
Again, $J_{D_4^{(1)}}^2$ is invariant under the action of $D_4^{(1)}$ on $\mathbb{T}^2$, and $J_{D_4^{(1)}}^2$ can be written in terms of the $D_4^{(1)}$-invariant elements $x$, $y$ as $J_{D_4^{(1)}}^2 = 16 \pi^4 (4-x^2)(4-y^2)$, where $4-p^2 \geq 0$ for all $p \in [-2,2]$. Thus we write $J_{D_4^{(1)}}$ in terms of $x$, $y$ as
$J_{D_4^{(1)}} = 4 \pi^2 \sqrt{(4-x^2)(4-y^2)}$.

\begin{Rem} \label{Rem:Z2^2<D4^1}
The group $\mathbb{Z}_2^{(2)}$ is a normal subgroup of $D_4^{(1)}$. The fundamental domain $C$ for $D_4^{(1)} = \mathbb{Z}_2^{(2)} \rtimes \mathbb{Z}_2^{(1)}$ in Figure \ref{fig:DomainC-D4^1} is obtained from the fundamental domain for $\mathbb{Z}_2^{(2)}$ in Figure \ref{fig:DomainC-Z2^2} by imposing one extra symmetry which comes from the additional $\mathbb{Z}_2^{(1)}$ action. If we denote $x,y$ for $\mathbb{Z}_2^{(2)}$ given in (\ref{eqn:x,y-Z2^2}) by $x_1,y_1$, and $x,y$ for $D_4^{(1)}$ given in (\ref{eqn:x,y-D4^1}) by $x_2,y_2$, then we see that $y_1=y_2$ whilst $x_2 = 2\mathrm{Re}(x_1)$. Then there is a homomorphism $\xi:\mathfrak{D}_{\mathbb{Z}_2^{(2)}} \rightarrow \mathfrak{D}_{D_4^{(1)}}$ from the domain $\mathfrak{D}_{\mathbb{Z}_2^{(2)}}$ of $\mathbb{Z}_2^{(2)}$ to the domain $\mathfrak{D}_{D_4^{(1)}}$ of $D_4^{(1)}$ such that $\xi(y_1)=y_2$ and $\xi(x_1)=2\mathrm{Re}(x_1)$, which maps the boundary of $\mathfrak{D}_{\mathbb{Z}_2^{(2)}}$ to part of the boundary of $\mathfrak{D}_{D_4^{(1)}}$. The rest of the boundary of $\mathfrak{D}_{D_4^{(1)}}$ is given by $\xi(v)$, where $v=(x,y) \in \mathfrak{D}_{\mathbb{Z}_2^{(2)}}$ such that $x \in \mathbb{R}$, i.e. $v = \Phi(t)$ for $t \in \mathbb{T}^2$ such that $t$ is fixed under the additional $\mathbb{Z}_2$ action in $D_4^{(1)}$.
Moreover, we see that $J_{D_4^{(1)}} = 2\mathrm{Re}(J_{\mathbb{Z}_2^{(2)}}) = \xi(J_{\mathbb{Z}_2^{(2)}})$.
\end{Rem}

\begin{Rem} \label{Rem:Z2^1=2.D4^1}
As noted in the introduction, the subgroup $\mathbb{Z}_2^{(1)}$ of $SL(2,\mathbb{Z})$ is also a normal subgroup of $D_4^{(1)} \cong \mathbb{Z}_2^{(1)} \rtimes \mathbb{Z}_2^{(2)}$. Since $P_+$ is not uniquely defined, the definition of formal character given in (\ref{eqn:character=S/S}) does not make sense. However, for a fixed choice $P_+$ of fundamental domain of $\mathbb{T}^2/\mathbb{Z}_2^{(1)}$ one can define a formal character for all $\lambda \in P_+$ simply by $\chi_{\lambda}(\omega_1,\omega_2) = S_{\lambda}(\theta_1,\theta_2)$ ($\, = C_{\lambda}(\theta_1,\theta_2)$ since both elements of $\mathbb{Z}_2^{(1)}$ have determinant one), where $\omega_j = e^{2\pi i \theta_j} \in \mathbb{T}$, $j=1,2$. Then one can construct representation graphs $\mathcal{G}_{\rho_j}^{\mathbb{Z}_2^{(1)}}$, $\mathcal{H}_{\rho_j}^{\mathbb{Z}_2^{(1)}}$ for $j=1,2$, where $\rho_1 = (1,0)$, $\rho_2 = (0,1)$ play the role of the fundamental generators (although one also needs the fact that $\chi_{(-\lambda_1,\lambda_2)} (\omega_1,\omega_2) = \chi_{(\lambda_1,\lambda_2)} (\overline{\omega_1},\omega_2)$ and/or $\chi_{(\lambda_1,-\lambda_2)} (\omega_1,\omega_2) = \chi_{(\lambda_1,\lambda_2)} (\omega_1,\overline{\omega_2})$ in order to generate all characters). Then $x:= x_{\rho_1} = \omega_1 + \overline{\omega_1}$ and $y:= x_{\rho_2} = \omega_2 + \overline{\omega_2}$, and the map $\Psi: (\omega_1,\omega_2) \mapsto (x,y)$ is a map from $\mathbb{T}^2$ to the joint spectrum $\mathfrak{D}^{D_4^{(1)}} = [-2,2] \times [-2,2]$ for $D_4^{(1)}$ described above. In this case, $\Psi$ is a two-to-one map from $P_+$ to $\mathfrak{D}^{D_4^{(1)}}$.
The Jacobian for the change of variables $\omega_1,\omega_2 \rightarrow x,y$ is again given by $J_{D_4^{(1)}}$, and the spectral measures (over $\mathfrak{D}^{D_4^{(1)}}$) for $\rho_1$, $\rho_2$ are given by twice the spectral measures (over $\mathfrak{D}^{D_4^{(1)}}$) for $D_4^{(1)}$ in (\ref{eqn:measure-rho-D4^1-H}), (\ref{eqn:measure-rho-D4^1-G}) below (the factor of two comes from the fact that $\Psi$ is now a two-to-one map from $P_+$ to $\mathfrak{D}^{D_4^{(1)}}$).
\end{Rem}

\subsection{Spectral measure for $\mathcal{H}^{D_4^{(1)}}_{\rho}$ for $SU(2) \times SU(2)$}

Here $\mathfrak{D} = [-2,2] \times [-2,2]$ is the joint spectrum $\sigma({}^{D_4^{(1)}} \hspace{-1mm} \Delta_{\rho_1},{}^{D_4^{(1)}} \hspace{-1mm} \Delta_{\rho_2})$ of the commuting self-adjoint operators ${}^{D_4^{(1)}} \hspace{-1mm} \Delta_{\rho_1}$, ${}^{D_4^{(1)}} \hspace{-1mm} \Delta_{\rho_2}$.
Then by integrating $|\Gamma| \, J_{D_4^{(1)}}(x,y)^{-1}$ over $y$ we obtain the spectral measure $\nu_{\rho_1}$ for ${}^{D_4^{(1)}} \hspace{-1mm} \Delta_{\rho_1}$. Integrating $|\Gamma| \, J_{D_4^{(1)}}(x,y)^{-1}$ over $x$ gives the same result, and thus $\nu_{\rho_1} = \nu_{\rho_2}$.
Since $\int_{-2}^2 \sqrt{4-y^2}^{-1} \mathrm{d}y = \pi$, the spectral measures $\nu_{\rho_j}$ (over $\chi_{\rho_j}(\mathbb{T}^2) = [-2,2]$) for the graphs $\mathcal{H}^{D_4^{(1)}}_{\rho_j}$, $j=1,2$, for $SU(2) \times SU(2)$ are both given by
\begin{equation} \label{eqn:measure-rho-D4^1-H}
\mathrm{d}\nu_{\rho_1}(x) = \mathrm{d}\nu_{\rho_2}(x) = (\pi\sqrt{4-x^2})^{-1} \, \mathrm{d}x, \qquad x \in [-2,2]
\end{equation}

\subsection{Spectral measure for $\mathcal{G}^{D_4^{(1)}}_{\rho}$ for $SU(2) \times SU(2)$}

The graphs $\mathcal{G}^{D_4^{(1)}}_{\rho_2}$ are both given by an infinite number of copies of the representation graph $\mathcal{G}^{\langle -1 \rangle}_{\rho}$ for $SU(2)$, and in particular the connected component of the distinguished vertex $\ast$, the vertex with lowest Perron-Frobenius weight which in this case is the apex vertex, i.e. the vertex in the bottom left corner in Figures \ref{Fig-Graph_D4^1-G10} and \ref{Fig-Graph_D4^1-G01}. Thus the spectral measures $\nu_{\rho_j}$ (over $\chi_{\rho_j}(\mathbb{T}^2) = [-2,2]$) for the graphs $\mathcal{G}^{D_4^{(1)}}_{\rho_j}$, $j=1,2$, for $SU(2) \times SU(2)$ are given by
\begin{equation} \label{eqn:measure-rho-D4^1-G}
\mathrm{d}\nu_{\rho_1}(x) = \mathrm{d}\nu_{\rho_2}(x) = \frac{1}{2\pi}\sqrt{4-x^2} \, \mathrm{d}x, \qquad x \in [-2,2]
\end{equation}

\section{$D_4^{(2)}$: $SO(4)$} \label{sect:D4^2}

\begin{figure}[tb]
\begin{minipage}[t]{7.9cm}
\begin{center}
  \includegraphics[width=35mm]{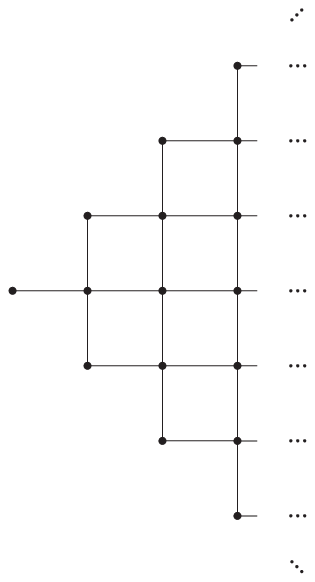}
 \caption{Infinite graph $\mathcal{G}^{D_4^{(2)}}_{\rho_1}$ for $SO(4)$} \label{Fig-Graph_D4^2-G10}
\end{center}
\end{minipage}
\hfill
\begin{minipage}[t]{7.9cm}
\begin{center}
  \includegraphics[width=35mm]{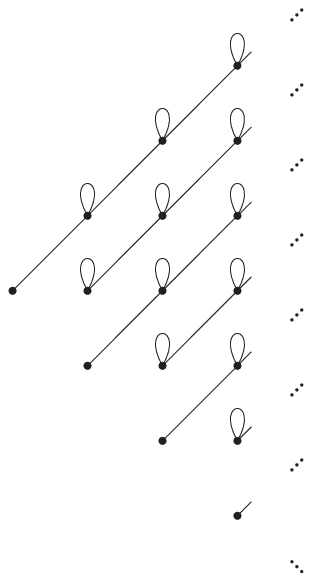}
 \caption{Infinite graph $\mathcal{G}^{D_4^{(2)}}_{\rho_2}$ for $SO(4)$} \label{Fig-Graph_D4^2-G11}
\end{center}
\end{minipage}
\end{figure}

\begin{figure}[tb]
\begin{minipage}[t]{7.9cm}
\begin{center}
  \includegraphics[width=40mm]{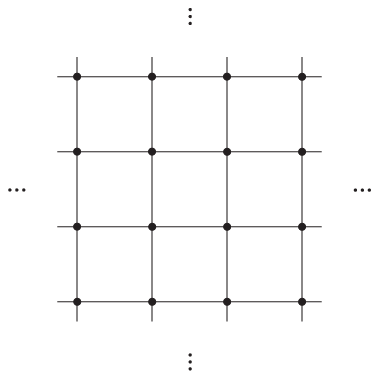}
 \caption{Infinite graph $\mathcal{H}^{D_4^{(2)}}_{\rho_1}$ for $SO(4)$} \label{Fig-Graph_D4^2-WG10}
\end{center}
\end{minipage}
\hfill
\begin{minipage}[t]{7.9cm}
\begin{center}
  \includegraphics[width=40mm]{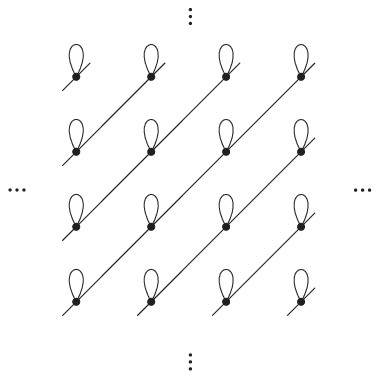}
 \caption{Infinite graph $\mathcal{H}^{D_4^{(2)}}_{\rho_2}$ for $SO(4)$} \label{Fig-Graph_D4^2-WG11}
\end{center}
\end{minipage}
\end{figure}

\begin{figure}[tb]
\begin{minipage}[t]{7cm}
\begin{center}
  \includegraphics[width=55mm]{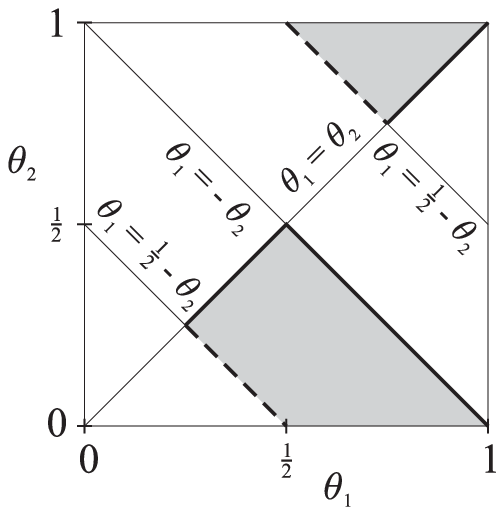}\\
 \caption{A fundamental domain $C$ \mbox{of $\mathbb{T}^2/D_4^{(2)}$} for $SO(4)$.} \label{fig:DomainC-D4^2}
\end{center}
\end{minipage}
\hfill
\begin{minipage}[t]{7cm}
\begin{center}
  \includegraphics[width=70mm]{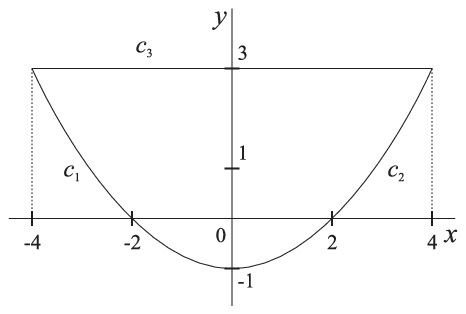}\\
 \caption{The domain $\mathfrak{D} = \Psi(C)$ \mbox{for $SO(4)$}.} \label{fig:DomainD-D4^2}
\end{center}
\end{minipage}
\end{figure}

The subgroup $D_4^{(2)} \in GL(2,\mathbb{Z})$ is generated by $-I$ and the matrix $T_2'$ of Section \ref{sect:intro}.
It is the Weyl group of the connected compact Lie group $SO(4)$.
We choose fundamental generators $\rho_1 = (1,0)$, $\rho_2 = (1,1)$.
The graphs $\mathcal{G}^{D_4^{(2)}}_{\rho}$, $\mathcal{H}^{D_4^{(2)}}_{\rho}$ for $\rho = \rho_1, \rho_2$ are illustrated in Figures \ref{Fig-Graph_D4^2-G10}-\ref{Fig-Graph_D4^2-WG11}.
Let
\begin{eqnarray}
x & := & x_{\rho_1} \;\; = \;\; \omega_1 + \omega_1^{-1} + \omega_2 + \omega_2^{-1} \;\; = \;\; 2\cos(2\pi \theta_1) + 2\cos(2\pi \theta_2), \label{eqn:x-D4^2} \\
y & := & x_{\rho_2} \;\; = \;\; 1 + \omega_1\omega_2 + \omega_1^{-1}\omega_2^{-1} \;\; = \;\; 1 + 2\cos(2\pi(\theta_1+\theta_2)), \label{eqn:y-D4^2}
\end{eqnarray}
and denote by $\Psi$ be the map $\Psi_{\rho_1,\rho_2}: (\omega_1,\omega_2) \mapsto (x,y)$.
A fundamental domain $C$ of $\mathbb{T}^2/D_4^{(2)}$ is illustrated in Figure \ref{fig:DomainC-D4^2}.
Then $\mathfrak{D} = \Psi(C)$ is illustrated in Figure \ref{fig:DomainD-D4^2}.
Now $(\theta_1,1/2-\theta_1)$ maps to $\Psi(\omega_1,-\overline{\omega_1}) = (0,-1)$ for all $\theta_1 \in \mathbb{T}$, the dashed line $\theta_2 = 1/2 - \theta_1$ in Figure \ref{fig:DomainC-D4^2} contracts to the single point $(0,-1)$ in $\mathfrak{D}$.
The boundary of $C$ given by $\theta_1 = \theta_2$, where $\theta_1 \in [1/4,1/2]$, yields the curve $c_1$ given by the parametric equations
$x = 4\cos(2\pi\theta_1)$, $y = 2\cos(4\pi\theta_1) + 1 = 4\cos^2(2\pi\theta_1) -1$. Similarly, the boundary of $C$ given by $\theta_1 = \theta_2$, where $\theta_1 \in [3/4]$, yields the curve $c_2$ given by the same parametric equations.
The boundary of $C$ given by $\theta_1 = -\theta_2$ yields the curve $c_3$ given by the parametric equations
$x = 4\cos(2\pi\theta_1)$, $y = 3$, where $\theta_1 \in [1/2,1]$.
As functions of $x \in [-4,4]$, the boundaries $c_1$, $c_2$ of $\mathfrak{D}$ are thus given by $y = x^2/4-1$ whilst $c_3$ is given by $y = 3$.
As functions of $y$, the boundaries $c_1$, $c_2$ are given by
$x = -2\sqrt{y+1}$, $x = 2\sqrt{y+1}$ respectively.
Under the change of variables (\ref{eqn:x-D4^2}), (\ref{eqn:y-D4^2}), the Jacobian is given by
$J_{D_4^{(2)}} = 16 \pi^2 (\cos(2 \pi \theta_2) + \cos(2 \pi (\theta_1+2\theta_2)) - \cos(2 \pi \theta_1) - \cos(2 \pi (2\theta_1+\theta_2)))$.
The Jacobian is real and vanishes in $\mathbb{T}^2$ only on the boundaries of the images of the fundamental domain $C$ under $D_4^{(2)}$.
Again, $J_{D_4^{(2)}}^2$ is invariant under the action of $D_4^{(2)}$ on $\mathbb{T}^2$, and can be written in terms of the $D_4^{(2)}$-invariant elements $x$, $y$ as $J_{D_4^{(2)}}^2 = 16 \pi^4 (y+1)(3-y)(4y-x^2+4)$, which is non-negative since $J_{D_4^{(2)}}$ is real. Thus we can write $J_{D_4^{(2)}}$ in terms of $x$, $y$ as
$J_{D_4^{(2)}} = 4 \pi^2 \sqrt{(y+1)(3-y)(4y-x^2+4)}$.

\begin{figure}[tb]
\begin{minipage}[t]{7.9cm}
\begin{center}
  \includegraphics[width=55mm]{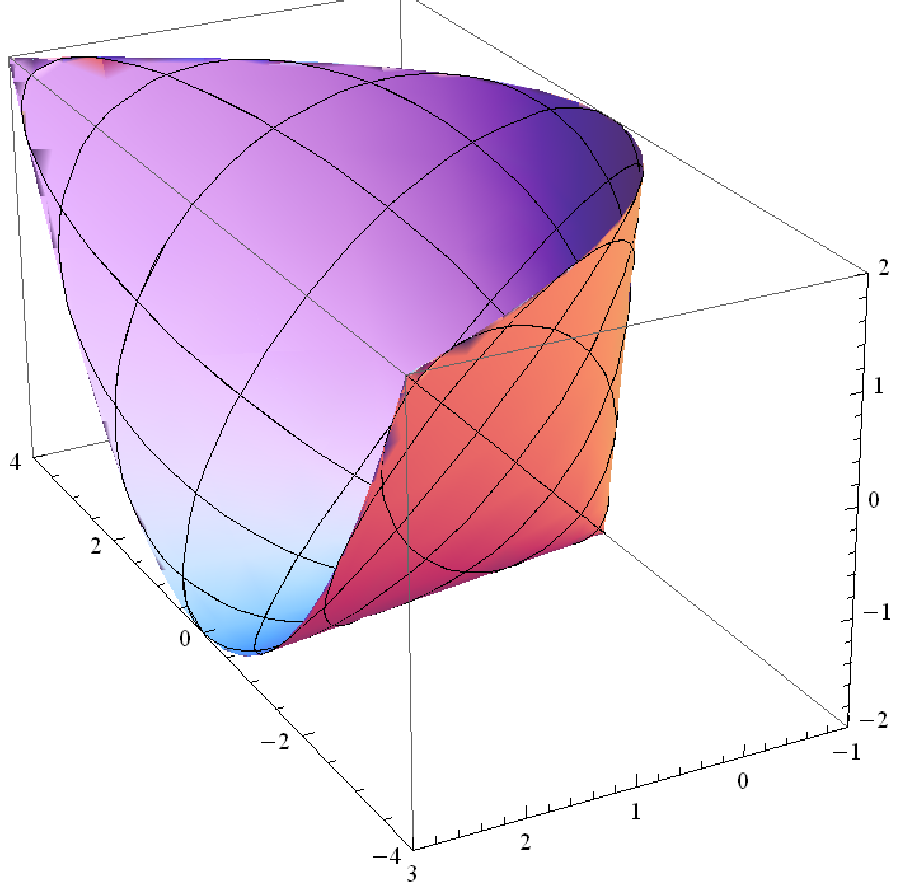}\\
 \caption{The domain $\widetilde{\mathfrak{D}} = \widetilde{\Phi}(C)$.} \label{fig:widetildeD-D4^2}
\end{center}
\end{minipage}
\hfill
\begin{minipage}[t]{7.9cm}
\begin{center}
  \includegraphics[width=55mm]{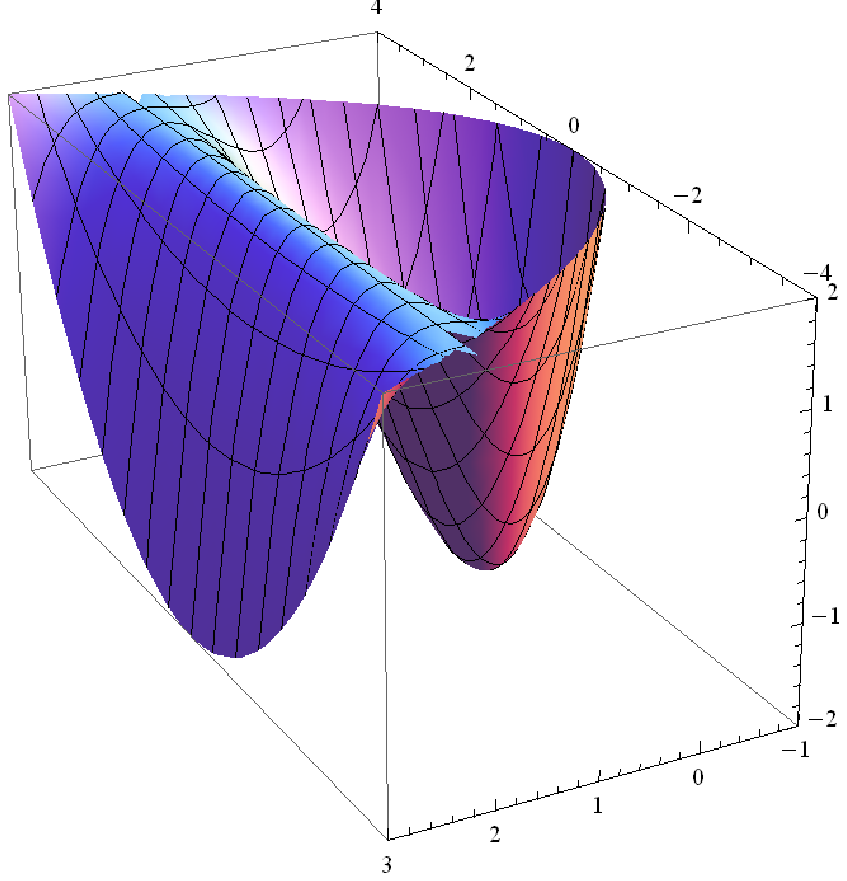}\\
 \caption{The surface $2\mathrm{Re}(J_{\mathbb{Z}_2^{(3)}}) = 0$.} \label{fig:ReJ-Z2^3}
\end{center}
\end{minipage}
\end{figure}

\begin{Rem} \label{Rem:Z2^3<D4^2}
The group $\mathbb{Z}_2^{(3)}$ is a normal subgroup of $D_4^{(2)}$. The fundamental domain $C$ for $D_4^{(2)} = \mathbb{Z}_2^{(3)} \rtimes \mathbb{Z}_2^{(1)}$ in Figure \ref{fig:DomainC-D4^2} is obtained from the fundamental domain for $\mathbb{Z}_2^{(3)}$ in Figure \ref{fig:DomainC-Z2^3} by imposing one extra symmetry which comes from the additional $\mathbb{Z}_2^{(1)}$ action. If we denote $x,y$ for $\mathbb{Z}_2^{(3)}$ given in (\ref{eqn:x,y-Z2^3}) by $x_1,y_1$, and $x,y$ for $D_4^{(2)}$ given in (\ref{eqn:x-D4^2}), (\ref{eqn:y-D4^2}) by $x_2,y_2$, then we see that $x_2 = 2\mathrm{Re}(x_1)$, whilst $y_2 = 2\mathrm{Re}(y_1)+1$, and there is a homomorphism $\xi:\mathfrak{D}_{\mathbb{Z}_2^{(3)}} \rightarrow \mathfrak{D}_{D_4^{(2)}}$ from the domain $\mathfrak{D}_{\mathbb{Z}_2^{(3)}}$ of $\mathbb{Z}_2^{(3)}$ to the domain $\mathfrak{D}_{D_4^{(2)}}$ of $D_4^{(2)}$ such that $\xi(x_1)=2\mathrm{Re}(x_1)$ and $\xi(y_1)=2\mathrm{Re}(y_1)+1$.
However, unlike the case of $\mathbb{Z}_2^{(2)} \lhd D_4^{(1)}$ described in Remark \ref{Rem:Z2^2<D4^1}, the respective Jacobians for $\mathbb{Z}_2^{(3)} \lhd D_4^{(2)}$ are not related by the relation $J_{D_4^{(2)}} = \xi(J_{\mathbb{Z}_2^{(3)}})$. In fact, it is not possible to express $\mathrm{Re}(J_{\mathbb{Z}_2^{(3)}})$ in terms of the $D_4^{(2)}$-invariant variables $x_2$, $y_2$. We introduce a third $D_4^{(2)}$-invariant variable $z_2 := \omega_1\omega_2^{-1} + \omega_1^{-1}\omega_2$, and we have a map $\widetilde{\Phi}:\mathbb{T}^2 \rightarrow \mathbb{R}^3$ given by $\widetilde{\Phi}(\omega_1,\omega_2) = (x_2,y_2,z_2)$ whose image is a surface $\widetilde{\Phi}(C) = \widetilde{\mathfrak{D}} \subset \mathbb{R}^3$, illustrated in Figure \ref{fig:ReJ-Z2^3}.
We can write $\mathrm{Re}(J_{\mathbb{Z}_2^{(3)}})$ in terms of the three $D_4^{(2)}$-invariant variables $x_2$, $y_2$, $z_2$ as
$$2\mathrm{Re}(J_{\mathbb{Z}_2^{(3)}}) = 4 \pi^2 (\omega_1\omega_2^2 - \omega_1^2\omega_2) = x_2^2y_2^2 - 4x_2^2y_2 + 4x_2^2 - 4y_2^3 + 12y_2^2 - 4z_2 - 8.$$
Then setting $2\mathrm{Re}(J_{\mathbb{Z}_2^{(3)}}) = 0$ we obtain a surface $\mathfrak{J}$ in $\mathbb{R}^3$, illustrated in Figure \ref{fig:widetildeD-D4^2}. Now $\mathfrak{J}$ intersects with the surface $\widetilde{\mathfrak{D}}$ only at the closed curve described by the curves $(x_2,(x_2^2-4)/4,2)$ and $(x_2,3,(x_2^2-8)/4)$, and the point $(0,0,-2)$, that is, $2\mathrm{Re}(J_{\mathbb{Z}_2^{(3)}})$ is zero only on this closed curve, which is the boundary of $\widetilde{\mathfrak{D}}$. The projection of this closed curve onto the $x_2,y_2$-plane is then precisely the boundary of the joint spectrum $\mathfrak{D} = \Phi(C)$ for $D_4^{(2)}$, as in Figure \ref{fig:DomainD-D4^2}.
\end{Rem}

\begin{Rem} \label{Rem:D4^1-D4^2}
As noted in Remark \ref{Rem:Z2^2-Z2^3}, the groups $\mathbb{Z}_2^{(2)}$ and $\mathbb{Z}_2^{(3)}$ are conjugate in $GL(2,\mathbb{R})$, with conjugating matrix $H$ as in Remark \ref{Rem:Z2^2-Z2^3}. Since $D_4^{(i)} = \mathbb{Z}_2^{(1)} \rtimes \mathbb{Z}_2^{(i+1)}$ for $i=1,2$, we see that $H$ is also the conjugating matrix for $D_4^{(1)}$ and $D_4^{(2)}$.
As in Remark \ref{Rem:Z2^2-Z2^3}, the origin of this relationship is the fact that the compact, connected Lie group $SU(2) \times SU(2)$ (which has Weyl group $D_4^{(1)}$) is a double cover of the compact, connected, but not simply connected Lie group $SO(4)$ (which has Weyl group $D_4^{(2)}$).
\end{Rem}

\subsection{Spectral measure for $\mathcal{H}^{D_4^{(2)}}_{\rho}$ for $SO(4)$}

\begin{figure}[tb]
\begin{minipage}[t]{9.9cm}
\begin{center}
  \includegraphics[width=90mm]{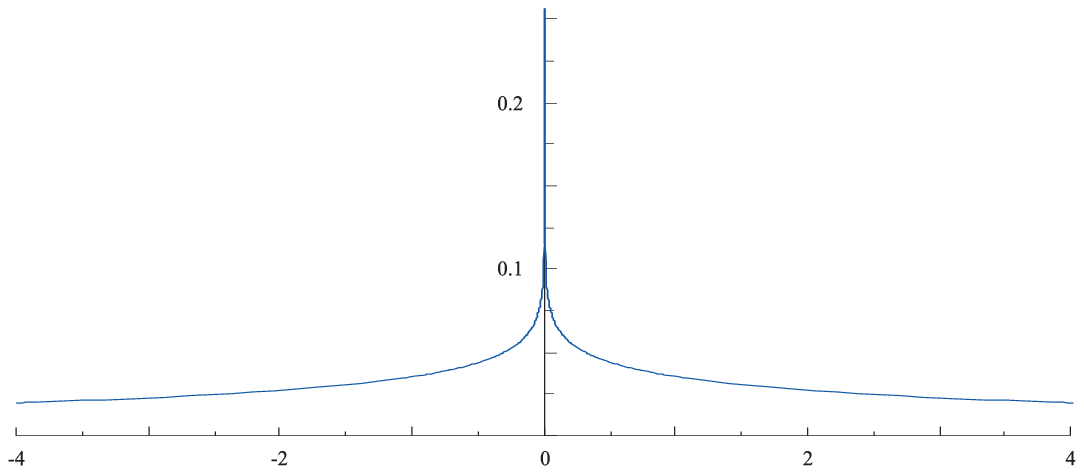}\\
 \caption{weight for $\nu_{\rho_1}$ for $SO(4)$.} \label{fig-Jx-D4^2}
\end{center}
\end{minipage}
\hfill
\begin{minipage}[t]{5.9cm}
\begin{center}
  \includegraphics[width=45mm]{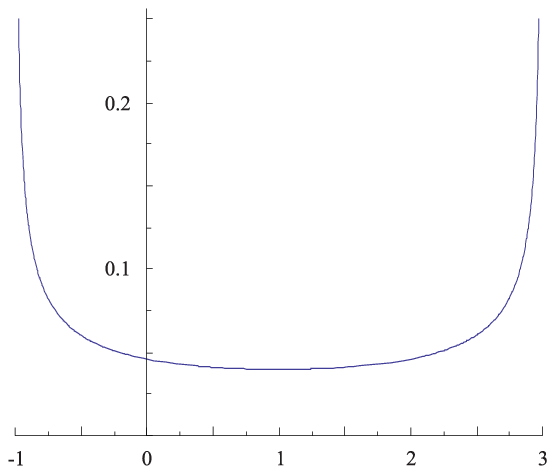}\\
 \caption{weight for $\nu_{\rho_2}$ for $SO(4)$.} \label{fig-Jy-D4^2}
\end{center}
\end{minipage}
\end{figure}

By integrating $|\Gamma| \, |J_{D_4^{(2)}}(x,y)^{-1}|$ over $y,x \in \mathfrak{D}$ respectively we obtain the spectral measure $\nu_{\rho_1}$, $\nu_{\rho_2}$ respectively for ${}^{D_4^{(2)}} \hspace{-1mm} \Delta_{\rho_1}$, ${}^{D_4^{(2)}} \hspace{-1mm} \Delta_{\rho_2}$.
More explicitly, using the expressions for the boundaries of $\mathfrak{D}$ given in Section \ref{sect:D4^2},
the spectral measure $\nu_{\rho_1}$ (over $[-4,4]$) is $\mathrm{d}\nu_{\rho_1}(x) = \mathfrak{J}_1(x) \, \mathrm{d}x$, where $\mathfrak{J}_1(x)$ is given by
$$\mathfrak{J}_1(x) = 4 \int_{(x^2-4)/4}^{3} |J_{D_4^{(2)}}(x,y)|^{-1} \, \mathrm{d}y.$$
The weight $\mathfrak{J}_1(x)$ is an integral of the reciprocal of the square root of a cubic in $y$, and thus can be written in terms of the complete elliptic integral $K(m) = \int_0^{\pi/2} \sqrt{1-m\sin^2\theta}^{-1} \mathrm{d}\theta$ of the first kind. Using \cite[equation 235.00]{byrd/friedman:1971}, we obtain
$$\mathfrak{J}_1(x) = \frac{1}{8\pi^2} \; K\left(\frac{16-x^2}{16}\right).$$
The weight $\mathfrak{J}_1(x)$ is illustrated in Figure \ref{fig-Jx-D4^2}.

The spectral measure $\nu_{\rho_2}$ (over $[-1,3]$) is $\mathrm{d}\nu_{\rho_2}(y) = \mathfrak{J}_2(y) \, \mathrm{d}y$, where $\mathfrak{J}_2(y)$ is given by
\begin{align*}
\mathfrak{J}_2(y) & = 4 \int_{-2\sqrt{y+1}}^{2\sqrt{y+1}} |J_{D_4^{(2)}}(x,y)|^{-1} \, \mathrm{d}x = \frac{1}{\pi^2\sqrt{(y+1)(3-y)}} \int_{-2\sqrt{y+1}}^{2\sqrt{y+1}} \frac{1}{\sqrt{4(y+1)-x^2}} \mathrm{d}x \\
& = \frac{1}{\pi\sqrt{(y+1)(3-y)}}.
\end{align*}
We observe that this is the weight which appears in Section \ref{sect:measure-rho-Z2^2-H}, but shifted $y \rightarrow y-1$.
The weight $\mathfrak{J}_2(y)$ is illustrated in Figure \ref{fig-Jy-D4^2}.

\begin{Thm} \label{thm:measure-rho-D4^2-H}
The spectral measures $\nu_{\rho_1}$, $\nu_{\rho_2}$ (over $\chi_{\rho_1}(\mathbb{T}^2) = [-4,4]$, $\chi_{\rho_2}(\mathbb{T}^2) = [-1,3]$ respectively) for the graphs $\mathcal{H}^{D_4^{(2)}}_{\rho_1}$, $\mathcal{H}^{D_4^{(2)}}_{\rho_2}$ respectively, for $SO(4)$ are given by
$$\mathrm{d}\nu_{\rho_1}(x) = \frac{1}{8\pi^2} \; K\left(\frac{16-x^2}{16}\right) \, \mathrm{d}x, \qquad \qquad
\mathrm{d}\nu_{\rho_2}(y) = \frac{1}{\pi\sqrt{(y+1)(3-y)}} \, \mathrm{d}y.$$
\end{Thm}

\subsection{Spectral measure for $\mathcal{G}^{D_4^{(2)}}_{\rho}$ for $SO(4)$}

We now consider the spectral measure (over $\chi_{\rho_j}(\mathbb{T}^2)$) for $\mathcal{G}^{D_4^{(2)}}_{\rho_j}$, $j=1,2$.
The adjacency matrix $\Delta_{\rho_1}$ for the graph $\mathcal{G}^{D_4^{(2)}}_{\rho_1}$ is self-adjoint, and thus its spectral measure is uniquely determined by its moments $\varphi(\Delta_{\rho_1}^m)$. It has spectrum $\sigma(\Delta_{\rho_1}) = \chi_{\rho_1}(\mathbb{T}^2) = [-4,4]$. We define a state $\varphi$ by $\varphi( \, \cdot \, ) = \langle \cdot \Omega, \Omega \rangle$, where $\Omega$ is vector in $\ell^2(\mathcal{G}^{D_4^{(2)}}_{\rho_1})$ corresponding to the distinguished vertex $\ast$, which is chosen to be the vertex with lowest Perron-Frobenius weight, which in this case is the apex vertex, i.e. the vertex with only one edge attached to it. Thus the $m^{\mathrm{th}}$ moment $\varphi(\Delta_{\rho_1}^m)$ counts the number of closed paths of length $m$ on $\mathcal{G}^{D_4^{(2)}}_{\rho_1}$ which start and end at $\ast$.

\begin{Lemma}
The number $\alpha_m$ of closed paths of length $m$ on $\mathcal{G}^{D_4^{(2)}}_{\rho_1}$ which start and end at $\ast$ is given
$$\alpha_m \;\; = \;\; \left\{ \begin{array}{cl} c_{m/2} & \textrm{ if } m \textrm{ is even}, \\ 0 & \textrm{ if } m \textrm{ is odd}. \end{array} \right.$$
\end{Lemma}

\emph{Proof}:
Since $\mathcal{G}^{D_4^{(2)}}_{\rho_1}$ is bipartite it is clear that the odd moments are all zero.
Any closed path of length $2n$ is given by a pair of paths which start at $\ast$ and which both end at the same vertex at the $n^{\mathrm{th}}$ level of the Bratteli diagram whose inclusion graph at each level is given by $\mathcal{G}^{D_4^{(2)}}_{\rho_1}$, and thus $\varphi(\Delta_{\rho_1}^m)$ counts the number of such pairs of paths. Any path of length $n$ on the Bratteli diagram is equivalent to a path in the lattice $\mathbb{N}^3$ from $(0,0,0)$ to $(i,j,n)$, for some $(i,j) \in \mathbb{N}^2$, with steps belonging to the set $\mathcal{S} = \{ (1,1,1), (1,-1,1), (-1,1,1), (-1,-1,1) \}$, where here we think of $\mathbb{N}$ as including 0. Then a closed path of length $2n$ is equivalent to a path in the lattice $\mathbb{N}^3$ from $(0,0,0)$ to $(0,0,2n)$ with steps belonging to the set $\mathcal{S}$. Any such path is given by a sequence of steps $(\epsilon_j^1,\epsilon_j^2,1)$ from $\mathcal{S}$, $j=1,2,\ldots,2n$. Thus we obtain two sequences $\epsilon^l = (\epsilon_1^l,\epsilon_2^l,\ldots,\epsilon_{2n}^l)$, for $l=1,2$, where $\epsilon_j^l \in \{\pm1\}$, $\sum_j \epsilon_j^l = 0$ and $\sum_{j=1}^k \epsilon_j^l \geq 0$ for each $k$ such that $\epsilon_k^l = -1$. Such a sequence is enumerated by the Catalan numbers $c_n$ \cite[Corollary 6.2.3(ii)]{stanley:1999}, and thus the number of closed paths of length $2n$ on $\mathcal{G}^{D_4^{(2)}}_{\rho_1}$, starting and ending at $\ast$, is given by the number of pairs of such sequences, which is $c_n^2$.
\hfill
$\Box$

\begin{Rem} \label{Rem:OEIS-conjecture}
This proves a conjecture from the Online Encyclopedia of Integer Sequences (OEIS) \cite{OEIS:2010}, namely that the number of walks on $\mathbb{N}^2$ starting and ending at $(0,0)$ and consisting of $2n$ steps taken from $\{(-1, -1), (-1, 1), (1, -1), (1, 1)\}$, is given by the squared Catalan numbers A001246. Moreover, from \cite[Proposition 8]{bousquet-melou/mishna:2010}, we also find that the squared Catalan numbers satisfy the following relation with multinomial coefficients $(m_1,m_2,m_3,m_4)! = (\sum_{i=1}^4 m_i)!/(m_1!m_2!m_3!m_4!)$:
$$c_n^2 = \sum_{k=0}^n (k,k,n-k,n-k)! + \sum_{k=0}^{n-1} (k,k+2,n-k-1,n-k-1)! -2\sum_{k=0}^{n-1} (k,k+1,n-k-1,n-k)!.$$
\end{Rem}

Thus the $2n^{\mathrm{th}}$ moment is $\varphi(\Delta_{\rho_1}^{2n}) = c_n^2$.
The following result is claimed in \cite{OEIS:2010}, but we have not found a proof of the claim in the literature. Therefore we give a proof here.

\begin{Lemma} \label{lemma:M(z)-c_n^2}
The generating function for the squared Catalan numbers is given by
$$M(z) = \frac{1}{\pi z}E(16z) - \frac{1-16z}{2\pi z}K(16z) - \frac{1}{4z},$$
where $E(m) = \int_0^{\pi/2} \sqrt{1-m\sin^2\theta} \mathrm{d}\theta$ is the complete elliptic integral of the second kind.
\end{Lemma}

\emph{Proof}:
The generating function $M(z)$ for the squared Catalan numbers satisfies the linear differential equation \cite{cleary/elder/rechnitzer/taback:2010}
\begin{equation} \label{eqn:diff_eqn-c_n^2}
z^2(1-16z) \frac{\mathrm{d}^2M}{\mathrm{d}z^2} + x(3-32z) \frac{\mathrm{d}M}{\mathrm{d}z} + (1-4z)M(z) \;\; = \;\; 1,
\end{equation}
where it is easy to check that the coefficient of $z^n$ satisfies the recurrence relation $(n+2)^2 c_{n+1}^2 = 4(2n+1)^2 c_n^2$ for the squared Catalan numbers. Now since (see e.g. \cite[$\S$22.736]{whittaker/watson:1927})
$$\frac{\mathrm{d}}{\mathrm{d}z}K(16z) = \frac{E(16z)-(1-16z)K(16z)}{2(1-16z)z}, \qquad \frac{\mathrm{d}}{\mathrm{d}z}E(16z) = \frac{E(16z)-K(16z)}{2z},$$
we have
\begin{eqnarray*}
\frac{\mathrm{d}}{\mathrm{d}z}M(z) & = & \frac{1-16z}{4\pi z^2}K(16z) - \frac{3}{4\pi z^2}E(16z) + \frac{1}{4z^2}, \\
\frac{\mathrm{d}^2}{\mathrm{d}z^2}M(z) & = & \frac{24z-1}{4\pi z^3}K(16z) + \frac{5}{4\pi z^3}E(16z) - \frac{1}{2z^3}.
\end{eqnarray*}
Thus a simple calculation shows that $M(z)$ satisfies the differential equation (\ref{eqn:diff_eqn-c_n^2}) for the squared Catalan numbers.
\hfill
$\Box$

A compactly supported probability measure $\mu$ on $\mathbb{R}$ can be recovered from the generating function of its moments by $\mu = -\lim_{\varepsilon \rightarrow 0} \textrm{Im}\left(G_{\mu}(t+i\varepsilon)\right)/\pi$ (see e.g. \cite{hiai/petz:2000}),
where $G_{\mu}(z) = \int_{\mathbb{R}} (z-t)^{-1} \mathrm{d}\mu(t)$ is the Cauchy transform for $\mu$, which has power series expansion about $z=\infty$ given by $G_{\mu}(z) = \sum_{n=0}^{\infty} m_n z^{-n-1}$, where $m_n$ are the moments of $\mu$. Thus $G_{\mu}(z)$ is related to the generating series for the moments by
$G_{\mu}(z) = M(z^{-1})/z$.

\begin{Thm} \label{thm:measure-rho-D4^2-G-x}
The spectral measure $\nu_{\rho_1}$ (over $\chi_{\rho_1}(\mathbb{T}^2) = [-4,4]$) for the graph $\mathcal{G}^{D_4^{(2)}}_{\rho_1}$ for $SO(4)$ is given by
$$\mathrm{d}\nu_{\rho_1}(x) = \frac{1}{2\pi^2 |x|} \; \left( 2x^2 E\left(\frac{x^2-16}{x^2}\right) - (x^2+16) K\left(\frac{x^2-16}{x^2}\right) \right) \, \mathrm{d}x.$$
\end{Thm}

\emph{Proof}:
The generating function for the moments of $\nu_{\rho_1}$ is $M_1(z) = M(z^2)$, where $M(z)$ is given in Lemma \ref{lemma:M(z)-c_n^2}, since the odd moments of $\nu_{\rho_1}$ are all zero. Thus the Cauchy transform for $M_1(z)$ is
$$G_{\nu_{\rho_1}}(z) \;\; = \;\; \frac{z}{\pi}E(16z^{-2}) - \frac{z^2-16}{2\pi z}K(16z^{-2}) - \frac{z}{4}.$$
Now, for $z \in [0,4]$,
\begin{eqnarray*}
\lefteqn{ \lim_{\varepsilon \rightarrow 0} \mathrm{Im} \left( \frac{t+i\varepsilon}{\pi}E(16(t+i\varepsilon)^{-2}) \right) \;\; = \;\; \lim_{\varepsilon \rightarrow 0} \mathrm{Im} \left( \frac{t+i\varepsilon}{\pi} \int_0^{\pi/2} \sqrt{1-16(t+i\varepsilon)^{-2} \sin^2\theta} \, \mathrm{d}\theta \right) } \\
& = & \mathrm{Im} \left( \frac{t}{\pi} \int_0^{\pi/2} \sqrt{1-16t^{-2} \sin^2\theta} \, \mathrm{d}\theta \right) \;\; = \;\; \mathrm{Im} \left( \frac{ti}{\pi} \int_0^{\pi/2} \sqrt{16t^{-2} \sin^2\theta-1} \, \mathrm{d}\theta \right) \\
& = & \frac{t}{\pi} \int_{\sin^{-1}|t/4|}^{\pi/2} \sqrt{16t^{-2} \sin^2\theta-1} \, \mathrm{d}\theta,
\end{eqnarray*}
since $16t^{-2} \sin^2\theta \geq 1$ for $\theta \geq \sin^{-1}|t/4|$.
Thus
\begin{eqnarray*}
\lefteqn{ \lim_{\varepsilon \rightarrow 0} \mathrm{Im} \left( \frac{t+i\varepsilon}{\pi}E(16(t+i\varepsilon)^{-2}) \right) } \\
& = & \frac{-ti}{\pi} \left( \int_0^{\pi/2} \sqrt{1-16t^{-2} \sin^2 \theta} \, \mathrm{d}\theta - \int_0^{\sin^{-1}|t/4|} \sqrt{1-16t^{-2} \sin^2 \theta} \, \mathrm{d}\theta \right) \\
& = & \frac{-ti}{\pi} \left( E(16t^{-2}) - E(\sin^{-1}|t/4|,16t^{-2}) \right).
\end{eqnarray*}
Similarly,
\begin{eqnarray*}
\lefteqn{ \lim_{\varepsilon \rightarrow 0} \mathrm{Im} \left( \frac{(t+i\varepsilon)^2-16}{2\pi(t+i\varepsilon)}K(16(t+i\varepsilon)^{-2}) \right) } \\
& = & \frac{-(t^2-16)i}{2\pi t} \left( K(16t^{-2}) - F(\sin^{-1}|t/4|,16t^{-2}) \right).
\end{eqnarray*}
Since, $\lim_{\varepsilon \rightarrow 0} \mathrm{Im} ((t+i\varepsilon)/4) = 0$, we obtain for $z \in [0,4]$ that
\begin{align*}
\mathrm{d}\nu_{\rho_1}(x) & = \frac{i}{2\pi^2 x} \; \left( 2x^2 \left( E(16x^{-2}) - E(\sin^{-1}|x/4|,16x^{-2}) \right) \right. \\
& \qquad \left. - (x^2-16) \left( K(16x^{-2}) - F(\sin^{-1}|x/4|,16x^{-2}) \right) \right) \, \mathrm{d}x,
\end{align*}
where $F(\phi,m)$, $E(\phi,m)$ are the elliptic integrals of the first, second kind respectively, defined by
$$F(\phi,m) = \int_0^{\phi} \frac{1}{\sqrt{1-m \sin^2 \theta}} \, \mathrm{d}\theta, \qquad E(\phi,m) = \int_0^{\phi} \sqrt{1-m \sin^2 \theta} \, \mathrm{d}\theta,$$
so that $K(m) = F(\pi/2,m)$, $E(m) = E(\pi/2,m)$ are the complete elliptic integrals of the first, second kind respectively.
Similarly for $z \in [-4,0]$ we obtain the same expression for $\mathrm{d}\nu_{\rho_1}(x)$ but multiplied by a factor of $-1$.
Then the result follows from the elliptic integral identities 111.09 in \cite{byrd/friedman:1971}.
\hfill
$\Box$

We now consider the graph $\mathcal{G}^{D_4^{(2)}}_{\rho_2}$ for the second fundamental generator $\rho_2$. The connected component of the distinguished vertex $\ast$ is the infinite graph $\mathcal{A}^{(\infty_{\mathrm{odd}})\ast} := \lim_{n \rightarrow \infty} \mathcal{A}^{(2n+1)\ast}$, whilst the connected component of the vertex $(1,0)$ immediately to its right is the infinite graph $\mathcal{A}^{(\infty_{\mathrm{even}})\ast} := \lim_{n \rightarrow \infty} \mathcal{A}^{(2n)\ast}$, where the graphs $\mathcal{A}^{(l)\ast}$ are $SU(3)$ $\mathcal{ADE}$ graphs which classify the conjugate modular invariant partition functions for $SU(3)$ integrable statistical mechanical models, and which also yield a non-negative integer matrix representation (nimrep) of the fusion rules of $SU(3)$ at level $l-3$ (see e.g. \cite{evans/pugh:2012iv} for more details about the other $SU(3)$ $\mathcal{ADE}$ graphs, $SU(3)$ modular invariants and nimrep theory). Spectral measures for the $\mathcal{A}^{(l)\ast}$ graphs were computed in \cite[Theorem 7]{evans/pugh:2009v}. The distinguished vertex $\ast$ with lowest Perron-Frobenius weight in each case is the vertex labelled 1 in \cite[Fig. 13]{evans/pugh:2009v}. However, the spectral measures determined in \cite[Theorem 7]{evans/pugh:2009v} are for the case where $\ast$ is the vertex labelled $\lfloor (l-1)/2 \rfloor$ in \cite[Fig. 13]{evans/pugh:2009v}, and not for the case where $\ast$ is the vertex with lowest Perron-Frobenius weight as claimed. We thus provide the correct result here for the spectral measure of $\mathcal{A}^{(l)\ast}$ where the distinguished vertex $\ast$ is chosen to be the vertex with lowest Perron-Frobenius weight:

\begin{Thm} \label{Thm:spec_measure-A(l)star-corrected}
The spectral measure of $\mathcal{A}^{(l)\ast}$, $l < \infty$, (over $\mathbb{T}$), where the distinguished vertex $\ast$ is chosen to be the vertex with lowest Perron-Frobenius weight, is
\begin{equation}
\mathrm{d}\varepsilon(u) = \alpha(u) \mathrm{d}_{l/2}u,
\end{equation}
where $\mathrm{d}_{l/2}u$ is the uniform measure over $l^{\mathrm{th}}$ roots of unity, and
$$\alpha(u) = \left\{ \begin{array}{cl}  2\mathrm{Im}(u)^2 & \textrm{ for } l \textrm{ even,} \\
                                        1-\mathrm{Re}(u) & \textrm{ for } l \textrm{ odd.} \end{array} \right.$$
\end{Thm}

\emph{Proof}:
From \cite{gaberdiel/gannon:2002} the eigenvectors of $\mathcal{A}^{(l)\ast}$ are $\psi^{\lambda}_a = 2 \sqrt{l^{-1}} \sin(2 \pi a \lambda/l)$, where $\lambda, a = 1,2,\ldots,\lfloor (l-1)/2 \rfloor$. For each eigenvalue $\beta^{(\lambda)}$ of $\mathcal{A}^{(l)\ast}$, the eigenvector entry corresponding to the distinguished vertex $\ast = \lfloor (l-1)/2 \rfloor$ is $\psi^{\lambda}_{\ast} = 2 \sqrt{l^{-1}} \sin(2 \pi \lfloor (l-1)/2 \rfloor \lambda/l)$.

In the case where $l=2n$ is even,
$$\psi^{\lambda}_{\ast} = \frac{2}{\sqrt{2n}} \sin\left(\frac{\pi}{n} (n-1) \lambda\right) \;\; = \;\; (-1)^{\lambda+1} \frac{2}{\sqrt{2n}} \sin\left(\frac{2\pi}{2n} \lambda\right) \;\; = \;\; (-1)^{\lambda+1} \frac{2}{\sqrt{l}} \sin\left(\frac{2\pi}{l} \lambda\right).$$
Then $|\psi^{\lambda}_{\ast}|^2 \;\; = \;\; 4 \sin^2(2\pi\lambda/l)/l$ and the result follows from \cite[Theorem 7]{evans/pugh:2009v}.

In the case where $l=2n+1$ is odd,
\begin{eqnarray*}
\psi^{\lambda}_{\ast} & = & \frac{2}{\sqrt{2n+1}} \sin\left(\frac{2\pi}{2n+1} n \lambda\right) \;\; = \;\; \frac{2}{\sqrt{2n+1}} \sin\left(\pi \left(1-\frac{1}{2n+1}\right) \lambda\right) \\
& = & (-1)^{\lambda+1} \frac{2}{\sqrt{2n+1}} \sin\left(\frac{\pi}{2n+1} \lambda\right) \;\; = \;\; (-1)^{\lambda+1} \frac{2}{\sqrt{l}} \sin\left(\frac{\pi}{l} \lambda\right) \\
& = & (-1)^{\lambda+1} \frac{1}{\sqrt{l}} \sqrt{ 1 - \cos\left(\frac{2\pi}{l} \lambda\right) }.
\end{eqnarray*}
Then $|\psi^{\lambda}_{\ast}|^2 \;\; = \;\; (1-\cos(2\pi\lambda/l))/l$ and the result follows as in \cite[Theorem 7]{evans/pugh:2009v} but with the weighting $2\sin^2(2\pi\lambda/l)$ replaced by $1-\cos(2\pi\lambda/l)$.
\hfill
$\Box$

\begin{Thm} \label{Thm:spec_measure-A(infty)star}
The spectral measure of the infinite graphs $\mathcal{A}^{(\infty_\mathrm{even})\ast}$, $\mathcal{A}^{(\infty_\mathrm{odd})\ast}$, (over $\mathbb{T}$), where the distinguished vertex $\ast$ is chosen to be the vertex with lowest Perron-Frobenius weight, is
\begin{equation}
\mathrm{d}\varepsilon(u) = \alpha(u) \mathrm{d}u,
\end{equation}
where $\mathrm{d}u$ is the uniform Lebesgue measure over $\mathbb{T}$, and
$$\alpha(u) = \left\{ \begin{array}{cl}  2\mathrm{Im}(u)^2 & \textrm{ for } \mathcal{A}^{(\infty_\mathrm{even})\ast}, \\
                                        1-\mathrm{Re}(u) & \textrm{ for } \mathcal{A}^{(\infty_\mathrm{odd})\ast}. \end{array} \right.$$
Over their spectrum $[-1,3]$, the corresponding spectral measures $\nu_{\mathrm{even}}$, $\nu_{\mathrm{odd}}$ of the infinite graphs $\mathcal{A}^{(\infty_\mathrm{even})\ast}$, $\mathcal{A}^{(\infty_\mathrm{odd})\ast}$ respectively, are
\begin{equation}
\mathrm{d}\nu_{\mathrm{even}}(y) = \frac{1}{2\pi}\sqrt{(3-y)(y+1)} \, \mathrm{d}y, \qquad \mathrm{d}\nu_{\mathrm{odd}}(y) = \frac{1}{2\pi} \frac{\sqrt{3-y}}{\sqrt{y+1}} \, \mathrm{d}y.
\end{equation}
\end{Thm}

\emph{Proof}:
For $\mathcal{A}^{(\infty_\mathrm{even})\ast}$ the result is given in \cite[$\S$7.3]{evans/pugh:2009v}. For $\mathcal{A}^{(\infty_\mathrm{odd})\ast}$, with the change of variables $y=2\cos(2\pi\theta)+1$ as in \cite[$\S$7.3]{evans/pugh:2009v}, where $u=e^{2\pi i \theta}$, we have $1-\cos(2\pi\theta) = (3-y)/2$, and the result follows in a similar way to that for $\mathcal{A}^{(\infty_\mathrm{even})\ast}$.
\hfill
$\Box$

Thus we obtain the following result:
\begin{Thm} \label{thm:measure-rho-D4^2-G-y}
The spectral measure $\nu_{\rho_2}$ (over $\chi_{\rho_2}(\mathbb{T}^2) = [-1,3]$) for the graph $\mathcal{G}^{D_4^{(2)}}_{\rho_2}$ for $SO(4)$ is given by
$$\mathrm{d}\nu_{\rho_2}(y) = \frac{1}{2\pi} \frac{\sqrt{3-y}}{\sqrt{y+1}} \, \mathrm{d}y.$$
\end{Thm}

\begin{figure}[tb]
\begin{minipage}[t]{9.9cm}
\begin{center}
  \includegraphics[width=90mm]{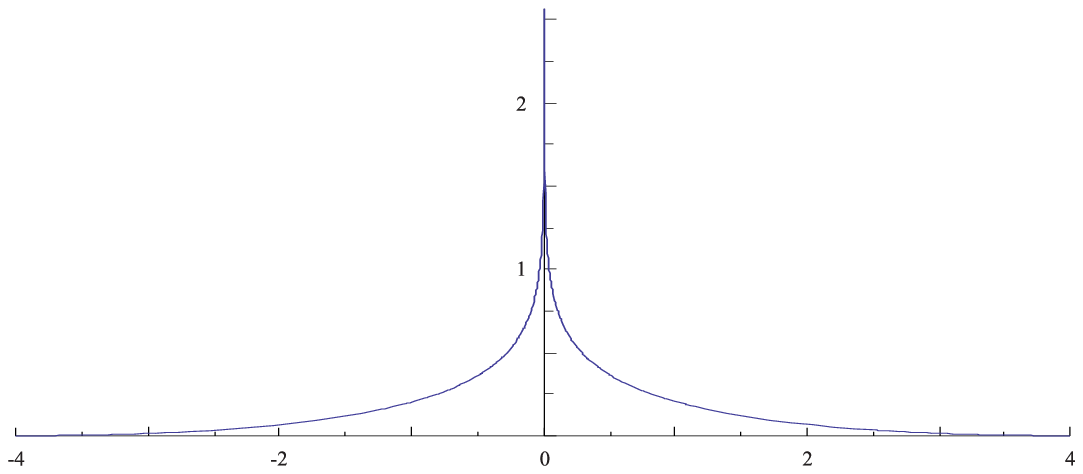}\\
 \caption{weight for $\nu_{\rho_1}$ for $SO(4)$.} \label{fig-Jx-D4^2-G}
\end{center}
\end{minipage}
\hfill
\begin{minipage}[t]{5.9cm}
\begin{center}
  \includegraphics[width=45mm]{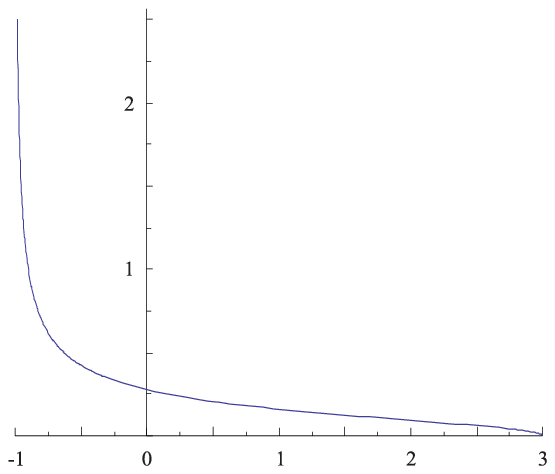}\\
 \caption{weight for $\nu_{\rho_2}$ \mbox{for $SO(4)$}.} \label{fig-Jy-D4^2-G}
\end{center}
\end{minipage}
\end{figure}

\section{$D_6^{(1)}$: $PSU(3)$} \label{sect:D6(1)}

\begin{figure}[tb]
\begin{minipage}[t]{7.9cm}
\begin{center}
  \includegraphics[width=60mm]{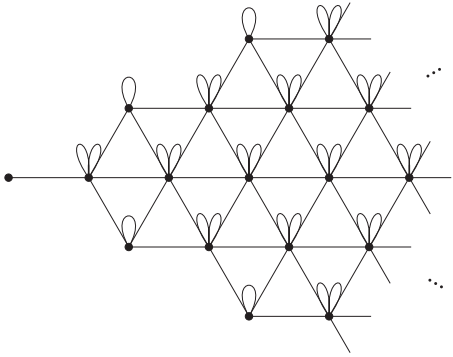}
 \caption{Infinite graph $\mathcal{G}^{D_6^{(1)}}_{\rho_1}$ for $PSU(3)$} \label{Fig-Graph_D6^1-G10}
\end{center}
\end{minipage}
\hfill
\begin{minipage}[t]{7.9cm}
\begin{center}
  \includegraphics[width=60mm]{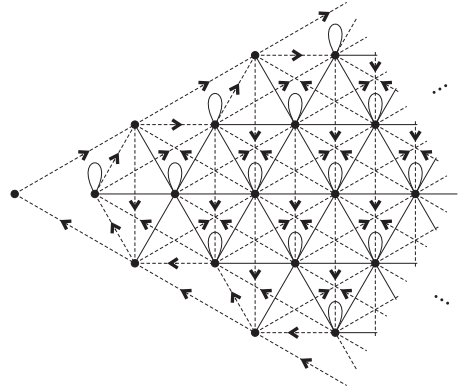}
 \caption{Infinite graph $\mathcal{G}^{D_6^{(1)}}_{\rho_2}$ for $PSU(3)$} \label{Fig-Graph_D6^1-G11}
\end{center}
\end{minipage}
\end{figure}

\begin{figure}[tb]
\begin{minipage}[t]{7.9cm}
\begin{center}
  \includegraphics[width=60mm]{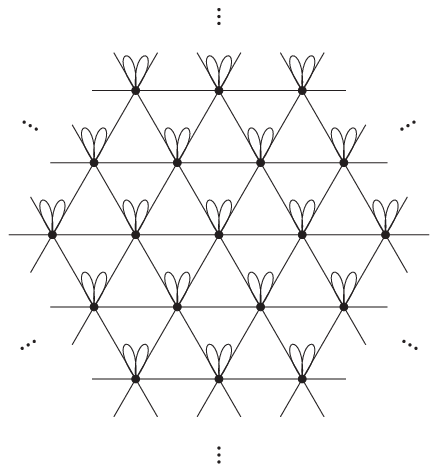}
 \caption{Infinite graph $\mathcal{H}^{D_6^{(1)}}_{\rho_1}$ for $PSU(3)$} \label{Fig-Graph_D6^1-WG10}
\end{center}
\end{minipage}
\hfill
\begin{minipage}[t]{7.9cm}
\begin{center}
  \includegraphics[width=60mm]{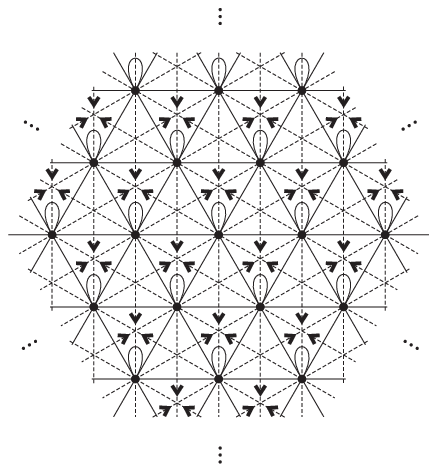}
 \caption{Infinite graph $\mathcal{H}^{D_6^{(1)}}_{\rho_2}$ for $PSU(3)$} \label{Fig-Graph_D6^1-WG11}
\end{center}
\end{minipage}
\end{figure}

The subgroup $D_6^{(1)} \in GL(2,\mathbb{Z})$ is generated by the matrices $T_3$, $T_2'$ of Section \ref{sect:intro}.
It is the Weyl group for the connected compact Lie group $PSU(3)$.
We choose fundamental generators $\rho_1 = (1,0)$, $\rho_2 = (1,1)$.
The graphs $\mathcal{G}^{D_6^{(1)}}_{\rho}$, $\mathcal{H}^{D_6^{(1)}}_{\rho}$ for $\rho = \rho_1, \rho_2$ are illustrated in Figures \ref{Fig-Graph_D6^1-G10}-\ref{Fig-Graph_D6^1-WG11}. The dashed lines in Figures \ref{Fig-Graph_D6^1-G11} and \ref{Fig-Graph_D6^1-WG11} are directed edges, whilst the solid lines are undirected edges.
For $\omega_j = e^{2\pi i \theta_j} \in \mathbb{T}$, $j=1,2$, let
\begin{align}
x & := x_{\rho_1} = 2 + \omega_1 + \omega_1^{-1} + \omega_2 + \omega_2^{-1} + \omega_1\omega_2^{-1} + \omega_1^{-1}\omega_2 \label{eqn:x-D6^1} \\
& = 2 + 2\cos(2\pi \theta_1) + 2\cos(2\pi \theta_2) + 2\cos(2\pi(\theta_1-\theta_2)), \nonumber \\
y & := x_{\rho_2} = 1 + \omega_1\omega_2 + \omega_1\omega_2^{-2} + \omega_1^{-2}\omega_2 + \omega_1 + \omega_1^{-1} + \omega_2 + \omega_2^{-1} + \omega_1\omega_2^{-1} + \omega_1^{-1}\omega_2 \label{eqn:y-D6^1} \\
& = 1 + e^{2\pi i (\theta_1+\theta_2)} + e^{2\pi i (\theta_1-2\theta_2)} + e^{2\pi i (-2\theta_1+\theta_2)} + 2\cos(2\pi \theta_1) + 2\cos(2\pi \theta_2) + 2\cos(2\pi(\theta_1-\theta_2)), \nonumber
\end{align}
and denote by $\Psi$ be the map $\Psi_{\rho_1,\rho_2}: (\omega_1,\omega_2) \mapsto (x,y)$.

\begin{figure}[tb]
\begin{minipage}[t]{5.7cm}
\begin{center}
  \includegraphics[width=55mm]{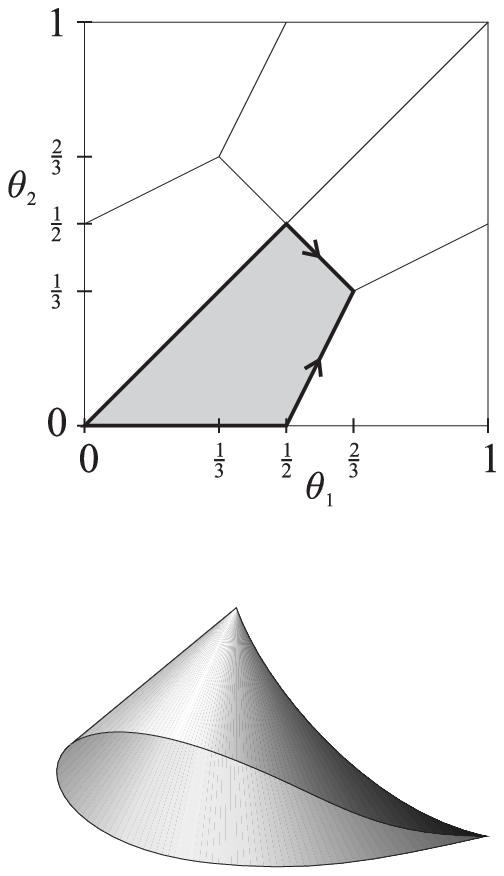}\\
 \caption{$PSU(3)$: A fundamental domain $C$ of $\mathbb{T}^2/D_6^{(1)}$ (top) and the surface $\mathbb{T}^2/D_6^{(1)}$ (bottom).} \label{fig:DomainC-D6^1}
\end{center}
\end{minipage}
\hfill
\begin{minipage}[t]{10.1cm}
\begin{center}
  \includegraphics[width=100mm]{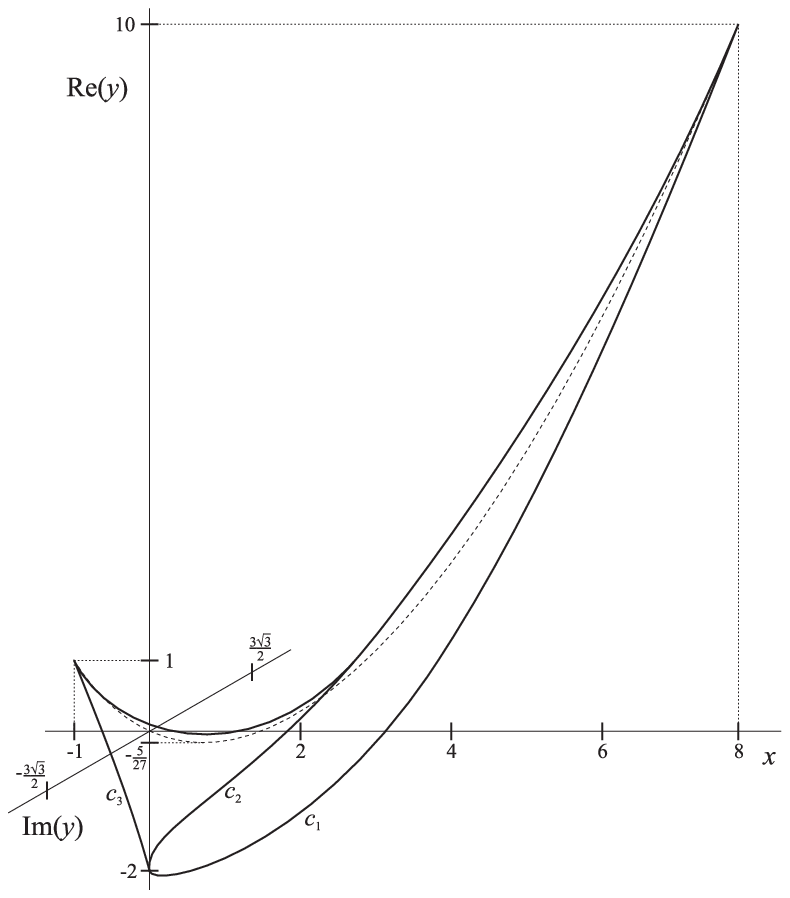}\\
 \caption{The domain $\mathfrak{D} = \Psi(C)$ for $PSU(3)$.} \label{fig:DomainD-D6^1}
\end{center}
\end{minipage}
\end{figure}

A fundamental domain of $\mathbb{T}^2/D_6^{(1)}$ is illustrated in Figure \ref{fig:DomainC-D6^1}, where the boundaries marked by arrows are identified. Then $\mathbb{T}^2/D_6^{(1)}$ is a surface, illustrated by the bottom figure in Figure \ref{fig:DomainC-D6^1}. It is the surface of a cone, with two singular points, one at the apex of the cone, and one on its boundary at the bottom, which is the only boundary of the surface.
The boundary of $C$ given by $\theta_2 = 0$, where $\theta_1 \in [0,1/2]$, yields the curve $c_1$ given by the parametric equations
$x = 4\cos(2\pi\theta_1) + 4$, $y = e^{-4\pi i \theta_1} + 2e^{2\pi i \theta_1} + 4\cos(4\pi\theta_1) + 3$,
where $x \in [0,8]$ for $\theta_1 \in [0,1/2]$.
We have
\begin{align}
\mathrm{Re}(y) & = (y+\overline{y})/2 = \cos(4\pi\theta_1) + 6\cos(2\pi\theta_1) + 3 = 2\cos^2(2\pi\theta_1) + 6\cos(2\pi\theta_1) + 3 \nonumber \\
& = (x^2 + 4x - 16)/8, \label{eqn:y-x-D6^1-1} \\
\mathrm{Im}(y) & = (y-\overline{y})/2i = \sin(4\pi\theta_1) - 2\sin(2\pi\theta_1) = 2\sin(2\pi\theta_2)(\cos(2\pi\theta_1)-1) \nonumber \\
& = \pm\sqrt{x(8-x)^3}/8, \label{eqn:y-x-D6^1-2}
\end{align}
where $x(8-x)^3 \geq 0$ for $x \in [0,8]$. Since $\mathrm{Im}(\Psi(e^{2\pi i/3},1)) > 0$, the boundary $c_1$ of $\mathfrak{D}$ is given by $8y = x^2 + 4x - 16 + i\sqrt{x(8-x)^3}$.
Similarly, the boundary of $C$ given by $\theta_1 = \theta_2$, where $\theta_1 \in [0,1/2]$, yields the curve $c_2$ given by the same parametric equations except with $y \leftrightarrow \overline{y}$. Thus the boundaries $c_1$, $c_2$ of $\mathfrak{D}$ are given by
$8y = x^2 + 4x - 16 + i\sqrt{x(8-x)^3}$, $8y = x^2 + 4x - 16 - i\sqrt{x(8-x)^3}$ respectively.
The point $(\theta_1,\theta_2) = (2/3,1/3)$ on the boundary of $C$ maps to $\Psi(e^{-2\pi i/3},e^{2\pi i/3}) = (-1,1)$, which is a singular point in $\mathfrak{D}$.
The boundaries of $C$ given by $\theta_2 = -\theta_1$ and $\theta_2 = 2\theta_1$, for $\theta_1 \in [1/2,2/3]$, are identified, and yield the curve $c_3$ given by the parametric equations
$x = 2\cos(4\pi\theta_1) + 4\cos(2\pi\theta_1) + 2 = 4\cos^2(2\pi\theta_1) + 4\cos(2\pi\theta_1)$ and
$y = 2\cos(6\pi\theta_1) + 2\cos(4\pi\theta_1) + 4\cos(2\pi\theta_1) + 1 = 8\cos^3(2\pi\theta_1) + 4\cos^2(2\pi\theta_1) - 2\cos(2\pi\theta_1)$.
Then writing $\cos(2\pi \theta_1)$ in terms of $x$, we obtain that $c_3$ maps to $y = -2x-1 \pm (x+1)^{3/2}$ in $\mathfrak{D}$, where $\pm$ is taken to be $-$. On $c_3$, $y \in \mathbb{R}$. The curve $c_3$ lies on the interior of the surface $\mathfrak{D}$, except at the point $(0,-2)$ on $c_1$, $c_2$, and the singular point $(-1,1)$.
Then $\mathfrak{D} = \Psi(C)$ is the surface of the cone illustrated in Figure \ref{fig:DomainD-D6^1}, with boundary given by the curves $c_1$, $c_2$. The equation of the surface is
\begin{equation} \label{eqn:surfaceD-D6^1}
\mathrm{Im}(y)^2=x^3-\mathrm{Re}(y)^2-4x\mathrm{Re}(y)-x^2-2\mathrm{Re}(y)-x.
\end{equation}
Away from the curve $c_3$, $y \in \mathbb{R}$ only for $\theta_1 = 2\theta_2$ in the fundamental domain $C$, where $\theta_2 \in [0,1/3]$. Such $y$ is given by the dashed curve in Figure \ref{fig:DomainD-D6^1}. This curve is also given by $\theta_2 = -\theta_1$ for $\theta_1 \in [2/3,1]$ (i.e. in $T_3(C)$), and thus it is given by the same equation $y = -2x-1 \pm (x+1)^{3/2}$, where now $\pm$ is taken to be $+$.

We can also write the curves $c_i$ as functions of $y$. For $c_1$ and $c_2$, we have from (\ref{eqn:y-x-D6^1-1}) and $x \geq -1$ that $x = 2(-1+\sqrt{5+\mathrm{Re}(y)})$.
For $\theta_2 = -\theta_1$, $\theta_2 = 2\theta_1$ and $\theta_1 = 2\theta_2$, i.e. for $y \in \mathbb{R}$, the parametric equation for $y$ is a cubic in $\cos(2\pi \theta_1)$. Solving this cubic we obtain
\begin{align*}
x & = 4p_2(y)+4p_2(y)^2, \qquad y \in [-2,1], \\
x & = 4p_3(y)+4p_3(y)^2, \qquad y \in [-5/27,1], \\
x & = 4p_1(y)+4p_1(y)^2, \qquad y \in [-5/27,10],
\end{align*}
where $p_i$ is given by $6p_i(y) = -1 + 2^{-1/3} \epsilon_i P + 2^{-1/3} 4\overline{\epsilon_i} P^{-1}$, for $\epsilon_j = e^{2 \pi i (j-1)/3}$ , $P = (27y-11 + \sqrt{3^3(27y^2-22y-5)})^{1/3}$.
The first equation above is the equation for the curve $c_3$, whilst the latter two are the equations for the dashed curve in Figure \ref{fig:DomainD-Z2^3}.
The surface $\mathfrak{D}$ is the joint spectrum $\sigma({}^{D_6^{(1)}} \hspace{-1mm} \Delta_{\rho_1},{}^{D_6^{(1)}} \hspace{-1mm} \Delta_{\rho_2})$ of the commuting normal operators ${}^{D_6^{(1)}} \hspace{-1mm} \Delta_{\rho_1}$, ${}^{D_6^{(1)}} \hspace{-1mm} \Delta_{\rho_2}$. In fact, ${}^{D_6^{(1)}} \hspace{-1mm} \Delta_{\rho_1}$ is a self-adjoint operator.

For the change of variables (\ref{eqn:x-D6^1}), (\ref{eqn:y-D6^1}), the Jacobian is given by
$J_{D_6^{(1)}} = -4 \pi^2 \omega_1^{-3}\omega_2^{-3} (1-\omega_1) (1-\omega_2) (\omega_1-\omega_2) (\omega_1+\omega_2+\omega_1\omega_2)^3$.
The Jacobian is complex and vanishes in $\mathbb{T}^2$ only on the images under $D_6^{(1)}$ of the boundaries of the fundamental domain corresponding to the curves $c_1$ and $c_2$ in $\mathfrak{D}$, and the singular point $(2/3,1/3)$. Note that $J_{D_6^{(1)}} \neq 0$ on the lines $\theta_2 = -\theta_1$ and $\theta_2 = 2\theta_1$, for $\theta_1 \in [1/2,2/3]$ -- which correspond to the curve $c_3$ in $\mathfrak{D}$ --  except at the endpoints $(1/2,1/2)$ and $(2/3,1/3)$, although $\mathrm{Re}(J_{D_6^{(1)}}) = -4 \pi^2 (\cos(2\pi (\theta_1 - 3\theta_2)) + \cos(2\pi (2\theta_1 - 3\theta_2)) - \cos(2\pi (3\theta_1 - 2\theta_2)) + \cos(2\pi (3\theta_1 - \theta_2)) - \cos(2\pi (2\theta_1 + \theta_2)) + \cos(2\pi (\theta_1 + 2\theta_2))) = 0$ on these lines.
Again, $J_{D_6^{(1)}}^2$ is invariant under the action of $D_6^{(1)}$ on $\mathbb{T}^2$, and can be written in terms of the $D_6^{(1)}$-invariant elements $x$, $y$ as $J_{D_6^{(1)}}^2 = -16 \pi^4 (y+2x+1)(16+24x-13x^2+2x^3+16y-4xy-x^2y+4y^2) = -\pi^4 (y+2x+1)(8y-x^2-4x+16+i\sqrt{x(8-x)^3})(8y-x^2-4x+16-i\sqrt{x(8-x)^3})$. Thus we can write $J$ in terms of $x$, $y$ as
$J_{D_6^{(1)}} = 4 \pi^2 i \sqrt{(y+2x+1)(16+24x-13x^2+2x^3+16y-4xy-x^2y+4y^2)} = \pi^2 i \sqrt{(y+2x+1)(8y-x^2-4x+16+\sqrt{Q})(8y-x^2-4x+16-\sqrt{Q})}$,
where $Q=x(x-8)^3$.
We have $Q \geq 0$ for $x \in [-1,0]$, whereas for $x \in [0,8]$, $Q \leq 0$ so that $\sqrt{Q}$ is purely imaginary.
The Jacobian $J_{D_6^{(1)}}$ is a cubic in $y$, with two of the three roots appearing as the equations of the boundaries $c_1$, $c_2$ of $\mathfrak{D}$. The third root intersects with $\mathfrak{D}$ only at the singular point $(-1,1)$.

\subsection{Spectral measure for $\mathcal{H}^{D_6^{(1)}}_{\rho}$ for $PSU(3)$} \label{sect:measure-rho-D6^1-H}

\begin{figure}[tb]
\begin{center}
  \includegraphics[width=80mm]{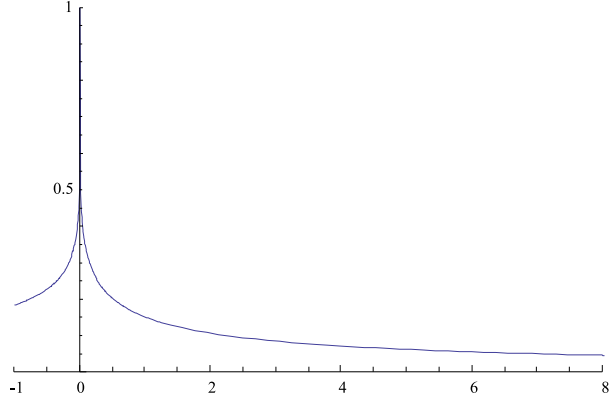}\\
 \caption{weight $\mathfrak{J}_1(x)$ for $\nu_{\rho_1}$ for $PSU(3)$.} \label{fig-Jx-D6_1}
\end{center}
\end{figure}

We first determine the spectral measure $\nu_{\rho_1}$ (over $\sigma_{\rho_1} = [-1,8]$) for ${}^{D_6^{(1)}} \hspace{-1mm} \Delta_{\rho_1}$.
For the change of variables $\mathbb{T}^2 \ni (e^{2\pi i \theta_1},e^{2\pi i \theta_2}) \mapsto (x,y_1)$, where $y_1=\mathrm{Re}(y)$, the Jacobian $J_{x,y_1}$ is given by
\begin{align*}
J_{x,y_1}(\theta_1,\theta_2) &= 4\pi^2 \left(\cos(2\pi(\theta_1+2\theta_2)) + \cos(2\pi(2\theta_1-3\theta_2)) + \cos(2\pi(3\theta_1-\theta_2)) \right. \\
& \qquad \left. - \cos(2\pi(2\theta_1+\theta_2)) - \cos(2\pi(3\theta_1-2\theta_2)) - \cos(2\pi(\theta_1-3\theta_2)) \right), \\
J_{x,y_1}(x,y_1) &= 4\pi^2\sqrt{(x^2+4x-16-8y_1)(y_1^2-x^3+x^2+x+2y_1+4xy_1)}.
\end{align*}
The map $\mathbb{T}^2 \ni (\omega_1,\omega_2) \mapsto (x,y_1)$ is a twelve-to-one map, and $J_{x,y_1}(x,y_1)$ is invariant under the dihedral group $D_{12}$ which contains $D_6^{(1)}$ as a normal subgroup.
Then by integrating $|D_{12}| \, |J_{x,y_1}(x,y_1)^{-1}|$ over $y_1$ we obtain the spectral measure $\nu_{\rho_1}$ for ${}^{D_6^{(1)}} \hspace{-1mm} \Delta_{\rho_1}$.
Thus the spectral measure $\nu_{\rho_1}$ (over $[-1,8]$) is $\mathrm{d}\nu_{\rho_1}(x) = \mathfrak{J}_1(x) \, \mathrm{d}x$, where $\mathfrak{J}_1(x)$ is given by
\begin{align*}
12 \int_{-2x-1-\sqrt{x+1}^3}^{-2x-1+\sqrt{x+1}^3} |J_{x,y_1}(x,y_1)|^{-1} \, \mathrm{d}y_1 & \textrm{ for } x \in [-1,0], \\
12 \int_{(x^2+4x-19)/8}^{-2x-1+\sqrt{x+1}^3} |J_{x,y_1}(x,y_1)|^{-1} \, \mathrm{d}y_1 & \textrm{ for } x \in (0,8].
\end{align*}
The weight $\mathfrak{J}_1(x)$ is an integral of the reciprocal of the square root of a cubic in $y$, and thus can be written in terms of the complete elliptic integral $K(m)$ of the first kind. Using \cite[equation 235.00]{byrd/friedman:1971}, we obtain
$$\mathfrak{J}_1(x) = \left\{ \begin{array}{cl}
\displaystyle \frac{6}{\pi^2\sqrt{8-20x-x^2+8\sqrt{x+1}^3}} \; K(v(x)) & \textrm{ for } x \in [-1,0], \\
\displaystyle \frac{3}{2\pi^2\sqrt[4]{x+1}^3} \; K(v(x)^{-1}) & \textrm{ for } x \in (0,8],
\end{array} \right. ,$$
where $v(x)=16\sqrt{x+1}^3/(8-20x-x^2+8\sqrt{x+1}^3)$.
The weight $\mathfrak{J}_1(x)$ is illustrated in Figure \ref{fig-Jx-D6_1}.

\begin{figure}[tb]
\begin{minipage}[t]{9cm}
\begin{center}
  \includegraphics[width=80mm]{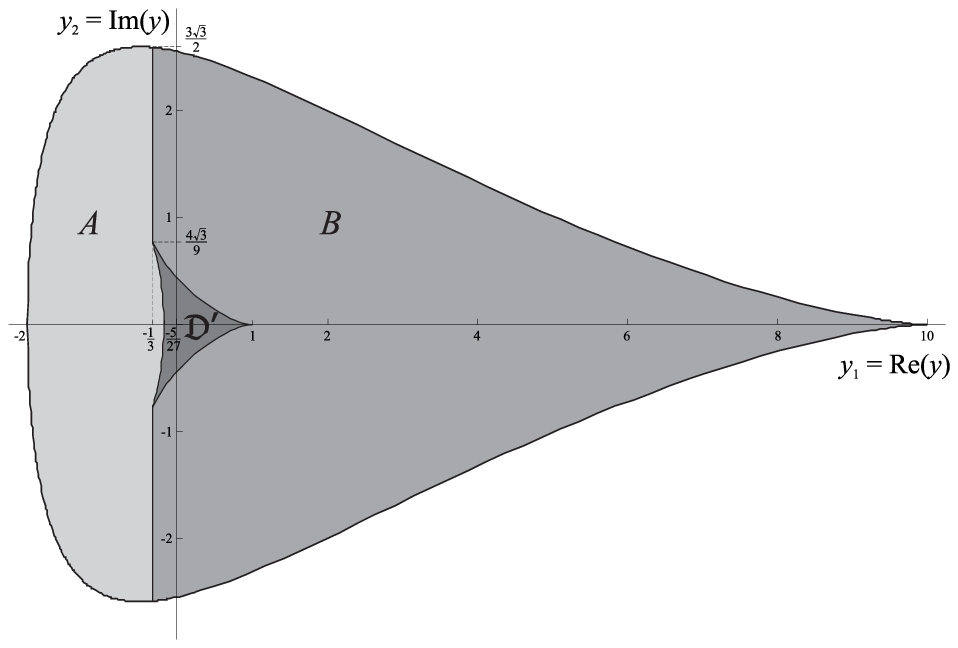}\\
 \caption{$\mathfrak{D}_{1,\rho_2}$ for $PSU(3)$.} \label{Fig-Dy-D6_1}
\end{center}
\end{minipage}
\hfill
\begin{minipage}[t]{7cm}
\begin{center}
  \includegraphics[width=70mm]{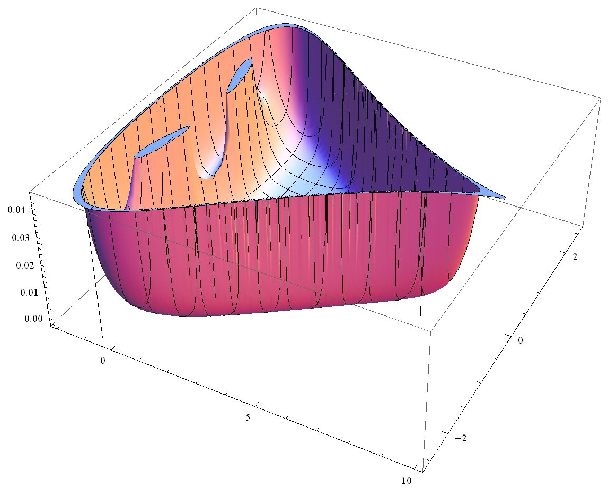}\\
 \caption{weight $\mathfrak{J}_2(y_1,y_2)$ for $\nu_{\rho_2}$ for $PSU(3)$.} \label{fig-Jy-D6_1}
\end{center}
\end{minipage}
\end{figure}

We now consider the spectral measure $\nu_{\rho_2}$ (over $\chi_{\rho_2}(\mathbb{T}^2)$) for ${}^{D_6^{(1)}} \hspace{-1mm} \Delta_{\rho_2}$. The spectrum $\sigma_{\rho_2} = \chi_{\rho_2}(\mathbb{T}^2)$ of $\rho_2$ is given by the region $\mathfrak{D}_{1,\rho_2} \subset \mathbb{C}$, illustrated in Figure \ref{Fig-Dy-D6_1}, whose boundary can be obtained from equations (\ref{eqn:y-x-D6^1-1}), (\ref{eqn:y-x-D6^1-2}) as $\mathrm{Im}(y)=\pm \sqrt{(a(y)-1)(5-a(y))^3}$, where $a(y)=\sqrt{5+2\mathrm{Re}(y)}$.
Let $y_1,y_2$ denote the real and imaginary parts of $y$ respectively.
The equation (\ref{eqn:surfaceD-D6^1}) of the surface $\mathfrak{D}$ is cubic in $x$, thus for any pair $(y_1,y_2)$ there are at most 3 values of $x$ for which $(x,y_1,y_2) \in \mathfrak{D}$.
Denote by $\mathfrak{D}'$ the subregion of $\mathfrak{D}_{1,\rho_2}$ illustrated in Figure \ref{Fig-Dy-D6_1}, whose boundaries are obtained when the discriminant $-256(1+3y_1)^3 + (11+90y_1+27y_1^2+27y_2^2)^2$ of the cubic polynomial (\ref{eqn:surfaceD-D6^1}) in $x$ is zero, i.e. by $y_2=\frac{\varepsilon}{3\sqrt{3}}\sqrt{-11-90y_1-27y_1^2+\varepsilon'16\sqrt{1+3y_1}^3}$ for $\varepsilon,\varepsilon'\in\{\pm1\}$.
Consider the restriction $\Upsilon = \Psi_{1,\rho_2}|_C$ to the fundamental domain $C$ of the map $\Psi_{1,\rho_2}: \mathbb{T}^2 \rightarrow \mathfrak{D}_{1,\rho_2}$. For $y \in \mathfrak{D}'$, its pre-image $\Upsilon^{-1}(y)$ in $C$ contains three distinct points, whilst for $y \in \mathfrak{D}_{1,\rho_2} \setminus \mathfrak{D}'$ its pre-image $\Upsilon^{-1}(y)$ is unique.

For the change of variables $\mathbb{T}^2 \ni (e^{2\pi i \theta_1},e^{2\pi i \theta_2}) \mapsto (y_1,y_2)$ the Jacobian $J_{y_1,y_2}$ is given by $J_{y_1,y_2}(\theta_1,\theta_2) = 4\pi^2 (2(\sin(2\pi(\theta_1-\theta_2))-\sin(2\pi\theta_1)+\sin(2\pi\theta_2)) + 4(\sin(4\pi(\theta_1-\theta_2))-\sin(4\pi\theta_1)+\sin(4\pi\theta_2)) + 3(\sin(6\pi(\theta_1-\theta_2))-\sin(6\pi\theta_1)+\sin(6\pi\theta_2)) + \sin(2\pi(\theta_1+2\theta_2)) + \sin(2\pi(2\theta_1-3\theta_2)) + \sin(2\pi(3\theta_1-2\theta_2)) - \sin(2\pi(2\theta_1+\theta_2)) - \sin(2\pi(\theta_1-3\theta_2)) + \sin(2\pi(3\theta_1-\theta_2)))$, which can be written in terms of the three $D_6^{(1)}$-invariant variables $x$, $y_1$, $y_2$ as
$$J_{y_1,y_2}(x,y_1,y_2) = 2\sqrt{2}\pi^2\sqrt{b(x,y_1,y_2)},$$
where $b(x,y_1,y_2) = 16 - 9 x + 34 x^2 - 2 y_1 + 178 x y_1 + 24 x^2 y_1 + 39 y_1^2 - 97 x y_1^2 - 10 x^2 y_1^2 + 4 y_1^3 + 24 x y_1^3 - 9 y_1^4 - 69 y_2^2 + 143 x y_2^2 - 42 x^2 y_2^2 - 124 y_1 y_2^2 + 24 x y_1 y_2^2 - 18 y_1^2 y_2^2 - 9 y_2^4$.
The solutions of (\ref{eqn:surfaceD-D6^1}) for $x$ are
\begin{equation} \label{eqn:xj-D6^1}
6x_j=2 - 2^{2/3}\epsilon_i P_x - 8 2^{1/3} \overline{\epsilon_i} (1 + 3 y_1)P_x^{-1}
\end{equation}
where $\epsilon_j=e^{2\pi i (j-1)}$, $j=1,2,3$, and
$$P_x= \left( -11 - 90 y_1 - 27 y_1^2 - 27 y_2^2 + \sqrt{-256 (1 + 3 y_1)^3 + (11 + 90 y_1 + 27 y_1^2 + 27 y_2^2)^2} \right)^{1/3}.$$
Then the spectral measure $\nu_{\rho_2}$ (over $\sigma_{\rho_2} = \mathfrak{D}_{1,\rho_2}$) is $\mathrm{d}\nu_{\rho_2}(y) = \mathfrak{J}_2(y_1,y_2) \, \mathrm{d}y_1 \, \mathrm{d}y_2$, where $\mathfrak{J}_2(y_1,y_2)$ is given by
\begin{equation} \label{eqn:J2(y)-D6^1}
\mathfrak{J}_2(y_1,y_2) = \left\{ \begin{array}{cl}
|J_{y_1,y_2}(x_1,y_1,y_2)|^{-1} & \textrm{ for } y \in A, \\
|J_{y_1,y_2}(x_2,y_1,y_2)|^{-1} & \textrm{ for } y \in B, \\
(\sum_{j=1}^3 |J_{y_1,y_2}(x_j,y_1,y_2)|)^{-1} & \textrm{ for } y \in \mathfrak{D}',
\end{array} \right.
\end{equation}
where region $A$ in Figure \ref{Fig-Dy-D6_1} is given by $y \in \mathfrak{D}_{1,\rho_2}$ such that $y_1 < -1/3$, and $27y_2^2 < -11-90y_1-27 y_1^2-16\sqrt{1+3y_1}^3$ for $-1/3 \leq y_1 < -5/27$; region $B$ by $27y_2^2 > -11-90y_1-27 y_1^2+16\sqrt{1+3y_1}^3$ for $-1/3 \leq y_1 \leq 1$, and $y_1 > 1$.
The weight $\mathfrak{J}_2(y_1,y_2)$ is illustrated in Figure \ref{fig-Jy-D6_1}.

\begin{Thm} \label{thm:measure-rho-D6^1-H}
The spectral measures $\nu_{\rho_1}$, $\nu_{\rho_2}$ (over $\chi_{\rho_1}(\mathbb{T}^2) = [-1,8]$, $\chi_{\rho_2}(\mathbb{T}^2) = \mathfrak{D}_{1,\rho_2}$ respectively) for the graphs $\mathcal{H}^{D_6^{(1)}}_{\rho_1}$, $\mathcal{H}^{D_6^{(1)}}_{\rho_2}$ respectively, for $PSU(3)$ are given by
\begin{align*}
\mathrm{d}\nu_{\rho_1}(x) &= \left\{ \begin{array}{cl}
\displaystyle \frac{6}{\pi^2\sqrt{8-20x-x^2+8\sqrt{x+1}^3}} \; K(v(x)) \, \mathrm{d}x & \textrm{ for } x \in [-1,0], \\
\displaystyle \frac{3}{2\pi^2\sqrt[4]{x+1}^3} \; K(v(x)^{-1}) \, \mathrm{d}x & \textrm{ for } x \in (0,8],
\end{array} \right., \\
\mathrm{d}\nu_{\rho_2}(y) &= \mathfrak{J}_2(y_1,y_2) \, \mathrm{d}y_1 \, \mathrm{d}y_2,
\end{align*}
where $\mathfrak{J}_2(y_1,y_2)$ is given by (\ref{eqn:J2(y)-D6^1}).
\end{Thm}

\begin{Rem} \label{Rem:D6^1-D6^2}
Similar remarks may be made regarding the relationship between the groups $D_6^{(1)}$ (the Weyl group for $PSU(3)$) and $D_6^{(2)}$ (the Weyl group for $SU(3)$) as were made in Remarks \ref{Rem:Z2^2-Z2^3}, \ref{Rem:D4^1-D4^2}, for the groups $\mathbb{Z}_2^{(2)}$, $\mathbb{Z}_2^{(3)}$ and $D_4^{(1)}$, $D_4^{(2)}$ respectively. The groups $D_6^{(1)}$ and $D_6^{(2)}$ are conjugate in $GL(2,\mathbb{R})$, with conjugating matrix $H$ given by $H=\left( \begin{array}{cc} 1 & -2 \\ 2 & -1 \end{array} \right)$ which has inverse $H^{-1}=\left( \begin{array}{cc} -1/3 & 2/3 \\ -2/3 & 1/3 \end{array} \right) = -\frac{1}{3} H$.
As in Remark \ref{Rem:Z2^2-Z2^3}, the origin of this relationship is the fact that the compact, connected Lie group $SU(3)$ (which has Weyl group $D_6^{(2)}$) is a triple cover of the compact, connected, but not simply connected Lie group $PSU(3)$. Thus not all irreducible representations of $SU(3)$ are irreducible representations of $PSU(3)$, but only those of triality zero, that is, those indexed by $(\lambda_1,\lambda_2) \in \mathbb{Z}^2$ such that $\lambda_1 - \lambda_2 \equiv 0$ mod(3).
\end{Rem}

\begin{Rem} \label{Rem:D6^1<D12}
There is also a relation between the Lie groups $PSU(3)$ and $G_2$, in that $D_6^{(1)}$ (the Weyl group for $PSU(3)$) is a normal subgroup of $D_{12}$ (the Weyl group for $G_2$). This relation is reflected in the spectral measures for these Lie groups. The spectral measures associated to $G_2$ were studied in \cite{evans/pugh:2012i}. The fundamental domain $F$ for $D_{12} = D_6^{(1)} \rtimes \mathbb{Z}_2$ in \cite[Figure 8]{evans/pugh:2012i} is obtained from the fundamental domain for $D_6^{(1)}$ in Figure \ref{fig:DomainC-D6^1} by imposing one extra symmetry which comes from the additional $\mathbb{Z}_2$ action. Suppose we denote by
$x',y'$ the variables for $D_{12}$ given by
\begin{align*}
x' & := 1 + \omega_1 + \omega_1^{-1} + \omega_2 + \omega_2^{-1} + \omega_1\omega_2^{-1} + \omega_1^{-1}\omega_2 \\
y' & := x' + 1 + \omega_1\omega_2 + \omega_1^{-1}\omega_2^{-1} + \omega_1^2\omega_2^{-1} + \omega_1^{-2}\omega_2 + \omega_1\omega_2^{-2} + \omega_1^{-1}\omega_2^2.
\end{align*}
Then we see that $x'=x-1$ whilst $y' = 2\mathrm{Re}(y)-x+2$, and there is a 2-to-1 homomorphism $\xi:\mathfrak{D}_{D_6^{(1)}} \rightarrow \mathfrak{D}_{D_{12}}$ from the domain $\mathfrak{D}_{D_6^{(1)}}$ of $D_6^{(1)}$ to the domain $\mathfrak{D}_{D_{12}}$ of $D_{12}$ such that $\xi(x)=x-1$ and $\xi(y)=2\mathrm{Re}(y)-x+2$, which maps the boundary of $\mathfrak{D}_{D_6^{(1)}}$ to part of the boundary of $\mathfrak{D}_{D_{12}}$. The rest of the boundary of $\mathfrak{D}_{D_{12}}$ is given by $\xi(v)$, where $v=(x,y) \in \mathfrak{D}_{D_6^{(1)}}$ such that $y \in \mathbb{R}$, i.e. $v = \Phi(t)$ for $t \in \mathbb{T}^2$ such that $t$ is fixed under the additional $\mathbb{Z}_2$ action in $D_{12}$.
Moreover, $J_{D_{12}} = 2\mathrm{Re}(J_{D_6^{(1)}}) = \xi(J_{D_6^{(1)}})$ (c.f. Remark \ref{Rem:Z2^2<D4^1}).
\end{Rem}

\subsection{Spectral measure for $\mathcal{G}^{D_6^{(1)}}_{\rho}$ for $PSU(3)$}

In this section we make the assumption that Conjecture \ref{conj:spectral_measure-orbitS} in Section \ref{sect:measures_over_T2-G} holds for the joint spectral measure of the pair of graphs $\mathcal{H}^{D_6^{(1)}}_{\rho_1}$, $\mathcal{H}^{D_6^{(1)}}_{\rho_2}$. Then we determine from this the spectral measures $\nu_{\rho_1}$, $\nu_{\rho_2}$ (over $\chi_{\rho_1}(\mathbb{T}^2) = [-1,8]$, $\chi_{\rho_2}(\mathbb{T}^2) = \mathfrak{D}_{1,\rho_2}$ respectively) for the graphs $\mathcal{H}^{D_6^{(1)}}_{\rho_1}$, $\mathcal{H}^{D_6^{(1)}}_{\rho_2}$ respectively, for $PSU(3)$. We have $S_{(1,0)}(\theta) = \omega_1-\omega_1^{-1}-\omega_2+\omega_2^{-1}+\omega_1^{-1}\omega_2-\omega_1\omega_2^{-1}$.

\begin{figure}[tb]
\begin{minipage}[t]{9cm}
\begin{center}
  \includegraphics[width=80mm]{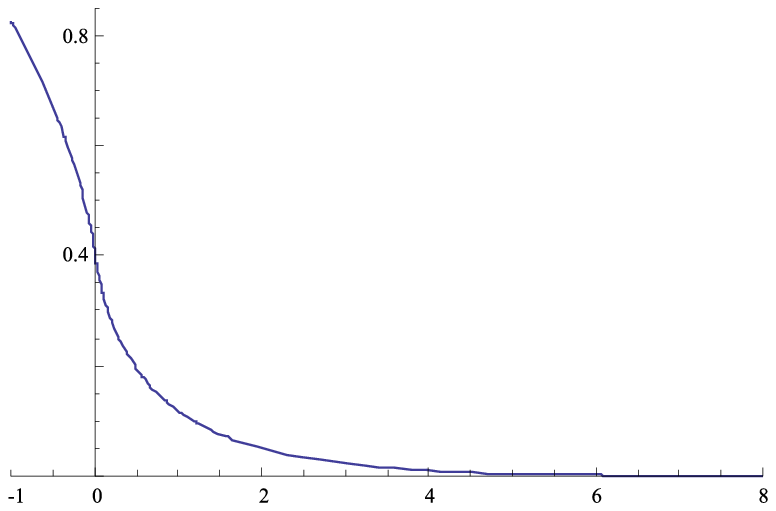}\\
 \caption{weight $\mathfrak{J}_3(x)$ for $\nu_{\rho_1}$ for $PSU(3)$.} \label{fig-Jx-D6_1-G}
\end{center}
\end{minipage}
\hfill
\begin{minipage}[t]{7cm}
\begin{center}
  \includegraphics[width=70mm]{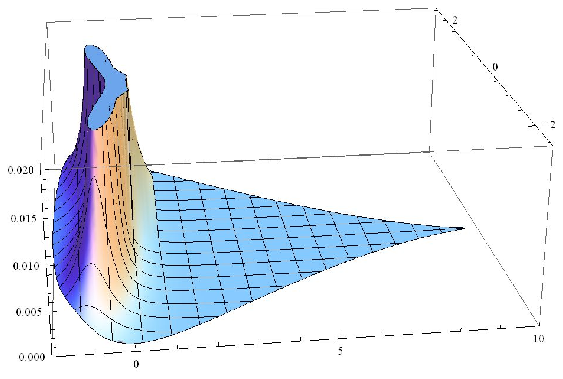}\\
 \caption{weight $\mathfrak{J}_4(y_1,y_2)$ for $\nu_{\rho_2}$ for $PSU(3)$.} \label{fig-Jy-D6_1-G}
\end{center}
\end{minipage}
\end{figure}

We consider first the spectral measure $\nu_{\rho_1}$. Now $S_{\varrho}(\theta) = S_{(1,0)}(\theta)$ may be written in terms of the variables $x$, $y_1$ as $S_{(1,0)}(x,y_1) = \sqrt{8y_1-x^2-4x+16}$. Then the spectral measure $\nu_{\rho_1}$ (over $[-1,8]$) is given by $\mathrm{d}\nu_{\rho_1}(x) = \mathfrak{J}_3(x) \, \mathrm{d}x$, where $\mathfrak{J}_3(x)$ is given by
\begin{align*}
4\sqrt{2} \int_{-2x-1-\sqrt{x+1}^3}^{-2x-1+\sqrt{x+1}^3} \frac{\sqrt{8y_1-x^2-4x+16}}{\sqrt{x^3-x^2-x-y_1^2-2y_1-4xy_1}} \, \mathrm{d}y_1 & \qquad \textrm{ for } x \in [-1,0], \\
4\sqrt{2} \int_{(x^2+4x-19)/8}^{-2x-1+\sqrt{x+1}^3} \frac{\sqrt{8y_1-x^2-4x+16}}{\sqrt{x^3-x^2-x-y_1^2-2y_1-4xy_1}} \, \mathrm{d}y_1 & \qquad \textrm{ for } x \in (0,8].
\end{align*}
The weight $\mathfrak{J}_3(x)$ can be written in terms of the complete elliptic integrals $K(m)$, $E(m)$ of the first, second kind respectively. Using \cite[equation 235.05, 235.06]{byrd/friedman:1971}, we obtain that $\mathfrak{J}_3(x)$ is given by
\begin{eqnarray*}
& \displaystyle \frac{\sqrt{8-20x-x^2+8\sqrt{x+1}^3}}{\pi^2} \; E(v(x)) & \textrm{ for } x \in [-1,0], \\
& \displaystyle \frac{1}{4\pi^2} \, \left( 16\sqrt[4]{x+1}^3 E(v(x)^{-1}) + \frac{\sqrt{8-20x-x^2-8\sqrt{x+1}^3}}{\sqrt[4]{x+1}^3} K(v(x)^{-1}) \right) & \textrm{ for } x \in (0,8],
\end{eqnarray*}
where $v(x)=16\sqrt{x+1}^3/(8-20x-x^2+8\sqrt{x+1}^3)$.
The weight $\mathfrak{J}_3(x)$ is illustrated in Figure \ref{fig-Jx-D6_1-G}.
It has been verified using Mathematica that the first 10 moments of this measure are correct.

We now turn to the spectral measure $\nu_{\rho_2}$ (over $\sigma_{\rho_2} = \mathfrak{D}_{1,\rho_2}$), which is given by $\mathrm{d}\nu_{\rho_2}(y) = \mathfrak{J}_4(y_1,y_2) \, \mathrm{d}y_1 \, \mathrm{d}y_2$, where $\mathfrak{J}_4(y_1,y_2)$ is given by
\begin{equation}
\mathfrak{J}_4(y_1,y_2) = \left\{ \begin{array}{cl}
|S_{(1,0)}(x_1,y_1)|^2/|J_{y_1,y_2}(x_1,y_1,y_2)| & \textrm{ for } y \in A, \\
|S_{(1,0)}(x_2,y_1)|^2/|J_{y_1,y_2}(x_2,y_1,y_2)| & \textrm{ for } y \in B, \\
\sum_{j=1}^3 |S_{(1,0)}(x_j,y_1)|^2/|J_{y_1,y_2}(x_j,y_1,y_2)| & \textrm{ for } y \in \mathfrak{D}',
\end{array} \right.
\end{equation}
where the regions $A$, $B$ and $\mathfrak{D}'$ are as in Section \ref{sect:measure-rho-D6^1-H}, and $x_j$ are as given by (\ref{eqn:xj-D6^1}).
The weight $\mathfrak{J}_4(y_1,y_2)$ is illustrated in Figure \ref{fig-Jy-D6_1-G}.

\bigskip \bigskip

\begin{footnotesize}
\noindent{\it Acknowledgement.}

The second author was supported by the Coleg Cymraeg Cenedlaethol.
\end{footnotesize}

\end{document}